\documentclass[pdflatex,sn-mathphys-ay]{sn-jnl}

\usepackage{graphicx}%
\usepackage{multirow}%
\usepackage{amsmath,amssymb,amsfonts}%
\usepackage{amsthm}%
\usepackage{mathrsfs}%
\usepackage[title]{appendix}%
\usepackage{xcolor}%
\usepackage{textcomp}%
\usepackage{manyfoot}%
\usepackage{booktabs}%
\usepackage{algorithm}%
\usepackage{algorithmicx}%
\usepackage{algpseudocode}%
\usepackage{listings}%
\usepackage{cleveref}
\usepackage{mathtools}

\theoremstyle{thmstyleone}%
\newtheorem{theorem}{Theorem}
\newtheorem{proposition}{Proposition}
\newtheorem{lemma}{Lemma}
\newtheorem{corollary}{Corollary}

\theoremstyle{thmstyletwo}%

\theoremstyle{thmstylethree}%

\crefname{theorem}{Theorem}{Theorems}
\Crefname{theorem}{Theorem}{Theorems}

\crefname{lemma}{Lemma}{Lemmas}
\Crefname{lemma}{Lemma}{Lemmas}

\crefname{proposition}{Proposition}{Propositions}
\Crefname{proposition}{Proposition}{Propositions}

\crefname{corollary}{Corollary}{Corollaries}
\Crefname{corollary}{Corollary}{Corollaries}

\crefname{definition}{Definition}{Definitions}
\Crefname{definition}{Definition}{Definitions}

\crefname{remark}{Remark}{Remarks}
\Crefname{remark}{Remark}{Remarks}

\crefname{example}{Example}{Examples}
\Crefname{example}{Example}{Examples}

\begin{document}

\title{An Improvement-Path Framework and an Exact Algorithm for Single-Machine Scheduling with Release Times}


\author[1]{\fnm{Xiaoyang} \sur{ Duan}}\email{xyduan@mail.sdu.edu.cn}

\author*[1]{\fnm{Peixin } \sur{Zhao}}\email{pxzhao@sdu.edu.cn}

\affil*[1]{\orgdiv{School of Management}, \orgname{Shandong University}, \orgaddress{\street{ 27  Shanda South Road}, \city{Jinan}, \postcode{250100}, \state{Shandong}, \country{China}}}


\abstract{
This paper studies the non-preemptive single-machine scheduling problem with heterogeneous release times and processing times, with the objective of minimizing total waiting time. The problem is known to be NP-hard. By modeling machine idle time as negative waiting time, we establish a unified waiting-time formulation that compresses the original four-dimensional description into a two-dimensional representation depending only on processing times and waiting times, thereby substantially simplifying the structural analysis of the problem.

Based on this representation, ideal directions are defined from the structure induced by an optimal solution, and it is shown that every non-optimal schedule admits an ideal improvement path. It is further proved that improvable ideal directions require no coordination and can therefore be treated as independent improvement units. The theoretical analysis also shows that the realizability of every improvable ideal direction can be determined in finitely many steps, thereby laying the foundation for an overall solution framework.

For each improvable ideal direction, by analyzing the propagation of the decrease flow and the increase flow under local adjustments, we prove that queue discontinuity is the unique structural obstacle to the realization of ideal improvement paths. A complete characterization of the resulting local optima is then established, and corresponding repair operators are designed.

On this basis, an iterative repair algorithm under the improvement-path framework is developed. It is proved to terminate in finitely many steps, return a globally optimal schedule, and admit an explicit upper bound on its time complexity. This work provides a new analytical perspective and solution framework for this class of NP-hard scheduling problems.}



\maketitle

\section{Introduction}\label{sec1}

Single-machine scheduling is one of the fundamental problems in operations research and also provides an important theoretical basis for more complex scheduling models, such as flow-shop scheduling and job-shop scheduling. In practical applications including computational resource allocation, healthcare service scheduling, and production planning, job completion times directly affect system responsiveness and service quality. Consequently, optimization problems related to completion times have long been a central topic in scheduling theory.

For the non-preemptive single-machine scheduling problem with heterogeneous release times and processing times, many classical additive completion-time-related objectives, such as total completion time, total waiting time, and the number of tardy jobs, have been shown to be NP-hard. The computational difficulty mainly arises from the coupling between release times and machine idle times. This coupling propagates globally throughout the schedule, making it difficult to assess the impact of local scheduling adjustments on the overall objective independently. Existing studies on this problem mainly fall into two categories. The first consists of exact algorithms based on frameworks such as Branch-and-Bound and Dynamic Programming, whose time complexity is typically exponential. The second consists of heuristic approaches based on dispatching rules or neighborhood search. Although computationally efficient, these methods generally lack rigorous optimality guarantees.

To overcome these limitations, the problem is studied from the structural perspective of the solution space, with emphasis on why a non-optimal schedule may fail to improve further and on the internal mechanism governing its evolution toward an optimal solution. By modeling machine idle time as negative waiting time, a unified waiting-time representation is established, which compresses the original four-dimensional information into a two-dimensional structure and thereby substantially reduces the complexity of the analysis. On this basis, ideal directions induced by an optimal solution are introduced, and improvable ideal directions are shown to admit no coordination, so that each can be treated as an independent improvement unit.

For each improvable ideal direction, the propagation of waiting-time variations along the queue is analyzed, revealing queue discontinuity as the key structural obstacle to its realization. Corresponding repair mechanisms are then constructed around queue discontinuity, together with a complete characterization of the resulting local optima. Building on this structural analysis, an iterative method based on structural repair is developed, under which the schedule is continuously improved through the repeated elimination of structural violations and eventually reaches a global optimum. Finite termination of the algorithm is further established.

\subsection{Related Work}

This subsection reviews research related to scheduling theory. The study of scheduling can be traced back to the pioneering work of \cite{Jackson1955}, who proposed one of the earliest earliest due date (EDD) rules under preemption and proved its optimality for minimizing maximum tardiness. Building on this line of research, \cite{Smith1956} introduced the shortest processing time (SPT) rule and proved that, when all jobs have identical release times, SPT is optimal for minimizing the total completion time $\sum C_j$. This principle was later extended to the weighted setting through the weighted shortest processing time (WSPT) rule. Subsequently, \cite{Lenstra1977} systematically characterized the computational complexity of a large class of scheduling problems, while \cite{Graham1979} proposed the well-known three-field notation, which established a unified framework for describing and classifying scheduling problems according to the machine environment, job characteristics, and objective function. A comprehensive survey of the development of scheduling theory can be found in \citep{Potts2009}.

In the area of single-machine scheduling, the development of exact algorithms has remained a central research direction. A variety of advanced algorithmic frameworks, including branch and bound (B\&B), Lagrangian relaxation, and dynamic programming (DP), have made substantial contributions to the field. \cite{Lawler1969} was the first to introduce dynamic programming into scheduling, successfully solving the weighted tardiness minimization problem with a common due date. Building on this foundation, \cite{Potts1982, Potts1992} proposed a pseudo-polynomial exact algorithm with time complexity $\mathcal{O}(nUB)$, which partially alleviated the computational bottleneck in solving larger instances. In the context of real-time scheduling, \cite{Lendl2023} studied rescheduling problems arising from unpredictable events such as the arrival of new orders and developed a pseudo-polynomial exact dynamic programming algorithm. Under the assumption that the resource consumption cost of a job is a non-increasing linear function of its release time, \cite{Ventura2002} proposed an exact dynamic programming algorithm for small- and medium-sized instances. To address uncertainty in practical environments, \cite{Malheiros2024} developed a robust dynamic programming framework for handling uncertain release times, due dates, and processing times. On the other hand, \cite{Tanaka2009, Tanaka2013} introduced the subsumed dynamic programming (SSDP) technique, which, under the strict no-idle-time assumption, extended the solvable instance size to 300 jobs.

The development of branch-and-bound methods has also been of major importance. \cite{Carlier1982} and \cite{Fisher1976} were among the first to incorporate Lagrangian relaxation into the branch-and-bound framework, achieving efficient pruning through the construction of tight lower bounds. More recently, \cite{Pessoa2021} proposed a hybrid strategy combining column generation and dual relaxation for batch scheduling, and for the first time solved exactly a batch scheduling problem with total weighted waiting time (TWT) as the objective. Overall, most existing exact algorithms for single-machine scheduling problems exhibit exponential-time complexity.

\subsection{Main Contributions}

The main contributions of this paper can be summarized as follows:

\begin{enumerate}
    \item A unified waiting-time representation is proposed. By modeling machine idle time as negative waiting time, a unified waiting-time formulation is established, which compresses the original four-dimensional state space into a two-dimensional structural representation.

    \item An analytical framework based on ideal improvement paths is developed. Ideal directions are defined from structures induced by an optimal solution. It is shown that improvable ideal directions require no coordination and that their realizability can be determined in finitely many steps. Furthermore, queue discontinuity is identified as the unique structural obstacle to the realization of improvable ideal directions, and a complete characterization of the resulting local optima is established.

    \item An iterative algorithm based on structural repair is developed. The algorithm repeatedly identifies and repairs structural obstacles, thereby continuously driving the current solution along ideal improvement paths, and is proved to reach a global optimum in finitely many iterations.
\end{enumerate}

\subsection{Organization of the Paper}

The remainder of this paper is organized as follows. Section 2 presents the problem definition and several basic properties, and rigorously reduces the original four-dimensional representation to a two-dimensional form through mathematical derivation. Section 3 analyzes the structural obstacles that prevent the realization of ideal improvement directions and develops a complete classification of local optimality types. Section 4 provides a detailed description of the proposed algorithm and proves its optimality and time complexity. Section 5 concludes the paper and discusses directions for future research.

\section{Problem Formulation and Fundamental Properties}

This paper studies the non-preemptive single-machine scheduling problem with heterogeneous release times and processing times, and takes total waiting time as the representative objective for structural analysis. Consider a set of $n$ independent jobs to be processed on a single machine. Each job $i$ is associated with a processing time $p_i$ and a release time $r_i$. For a given schedule, let $S_i$ denote the starting time of job $i$, and let $C_i=S_i+p_i$ denote its completion time. In this paper, we focus on total waiting time, where the waiting time of job $i$ is defined as $w_i=S_i-r_i$.

Since $C_i=S_i+p_i$, we have
\begin{align*}
    \sum C_i
    &= \sum (S_i+p_i) \\
    &= \sum (r_i+w_i+p_i) \\
    &= \sum (r_i+p_i)+\sum w_i,
\end{align*}
where $\sum (r_i+p_i)$ is a constant. Therefore, under release-time and non-preemptive scheduling constraints, minimizing the total waiting time $\sum_i w_i$ is equivalent to minimizing the total completion time $\sum_i C_i$. Moreover, Lenstra et al.~\citep{Lenstra1977} proved that minimizing total completion time in the non-preemptive single-machine scheduling problem with release times is NP-hard. Hence, minimizing total waiting time under the same constraints is also NP-hard.

To provide a unified characterization of machine idle time, we introduce an extended representation of waiting time in which $w_i<0$ is allowed. This extension is used solely for analysis and does not change the feasible solution set induced by the schedule. In any feasible schedule, we always have $S_i \ge r_i$, and hence $w_i \ge 0$, so that $\sum_i w_i=\sum_i \max(0,w_i)$. Under the extended representation, the objective function can equivalently be written as $\min \sum_i \max(0,w_i)$. In this setting, negative waiting time $w_i<0$ indicates that the machine is idle before job $i$ starts processing, and the idle length is given by $-\min(0,w_i)$.

Consider an interval during which the machine processes jobs continuously without idle time, where the first processed job is $i$, followed by jobs $i+1,i+2,\dots,j$ in sequence. Under the non-preemptive constraint, the start time of job $j$ can be expressed as
\[
S_j = r_i^o + \sum_{k=i}^{j-1} p_k, \qquad (\forall j>i),
\]
where $r_i^o$ denotes the starting time of this interval, that is, the earliest time at which the machine begins continuously processing this block. Therefore, the waiting time of job $j$ is
\[
w_j = r_i^o + \sum_{k=i}^{j-1} p_k - r_j.
\]

\begin{figure}[H]
\centering
\includegraphics[width=0.8\textwidth]{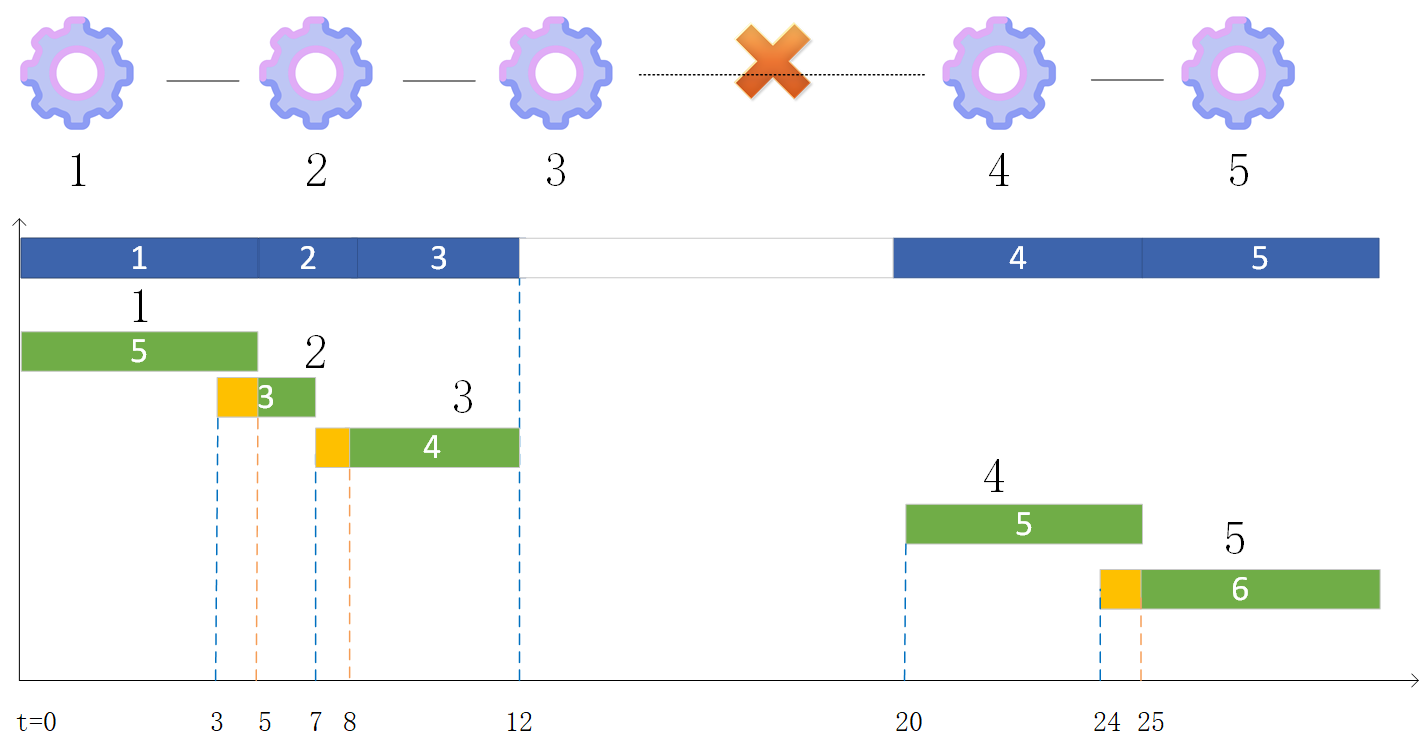}
\caption{This figure illustrates the calculation of queue waiting times and the update process of queue leaders. Consider a queue consisting of five jobs processed in the order $\{1,2,3,4,5\}$. Their processing times are $p_1=5$, $p_2=3$, $p_3=4$, $p_4=5$, and $p_5=6$ (green bars). Their release times are $r_1=0$, $r_2=3$, $r_3=7$, $r_4=20$, and $r_5=25$ (blue dashed lines).\\
Job 1 is the initial queue leader, and the waiting times are computed from it. At this point, $r_1^o=0$, and the waiting time of job 1 is $w_1=0$. By the definition, the waiting time of job 2 is $w_2=r_1^o+p_1-r_2=0+5-3=2$, so $w_2>0$, which means that job 2 belongs to the contiguous queue led by job 1. For job 3, we have $w_3=r_1^o+p_1+p_2-r_3=0+5+3-7=1$, so $w_3>0$, and the queue remains contiguous. The waiting time of job 4 is
$w_4=r_1^o+p_1+p_2+p_3-r_4=0+5+3+4-20=-8$. Since $w_4<0$, the queue leader is updated and $r_4^o=r_4=20$. Job 4 then becomes the new queue leader and starts a new contiguous queue. Finally, the waiting time of job 5 is
$w_5=r_4^o+p_4-r_5=20+5-24=1$. Since $w_5>0$, the queue remains contiguous.\\
Therefore, under this processing order, the waiting times are $w_1=0$, $w_2=2$, $w_3=1$, $w_4=-8$, and $w_5=1$ (yellow bars), forming two contiguous queues: $\{1,2,3\}$ and $\{4,5\}$. The corresponding lengths of continuous machine processing are 12 and 11, respectively (blue bars). The release time of job 4 forces the machine to remain idle for $-\min(0,-8)=8$ time units (white bar), while the actual waiting time of job 4 is $\max(0,-8)=0$.}
\label{waiting time calculation}
\end{figure}

We further define a \emph{contiguous queue} as a consecutive subsequence of jobs in processing order such that its first job, say $i$, satisfies $w_i \le 0$, and is referred to as the \emph{queue leader}, while all remaining jobs satisfy $w_k>0$ for $k>i$. For the entire job sequence, all positions satisfying $w_k \le 0$ are called \emph{breakpoints}. In particular, each breakpoint corresponds to the starting position of a contiguous queue, namely, to its queue leader, and the queue is said to exhibit a \emph{queue discontinuity} at that breakpoint. Figure~\ref{waiting time calculation} illustrates the waiting-time calculation process in detail through an example.

To accurately describe the processing order of jobs during the iterations, let $J=\{1,2,\ldots,n\}$ be the job set. We define
\[
\sigma_t:\{1,2,\ldots,n\}\to J,\qquad t\in\mathbb{Z}_{\ge 0},
\]
to denote the job permutation obtained at iteration $t$, where $\sigma_t(k)$ represents the job assigned to the $k$th processing position at iteration $t$. For any $t$, $\sigma_t$ is a permutation of the job set $J$.

The initial sequence $\sigma_0$ is sorted in nondecreasing order of release times:
\[
r_{\sigma_0(1)} \le r_{\sigma_0(2)} \le \cdots \le r_{\sigma_0(n)}.
\]

The sequence is updated iteratively by algorithm $\mathcal{A}$ according to
\[
\sigma_{t+1}=\mathcal{A}(\sigma_t), \qquad t=0,1,2,\ldots
\]

\subsection{Propagation of Waiting-Time Variations in Queues}
\begin{figure}[h]
    \centering
    \includegraphics[width=0.8\textwidth]{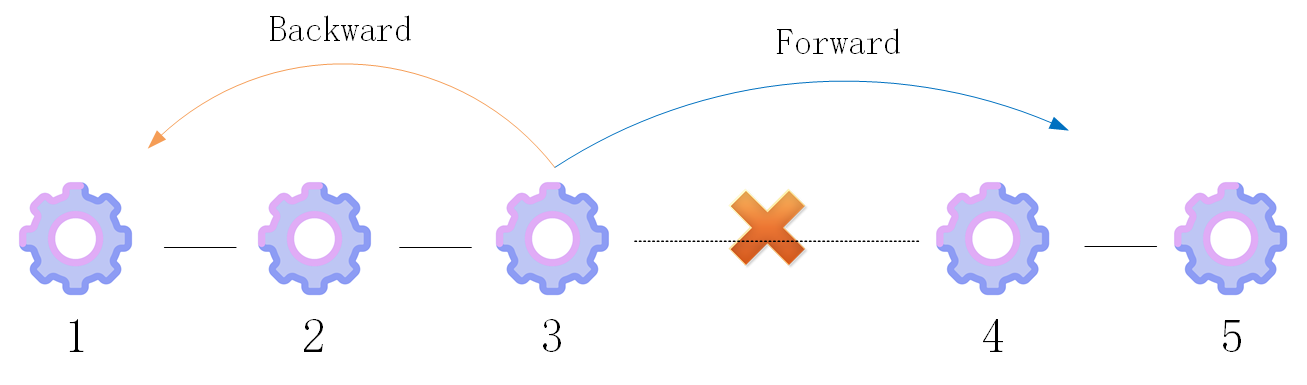}
    \caption{Illustration of position adjustment for job 3. If changing the processing position of job 3 can reduce the total waiting time, then the adjustment can be carried out in two basic ways: a \textit{forward move} (blue arrow) and a \textit{backward move} (orange arrow).}
    \label{forward and backward}
\end{figure}

As illustrated in Fig.~\ref{forward and backward}, consider a sequence 
$\sigma_t=\{\sigma_t(1),\dots,\sigma_t(n)\}$ and a job $\sigma_t(j)$ within it (job 3 is taken as an example in the figure). The positional adjustment of 
$\sigma_t(j)$ is defined as follows.  
(1) A \emph{forward move} relocates $\sigma_t(j)$ from its current position $j$ to a position $k>j$, while preserving the relative order of all remaining jobs (as illustrated by the rightward movement of job 3 in the figure).  
(2) A \emph{backward move} relocates $\sigma_t(j)$ from its current position $j$ to a position $k<j$, again without changing the relative order of the remaining jobs (as illustrated by the leftward movement of job 3).  Both types of moves can be equivalently represented as a \emph{remove-and-insert} operation: $\sigma_t(j)$ is first removed from its original position and then inserted into the target position $k$.

To facilitate the analysis, we further decompose each positional adjustment into the two elementary operations of removal and insertion, and examine their respective effects on the waiting times of jobs in the queue. Specifically, we study:  
(1) how the removal of $\sigma_t(j)$ affects the waiting times of the jobs that follow it; and  
(2) how the insertion of $\sigma_t(j)$ at a new position $k$ propagates waiting-time changes to the jobs located after position $k$.

These two elementary processes are characterized in \cref{Lemma1} and \cref{Lemma2}, respectively, where unified recursive formulas are established for the corresponding waiting-time propagation. Building on these results, \cref{Theorem1} and \cref{Theorem2} further derive the structural properties of forward and backward moves induced by the remove-and-insert operation, together with their effects on the total waiting time.

\begin{lemma}[Propagation of Waiting-Time Decrease]
\label{Lemma1}
Let the current job sequence be $\sigma_t$. Suppose that the waiting time of job $\sigma_t(j)$, where $j=1,2,\dots,n$, decreases by an initial amount $\Delta w_{\sigma_t(j)}^-=\Delta$ with $\Delta\geq 0$. Then the waiting-time changes of the subsequent jobs $\sigma_t(k)$, where $k=j,j+1,\dots,n$, satisfy the following recursion:
\begin{align*}
    \Delta w_{\sigma_t(k+1)}^- =
    \begin{cases}
        0, & \text{if } w_{\sigma_t(k)} \leq 0,\\
        \min\bigl(\Delta w_{\sigma_t(k)}^-,\, w_{\sigma_t(k)}\bigr), & \text{if } w_{\sigma_t(k)} > 0.
    \end{cases}
\end{align*}
\end{lemma}

\textbf{Remark.} The quantity $\Delta w_{\sigma_t(n+1)}^-$ has a concrete interpretation, which will be explained in \cref{Proposition1}.

The proof of \cref{Lemma1} is based on the recursive structure of the waiting-time formula. Specifically, we first analyze how a decrease in the waiting time of a given job affects the waiting time of its immediate successor, and then extend this local propagation recursively to all subsequent jobs. In this way, the recursive propagation rule for waiting-time decrease is established. The detailed proof is provided in the Appendix~\ref{app:Lemma1}.

\begin{lemma}[Propagation of Waiting-Time Increase]
\label{Lemma2}
Let the current job sequence be $\sigma_t$. Suppose that the waiting time of job $\sigma_t(j)$, where $j=1,2,\dots,n$, increases by an initial amount $\Delta w_{\sigma_t(j)}^+=\Delta$ with $\Delta\geq 0$. Then the waiting-time changes of the subsequent jobs $\sigma_t(k)$, where $k=j,j+1,\dots,n$, satisfy the following recursion:
\begin{equation*}
    \Delta w_{\sigma_t(k+1)}^+ =
    \begin{cases}
        \Delta w_{\sigma_t(k)}^+, & \text{if } w_{\sigma_t(k)} > 0,\\
        \max\bigl(0,\, \Delta w_{\sigma_t(k)}^+ + w_{\sigma_t(k)}\bigr), & \text{if } w_{\sigma_t(k)} \leq 0.
    \end{cases}
\end{equation*}
\end{lemma}

\textbf{Remark.} Likewise, the interpretation of $\Delta w_{\sigma_t(n+1)}^+$ will be given in \cref{Proposition1}.

The proof of \cref{Lemma2} follows the same line of argument as that of \cref{Lemma1}. The details are deferred to the Appendix	~\ref{app:Lemma2}.

\begin{figure}[h]
    \centering
    \includegraphics[width=0.8\linewidth]{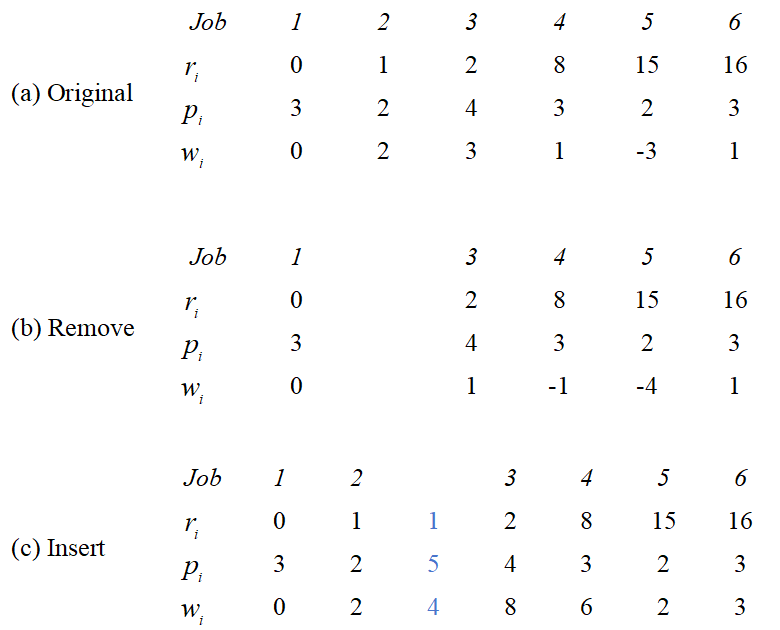}
    \caption{(a) The original queue. (b) The changes in waiting times after removing job 2 from the original queue in (a). (c) The changes in waiting times after inserting a new job, with release time 1 and processing time 5, between job 2 and job 3 in the original queue in (a).}
    \label{Removal and insert}
\end{figure}

\cref{Lemma1} and \cref{Lemma2} characterize the recursive propagation of waiting-time changes induced by removal and insertion, respectively. Specifically, \cref{Lemma1} describes how the removal of a job propagates a decrease in waiting times to the subsequent jobs, whereas \cref{Lemma2} describes how the insertion of a job propagates an increase in waiting times to the subsequent jobs. Illustrative examples are provided in Fig.~\ref{Removal and insert}.

A further issue is that waiting-time variation and objective-value variation are not aligned at the job level. Since the objective is
\[
\sum_{i=1}^{n}\max(0,w_i),
\]
the change in the waiting time of an individual job does not necessarily induce an equal change in the total waiting time. In particular, although the local variation of job $j$ is represented by $\Delta w_j^+$ or $-\Delta w_j^-$, the term that matters for the objective is the propagated variation on the subsequent job, namely $\Delta w_{j+1}^+$ or $-\Delta w_{j+1}^-$. Hence, there is an index shift between waiting-time propagation and objective variation.

This shift arises because a variation in $w_j$ no longer affects the objective once $w_j$ becomes nonpositive, whereas its propagated effect may still continue along the queue through the recursive relations. Equivalently, only the positive part entering $\max(0,w_j)$ can produce an effective change in the total waiting time. This observation is formalized in \cref{Proposition2}; the full proof is deferred to Appendix~\ref{app:Proposition2}.

\begin{proposition}
\label{Proposition2}
According to the recursive relations in \cref{Lemma1} and \cref{Lemma2}, let the current job sequence be $\sigma_t$, and suppose that the waiting time of job $\sigma_t(j)$ either increases by $\Delta w_{\sigma_t(j)}^+ = \Delta$ or decreases by $\Delta w_{\sigma_t(j)}^- = \Delta$, where $\Delta \geq 0$. Then the updated waiting time of $\sigma_t(j)$ is given by
\[
w_{\sigma_t(j)}^2
=
w_{\sigma_t(j)} + \Delta w_{\sigma_t(j)}^+
\qquad \text{or} \qquad
w_{\sigma_t(j)}^2
=
w_{\sigma_t(j)} - \Delta w_{\sigma_t(j)}^-.
\]
Accordingly, the actual variation in total waiting time is
\[
\Delta w = \Delta w_{\sigma_t(j+1)}^+
\qquad \text{or} \qquad
\Delta w = -\Delta w_{\sigma_t(j+1)}^-,
\]
where $\Delta w_{\sigma_t(j+1)}^+$ and $\Delta w_{\sigma_t(j+1)}^-$ are determined by the corresponding recursive relations.
\end{proposition}

The detailed proof is given in the Appendix~\ref{app:Proposition2}.

We refer to the phenomenon that a decrease in the waiting time of job $j$ causes decreases in the waiting times of subsequent jobs as the \emph{decrease flow}, denoted by
\[
\pi_{\sigma_t(j)}^-\bigl(\Delta w_{\sigma_t(j)}\bigr)
=
-\sum_{k=j+1}^{n+1}\Delta w_{\sigma_t(k)}^-.
\]

Similarly, we refer to the phenomenon that an increase in the waiting time of job $j$ causes increases in the waiting times of subsequent jobs as the \emph{increase flow}, denoted by
\[
\pi_{\sigma_t(j)}^+\bigl(\Delta w_{\sigma_t(j)}\bigr)
=
\sum_{k=j+1}^{n+1}\Delta w_{\sigma_t(k)}^+.
\]

Combining the decrease flow and the increase flow, we define the unified function
\begin{align*}
\pi_{\sigma_t(j)}\bigl(\Delta w_{\sigma_t(j)}\bigr)
=
\begin{cases}
\pi_{\sigma_t(j)}^+\bigl(\Delta w_{\sigma_t(j)}\bigr), & \text{if } \Delta w_{\sigma_t(j)} > 0,\\
\pi_{\sigma_t(j)}^-\bigl(-\Delta w_{\sigma_t(j)}\bigr), & \text{if } \Delta w_{\sigma_t(j)} \leq 0.
\end{cases}
\end{align*}

Building on the waiting-time propagation mechanisms established in \cref{Lemma1} and \cref{Lemma2}, we next study the effect of a positional adjustment of a single job on the entire queue. Specifically, \cref{Theorem1} characterizes the new waiting time of a job postponed by a forward move, together with the induced change in the total waiting time, whereas \cref{Theorem2} provides the corresponding characterization for a job advanced by a backward move. These results lift the local propagation laws of removal and insertion to a queue-level characterization of job positional adjustments, and provide the basis for the subsequent candidate-solution screening and local structural analysis.

\begin{theorem}
\label{Theorem1}
Given the current job sequence $\sigma_t$, suppose that an arbitrary job $\sigma_t(i)$ is moved forward to the position immediately after job $\sigma_t(k)$, where $k \in \{i+1,\dots,n\}$. Then the variation in total waiting time, denoted by $\Delta w$, can be computed as
\begin{align*}
\Delta w = I_1 + I_2 + \pi_{\sigma_t(k+1)}(I_2),
\end{align*}
where
\begin{align*}
I_1 = \sum_{j=i+1}^{k} \left( p_{\sigma_t(j)} - \min\left(0, w_{\sigma_t(j)}\right) \right) - \sum_{j=i+2}^{k+1} \Delta w_{\sigma_t(j)}^-,
\end{align*}
and
\begin{align*}
I_2 =
\begin{cases}
\begin{aligned}
\max\Bigl(&\min\left(0, w_{\sigma_t(i)}\right), \\
          &-\sum_{j=i+1}^{k} p_{\sigma_t(j)}, \\
          &\max_{j \in [i+1,k]} \left( p_{\sigma_t(i)} - w_{\sigma_t(j)} \right)\Bigr),
\end{aligned}
& \text{if } w_{\sigma_t(j)} > 0,\ \forall j \in [i+1,k], \\[1ex]
p_{\sigma_t(i)},
& \text{otherwise}.
\end{cases}
\end{align*}

Moreover, after relocating job $\sigma_t(i)$, its new waiting time $w_{\sigma_t(i)}^2$ can be expressed as
\begin{align*}
w_{\sigma_t(i)}^2
=
\sum_{j=i}^{k} p_{\sigma_t(j)}
-
\sum_{j=i+1}^{k} \min\left(0, w_{\sigma_t(j)}\right)
-
\Delta w_{\sigma_t(k+1)}^-
+
\max\left(0, w_{\sigma_t(i)}\right).
\end{align*}

Here, $\Delta w_{\sigma_t(j)}^-$ is determined by the recursion in \cref{Lemma1}, with initial value
\begin{align*}
\Delta w_{\sigma_t(i+1)}^- = p_{\sigma_t(i)} - \min\left(0, w_{\sigma_t(i)}\right).
\end{align*}
\end{theorem}

\begin{proof}[Proof sketch]
We first consider a unit forward move of a job in a fixed sequence. Using the recursive waiting-time relations, we compare the waiting times before and after the move and derive the resulting change in the total waiting time. We then extend the analysis to a two-position forward move and show that enlarging the moving distance only adds the processing-time and waiting-time information of the newly crossed jobs, while leaving the structural form of the objective-change expression unchanged. Since no new structural term arises, the same argument extends to an arbitrary forward move, yielding the general formula. The complete derivation is provided in Appendix~\ref{app:Theorem1}.
\end{proof}

\begin{theorem}
\label{Theorem2}
Given the current job sequence $\sigma_t$, suppose that an arbitrary job $\sigma_t(i)$ is moved backward to the position immediately before job $\sigma_t(k)$, where $k \in \{1,\dots,i-1\}$. Then the variation in total waiting time, denoted by $\Delta w$, can be computed as
\[
\Delta w = I_3 + I_4 + \pi_{\sigma_t(k)}(I_4).
\]

Here,
\[
I_3 = \sum_{j=k}^{i-1} \left( \min(0, w_{\sigma_t(j)}) - p_{\sigma_t(j)} \right) + \sum_{j=k}^{i-1} \Delta w_{\sigma_t(j)}^+,
\]
where $\Delta w_{\sigma_t(j)}^+$ is determined by the recursion in \cref{Lemma2}, with initial value
\begin{align*}
\Delta w_{\sigma_t(k)}^+
=
\max \biggl(
    p_{\sigma_t(i)},
    \sum_{j=k}^{i} p_{\sigma_t(j)}
    -
    \sum_{j=k}^{i-1} \min(0, w_{\sigma_t(j)})
    -
    w_{\sigma_t(i)}
\biggr).
\end{align*}

In addition,
\begin{align*}
I_4
=
\max \biggl(
    -p_{\sigma_t(i)},
    \sum_{j=k}^{i-1} \min(0, w_{\sigma_t(j)}),
    \sum_{j=k}^{i-1} p_{\sigma_t(j)} - \max(0, w_{\sigma_t(i)})
\biggr).
\end{align*}

Moreover, after relocating job $\sigma_t(i)$, its new waiting time $w_{\sigma_t(i)}^2$ is given by
\[
w_{\sigma_t(i)}^2
=
\sum_{j=k}^{i-1} \left( \min(0, w_{\sigma_t(j)}) - p_{\sigma_t(j)} \right) + w_{\sigma_t(i)}.
\]
\end{theorem}

The proof of \cref{Theorem2} is analogous to that of \cref{Theorem1} and is omitted for brevity. The complete derivation is deferred to Appendix~\ref{app:Theorem2}.

Building on the above analysis, the original four-dimensional description of the problem, involving processing times, release times, machine idle times, and waiting times, can be compressed into a two-dimensional representation involving only processing times and waiting times. This reformulation further yields a local criterion for the ordering of adjacent jobs. Specifically, the relationship between the processing times and waiting times of a job and its immediate successor determines whether the local order of the pair remains consistent with the nondecreasing order of release times. Otherwise, the local service order has shifted from the initial first-come-first-served (FCFS) rule to a last-come-first-served (LCFS) pattern. \cref{Proposition1} formalizes the relationship between waiting times and processing times under these two local ordering regimes, thereby providing a criterion for identifying whether a local exchange has occurred between two adjacent jobs. The full derivation is deferred to Appendix~\ref{app:Proposition1}.

\begin{proposition}
\label{Proposition1}
For the job sequence $\sigma_t$, if the release times satisfy$r_{\sigma_t(i)} \leq r_{\sigma_t(i+1)},$ then $\max(0,w_{\sigma_t(i)}) + p_{\sigma_t(i)} \geq \max(0,w_{\sigma_t(i+1)}).$ 

Otherwise, if $r_{\sigma_t(i)} > r_{\sigma_t(i+1)},$ then $\max(0,w_{\sigma_t(i)}) + p_{\sigma_t(i)} < \max(0,w_{\sigma_t(i+1)})$holds.
\end{proposition}

\section{An Analysis Framework Based on Feasible Improvement Paths}

For any non-optimal schedule $\sigma_t$, there exists at least one ideal forward or backward direction induced by an optimal solution. This section studies such directions from two perspectives: their mutual relationship, and the structural obstacles to their realization. An improvable ideal direction, however, need not be directly executable. Indeed, the waiting-time variation induced by a job move propagates along the queue and may create obstructions during the propagation process, thereby preventing the ideal improvement from being realized directly. It is therefore necessary to determine both whether different improvable ideal directions interact in a coordinated manner and what structural factors obstruct their realization.

We first show in \cref{independent} that improvable ideal directions admit no coordination relationship. Hence, each such direction can be treated as an independent improvement unit. We then analyze the structural obstacles arising in forward and backward moves, and prove in \cref{zhuai} that queue discontinuity is the unique structural obstacle to the realization of an improvable ideal direction.

Once a queue discontinuity arises, the original ideal move path is blocked. Restoring its feasibility therefore requires repairing the discontinuity so that the blocked ideal improvement can again be converted into an executable actual improvement. To this end, \cref{Lemma3} characterizes the structural conditions under which a job can repair a discontinuity, and shows in particular that the repairing job must satisfy the corresponding release-time constraints, thereby yielding a feasible criterion for discontinuity repair. Building on this result, \cref{finite ideal} proves that the realizability of every improvable ideal direction can be decided in finitely many steps, and \cref{Proposition3} further provides a complete characterization of the induced local-optimum types.

\subsection{Ideal Directions and Ideal Improvement Values}

Let $\sigma_t$ be the current schedule and $\sigma^*$ an optimal schedule, with $\sigma_t\neq \sigma^*$. Then there exists at least one job whose position in $\sigma_t$ differs from that in $\sigma^*$. Such a positional mismatch induces an ideal direction: if the job appears before its position in $\sigma^*$, moving it forward to that position defines an \emph{ideal forward direction}; if it appears after its position in $\sigma^*$, moving it backward to that position defines an \emph{ideal backward direction}.

Hence, every non-optimal schedule admits at least one ideal forward or ideal backward direction induced by some optimal schedule.

Let $M^*$ be such a direction, and let $F(\sigma)$ denote the total waiting time of schedule $\sigma$. Define the \emph{ideal improvement value} of $M^*$ by
\[
Z(M^*) = F(\sigma_t^{I}) - F(\sigma_t),
\]
where $\sigma_t^{I}$ is the schedule obtained by applying the ideal adjustment prescribed by $M^*$ to $\sigma_t$. If $Z(M^*)<0$, then $M^*$ is theoretically capable of improving the current schedule.

The existence of an ideal direction, however, does not imply that $Z(M^*)<0$. Thus, not every ideal direction has improvement potential. By \cref{Theorem1} and \cref{Theorem2}, any positional adjustment simultaneously induces an increase flow and a decrease flow, and the two effects act against each other in the objective. Consequently, only ideal directions satisfying certain structural conditions can reduce the total waiting time.

We next derive necessary conditions for a forward or backward move to possess such improvement potential.

\begin{lemma}[Criterion for the Forward Candidate Set]
\label{forward condition}
For a schedule $\sigma_t$, consider moving job $\sigma_t(i)$ forward from position $i$ to position $k$, where $k>i$. If this move satisfies the following necessary conditions:
\begin{itemize}
    \item $\sum_{j=i+1}^{k} \Bigl( p_{\sigma_t(j)} - \bigl(p_{\sigma_t(i)} - \min(w_{\sigma_t(i)}, 0)\bigr) \Bigr) \leq 0$;
    \item queue discontinuity occurs,
\end{itemize}
then the move has the potential to reduce the total waiting time. In this case, position $k$ is called a \emph{forward candidate solution} of job $\sigma_t(i)$, denoted by $\theta_{i,k}^f$.
\end{lemma}
\begin{proof}[Proof sketch]
The argument is based on a case-by-case analysis of the different combinations of $I_1$ and $I_2$ in \cref{Theorem1}. By further exploiting the recursive decrease-flow relation in \cref{Lemma1}, we derive a lower bound on the resulting change in the total waiting time. This lower bound immediately yields a necessary condition under which a forward move can be improving. The complete proof is given in Appendix~\ref{app:forward condition}.
\end{proof}

Based on Condition 1, we define the forward candidate set of job $\sigma_t(i)$ by
\[
\Theta_i^f=
\Bigl\{
k\in[i+1,n]\ \Bigm|\ 
\sum_{j=i+1}^{k}
\Bigl(
p_{\sigma_t(j)}-
\bigl(p_{\sigma_t(i)}-\min(w_{\sigma_t(i)},0)\bigr)
\Bigr)\le 0
\Bigr\},
\qquad i=1,\dots,n-1.
\]
Because queue discontinuities are resolved by the bottleneck breakthrough rule, it is sufficient to retain only the boundary candidates required to activate that rule in the candidate set. Because queue discontinuities are handled separately by the bottleneck breakthrough rule, $\theta_{1,2}^f$ only needs to be retained in the candidate set to initiate that rule.

The family of all forward candidate sets, together with the discontinuity positions, forms the complete family of forward candidate solutions:
\[
\Theta_f=\{\Theta_i^f\}_{i=1}^{n-1}\cup \{\, k \mid w_{\sigma_t(k)} < 0 \,\}.
\]

We next apply \cref{Theorem2} to identify the corresponding structure for backward candidate sets. The detailed derivation is provided in the Appendix	~\ref{app:backward condition}.

\begin{lemma}[Criterion for the Backward Candidate Set]
\label{backward condition}
For a schedule $\sigma_t$, consider moving job $\sigma_t(i)$ backward from position $i$ to position $j$, where $j \in [1,i-1]$. If position $j$ satisfies the necessary conditions
\begin{align*}
w_{\sigma_t(i)}
&\ge
\sum_{k=j}^{i-1}
\bigl(
p_{\sigma_t(k)}-\min(0,w_{\sigma_t(k)})
\bigr)
\;\land\;
p_{\sigma_t(i)}
\ge
\sum_{k=j}^{i-1}
\bigl(
-\min(0,w_{\sigma_t(k)})
\bigr),
\end{align*}
then position $j$ is called a \emph{backward candidate solution} of job $\sigma_t(i)$, denoted by $\theta_{i,j}^b$.
\end{lemma}

We define the backward candidate set associated with job $\sigma_t(i)$ by
\begin{align*}
\Theta_i^b = \{ j \mid j \in [j^*, i-1] \},
\end{align*}
where
\begin{align*}
j^* = \min \Bigl\{ j \in [1,i-1] \mid &
\bigl(
w_{\sigma_t(i)} \leq \sum_{k=j}^{i-1} \bigl( p_{\sigma_t(k)} - \min(0, w_{\sigma_t(k)}) \bigr)
\bigr)
\\
& \land
\bigl(
p_{\sigma_t(i)} \leq \sum_{k=j}^{i-1} \bigl( - \min(0, w_{\sigma_t(k)}) \bigr)
\bigr)
\Bigr\},\qquad i=2,\dots,n.
\end{align*}

Since the above two cumulative sums are monotonically decreasing with respect to $j$, once $j^*$ satisfies the conditions, they remain satisfied for any $j \ge j^*$. Hence, the backward candidate set forms a contiguous interval.

The family of all backward candidate sets forms the complete family of backward candidate solutions:
\[
\Theta_b=\{\Theta_i^b\}_{i=2}^{n}.
\]

By \cref{forward condition} and \cref{backward condition}, if a forward move or a backward move can reduce the total waiting time, then it must satisfy the corresponding necessary conditions. Since a negative ideal improvement value means that the move has the potential to reduce the total waiting time, it must satisfy the above conditions. Therefore, any improvement direction with negative ideal improvement value must be contained in either the forward candidate set or the backward candidate set. This guarantees that the candidate-set screening procedure does not exclude any potentially improvement direction.

\subsection{Independence of Improvable Ideal Directions}

This subsection analyzes the relationship among improvable ideal directions. By \cref{Lemma1}, if the decrease flows induced by multiple job moves act on the same region, then these flows must compete with each other, because the positive waiting times in that region admit only limited room for further reduction. By contrast, increase flows acting on the same region may accumulate and jointly raise the waiting times of the affected jobs. Therefore, the main structural difficulty of the problem lies in the combined effect of multiple increase-flow regions. We next analyze the structural role of such regions.

\begin{lemma}
\label{increase flow}
Let $P_I$ denote the increase-flow region induced by a job-position adjustment. If there exists no feasible adjustment within $P_I$ that can strictly reduce the total waiting time, then under the action of increase flow, any adjustment within $P_I$ can only mitigate the additional loss caused by the increase flow, but cannot serve as an independent source of improvement in the total waiting time.
\end{lemma}

\begin{proof}
Consider the increase-flow region $P_I$ induced by a given job-position adjustment. Suppose that there exists no feasible adjustment within $P_I$ that can strictly reduce the total waiting time. We show that, under the action of increase flow, any adjustment inside $P_I$ can only reduce the additional loss caused by the increase flow, but cannot generate a net improvement in the total waiting time by itself.

Take an arbitrary consecutive subinterval of $P_I$,
\[
P_a=\{\sigma_t(l),\dots,\sigma_t(k)\}.
\]
Its completion time can be written as
\[
C_{P_a}=C_{P_o}+\sum_{u\in P_a}p_u+\sum_{u\in P_a}\bigl(-\min(0,w_u)\bigr),
\]
where $C_{P_o}$ is the completion time of the job immediately preceding the set $P_a$. Since the set of jobs contained in $P_a$ is fixed, the term $\sum_{u\in P_a}p_u$ is constant. If an adjustment increases the total machine idle time inside $P_a$ by $\Delta>0$, then
\[
C_{P_a}'-C_{P_a}=\Delta.
\]

Therefore, the completion time of $P_a$ is delayed by exactly $\Delta$. Since the completion time of $P_a$ serves as the starting-time reference for the service of its succeeding region, this delay generates a new increase-flow input of magnitude $\Delta$ after $P_a$. By \cref{Lemma2}, this increase flow can only propagate downstream along the sequence, and can at most be partially or fully absorbed by subsequent negative waiting-time components.

This shows that attempting to reduce the original increase flow in the downstream region by increasing the machine idle time inside $P_a$ merely reinjects an equal amount of delay into the propagation process after $P_a$, rather than creating any new net gain. Combined with the premise of the lemma, namely, that there exists no feasible adjustment within $P_I$ that can strictly reduce the total waiting time, it follows that adjustments inside $P_I$ cannot constitute an independent source of net improvement in the total waiting time. At most, they can only mitigate the additional total waiting time caused by the increase flow.

Hence, if no feasible adjustment within $P_I$ can strictly reduce the total waiting time, then under the action of increase flow, any subsequent adjustment within $P_I$ can only mitigate the additional loss, but cannot serve as an independent source of net improvement in the total waiting time.
\end{proof}

Since the roles of decrease flow and increase flow are now clear, and any move induces at most one of them or both simultaneously, we are ready to study the relationship among improvable ideal directions.

\begin{theorem}[No coordination among improvable ideal directions]
\label{independent}
For any two improvable ideal directions, their net sources of improvement cannot be mutually enhanced through coordination. Equivalently, no coordination relationship exists among improvable ideal directions.
\end{theorem}

\begin{proof}
By \cref{Theorem1,Theorem2}, the execution of any ideal move induces at most a decrease flow, an increase flow, or both.

By \cref{increase flow}, if an increase-flow region contains no independently feasible adjustment that can strictly decrease the total waiting time, then any subsequent adjustment under the action of that increase flow can at most reduce the additional loss caused by the increase flow, but cannot serve as a net source of improvement in the total waiting time. Therefore, for any ideal direction that yields a strict improvement, its net gain can only come from the release of benefit in the corresponding decrease-flow region.

On the other hand, by \cref{Lemma1}, if the decrease flows induced by multiple job moves act on the same region, then these flows must constrain each other, because the positive waiting times in that region admit only limited room for further reduction. Hence, different decrease flows cannot simultaneously serve as additional benefit sources for one another. Equivalently, the benefit released by one improvable ideal direction through its decrease flow cannot generate a synergistic gain with that of another improvable ideal direction.

Therefore, no coordination relationship exists among improvable ideal directions.
\end{proof}

This theorem implies that each improvable ideal direction can be treated as an independent improvement unit, without relying on coordination with other improvable ideal directions. Therefore, if every improvable ideal direction can be realized in finitely many steps, then, in principle, a finite-step algorithm for obtaining an optimal solution exists.

\subsection{Structural Obstacles to Improvement Paths}

For any non-optimal schedule $\sigma_t$, there exists at least one ideal forward or backward improvement direction induced by an optimal solution. However, such an ideal direction is not necessarily directly executable. The core reason is that, during the move, changes in waiting times propagate along the queue and may cause the corresponding move to lose feasibility. It is therefore necessary to characterize the structural factors that obstruct the realization of ideal directions.

\begin{lemma}[Obstacle structures]
\label{zhuai}
Let the current schedule be $\sigma_t$, and suppose that there exists an ideal direction $M^*$ induced by some optimal solution such that its ideal improvement value satisfies $Z(M^*)<0$. If the corresponding ideal move fails to realize its ideal improvement when actually executed, then the structural obstruction can only arise from a discontinuity in the propagation of the decrease flow. More specifically, such a discontinuity can only be of one of the following two types:
\begin{enumerate}
    \item \textbf{Pre-existing discontinuity}: before the move, there exists a position $l$ in the propagation range of the decrease flow such that
    \[
    w_{\sigma_t(l)} \le 0;
    \]
    \item \textbf{Post-move discontinuity}: after the move, there exists a position $l$ in the propagation range of the decrease flow such that
    \[
    w_{\sigma_t(l)}^2 \le 0.
    \]
\end{enumerate}
Therefore, any structure that obstructs the realization of an improvable ideal direction must, and can only, take the form of one of the above two types of discontinuity.
\end{lemma}

\begin{proof}
By \cref{Theorem1,Theorem2}, a job move generates only two propagation effects, namely, the decrease flow and the increase flow. The decrease flow reduces the waiting times of subsequent jobs and is therefore the direct carrier of objective improvement, whereas the increase flow raises the waiting times of some jobs and thus creates additional loss. Hence, if an improvable ideal direction fails to achieve its ideal improvement in actual execution, the key issue is whether the decrease flow can be propagated as required by the ideal move.

We first consider the decrease flow. By \cref{Lemma1}, the propagation of the decrease flow at each position depends only on the original waiting time at that position and on the amount of decrease flow received there. In particular, if at some position $l$,
\[
w_{\sigma_t(l)} \le 0,
\]
then the decrease flow is absorbed at that position and cannot continue to propagate further downstream. As a consequence, the waiting-time reductions that should have been transmitted to the jobs after position $l$ are lost, and the actual improvement becomes strictly weaker than the ideal one. Therefore, if there already exists such a position in the propagation range before the move, this pre-existing discontinuity constitutes a structural obstruction to the realization of the ideal improvement.

Next, even if no discontinuity exists initially in the relevant propagation range, a new discontinuity may be created by the move itself. More precisely, suppose that for some position $l$ one has $w_{\sigma_t(l)}>0$ before the move, but
\[
w_{\sigma_t(l)}^2
=
w_{\sigma_t(l)}-\Delta w_{\sigma_t(l)}^-
\le 0
\]
after the move. Then position $l$ becomes a new cutoff point of the decrease flow. Consequently, the decrease flow can no longer propagate beyond $l$ as required by the ideal move, and part of the intended downstream improvement is lost. Hence, a post-move discontinuity is likewise a structural obstruction to the realization of the ideal improvement.

We then turn to the increase flow. The increase flow itself does not generate improvement; it only represents additional loss. By \cref{increase flow}, any subsequent adjustment within the increase-flow region can at most mitigate this additional loss, but cannot by itself become an independent source of net improvement. Therefore, even when the loss caused by the increase flow is not sufficiently eliminated, the underlying reason is still not the increase flow itself, but rather the failure of some corresponding decrease flow to propagate far enough to offset it. In other words, the failure of compensation in the increase-flow region is still reducible to an obstruction in the propagation of the decrease flow.

Combining the above observations, we conclude that any failure in realizing an improvable ideal direction must ultimately be caused by incomplete propagation of the decrease flow; and such a failure can occur if and only if the propagation range contains either a pre-existing discontinuity or a post-move discontinuity. This proves the claim.
\end{proof}

\subsection{Repair of Queue Discontinuities and a Complete Characterization of Local Optimality Types}

This subsection first characterizes the structural conditions that a job must satisfy in order to repair a queue discontinuity. Building on this result, we show that the realizability of every improvable ideal direction can be determined in finitely many steps, and finally establish a complete characterization of the local-optimum types.

\begin{lemma}[Repair Conditions at a Discontinuity]
\label{Lemma3}
Given a job sequence $\sigma_t$, suppose that the queue is discontinuous at job $\sigma_t(i)$, that is,
\[
w_{\sigma_t(i)} \le 0.
\]
Consider moving some job $\sigma_t(j)$, where $j>i$, to position $i$ in order to repair this discontinuity. If this move can indeed repair the discontinuity at position $i$, then it must satisfy
\[
r_{\sigma_t(j)} < r_{\sigma_t(i)}
\qquad\text{and}\qquad
w_{\sigma_t(j)} > 0.
\]
\end{lemma}

\begin{proof}[Proof sketch]
First, if one only adjusts the part of the sequence before job $\sigma_t(i)$, then job $\sigma_t(i)$ is affected only by the increase flow. For the increase flow, queue discontinuities do not obstruct the realization of the ideal move, and therefore no repair is needed there. What truly requires repair is the region affected by the decrease flow, which means that the repairing job must be chosen from positions after $\sigma_t(i)$.

Second, for any $j>i$, \cref{Theorem2} gives the expression for the waiting-time variation when job $\sigma_t(j)$ is moved to the position immediately before $i$. Combining this with \cref{Proposition1}, one can conclude that repairing the discontinuity requires
\[
r_{\sigma_t(j)} < r_{\sigma_t(i)}.
\]
The full derivation is given in the Appendix~\ref{app:Lemma3}.
\end{proof}

\cref{Lemma3} has characterized the necessary conditions for repairing a discontinuity, thereby restricting the feasible repairing jobs to a finite candidate set. Based on this finiteness, we next prove that the realizability of every improvable ideal direction can be determined in finitely many steps.

\begin{lemma}[Polynomial-time decidability of improvable ideal directions]
\label{finite ideal}
For any improvable ideal direction, its realizability can be decided in polynomial time. Moreover, if the direction can be restored as a strict improvement, then a corresponding realization path can also be constructed in polynomial time.
\end{lemma}

\begin{proof}
By \cref{Theorem1,Theorem2}, any job move can induce either the decrease flow, the increase flow, or both. By \cref{zhuai}, any obstruction to the realization of an improvable ideal direction must ultimately be reduced to a discontinuity in the decrease flow induced by that direction. By \cref{Lemma3}, the repair candidates for each such discontinuity can only come from a finite set of jobs. Furthermore, by \cref{forward condition,backward condition}, both the forward candidate set and the backward candidate set are of polynomial size. In the worst case, the total numbers of forward and backward candidates are each bounded by
\[
\sum_{i=1}^{n-1}(n-i)=O(n^2).
\]

For any candidate move, the number of explicit discontinuities induced by the associated decrease flow is at most $O(n)$, and the number of feasible repair candidates for each discontinuity is also at most $O(n)$. Hence, the total number of discontinuity-repair checks associated with a given direction is polynomially bounded. If a single-job repair is insufficient to restore the ideal direction, then the corresponding combined compensation structures are still constrained by the traversal range of the candidate set. By \cref{LemmaCoverB}, any potentially effective combined compensation structure will be exposed during the sequential traversal of the forward candidate set. Therefore, the algorithm does not need to enumerate all subsets of jobs within that range; it only needs to perform finitely many extensions and checks along candidate positions. Consequently, combined compensation does not introduce any additional exponential combinatorial search. On the other hand, by \cref{increase flow}, the increase flow cannot serve as an independent source of net improvement, and its treatment can ultimately still be reduced to checking repair candidates for the decrease flow.

In addition, the feasibility and objective-value change associated with each candidate move, discontinuity repair, and combined compensation only involve the propagation of waiting-time changes over the corresponding local interval. By the propagation relations given in \cref{Theorem1,Theorem2}, these waiting-time changes can be updated linearly along the local interval. Hence, the computational cost of each individual check is $O(n)$.

Therefore, the candidate-move checks, discontinuity-repair checks, and combined-compensation checks required for any improvable ideal direction all lie within a polynomial-size search space whose elements can be evaluated one by one. It follows that the realizability of any improvable ideal direction can be decided in polynomial time. Moreover, if the answer is affirmative, then the corresponding realization path can be recovered in polynomial time by following the same constructive checking procedure.
\end{proof}

Based on the above discussion, we next classify the local-optimum types and show that this classification is complete. Specifically, \cref{Proposition3} will show that every non-global-optimal local structure that lies on an improvement path but cannot directly realize a strict improvement

\begin{proposition}[Classification and completeness of local optimum types]
\label{Proposition3}
Let $\sigma_t$ be the current schedule and let $\sigma^*$ be an optimal schedule, with $\sigma_t\neq \sigma^*$. Suppose that there exists a job $J$ whose positions in $\sigma_t$ and $\sigma^*$ are $j$ and $j^*$, respectively. Let $Z$ be the ideal objective-value change induced by the improvable ideal direction of $J$ under $\sigma^*$, let $Z'$ be the objective-value change obtained by executing this direction directly on $\sigma_t$, and let $\bar Z$ be the objective-value change obtained after finitely many executable adjustments. If
\[
Z<0,\qquad Z'\ge 0,\qquad \bar Z<0,
\]
then the corresponding non-global-optimal local structure belongs uniquely to exactly one of the following two types:
\begin{enumerate}
    \item if $j^*>j$, then the blocked ideal direction is forward, and the structure is a forward-type local optimum;
    \item if $j^*<j$, then the blocked ideal direction is backward, and the structure is a backward-type local optimum.
\end{enumerate}
Moreover, no other independent local optimum type can arise from an obstructed ideal improvement.
\end{proposition}

\begin{proof}
Assume that $\sigma_t$ is not globally optimal and that there exists a job $J$ such that
\[
Z<0,\qquad Z'\ge 0,\qquad \bar Z<0.
\]
The condition $Z<0$ means that the ideal direction induced by $\sigma^*$ is theoretically improving, $Z'\ge 0$ means that it cannot be realized as a strict improvement directly from the current schedule, and $\bar Z<0$ means that it can still be converted into a strict improvement after finitely many executable adjustments.

Since the position of $J$ in $\sigma_t$ differs from that in $\sigma^*$, its ideal direction is uniquely determined: it is forward if $j^*>j$, and backward if $j^*<j$. Hence, the associated non-optimal local structure must already be tied to one definite ideal direction, and no third direction type is possible.

The condition $Z'\ge 0$ implies that this ideal direction is blocked by some structural obstacle in the current schedule. On the other hand, $\bar Z<0$ shows that the obstacle is removable, since after finitely many executable adjustments the same ideal direction can still be turned into a strict improvement. Therefore, the type of the local structure is determined by which ideal direction is blocked: a blocked forward ideal direction gives a forward-type local optimum, and a blocked backward ideal direction gives a backward-type local optimum. Any auxiliary adjustments appearing along the feasible path only help release the improvement potential of that same blocked ideal direction, and thus do not change the type classification.

Therefore, every non-global-optimal local structure satisfying
\[
Z<0,\qquad Z'\ge 0,\qquad \bar Z<0
\]
belongs uniquely to either the forward-type or the backward-type class. The classification is thus complete.
\end{proof}

By \cref{independent} and \cref{finite ideal}, improvable ideal directions are mutually independent, and the realizability of each such direction can be decided in finitely many steps. In addition, \cref{Proposition3} gives a complete classification of the local-optimum types. These results together yield a finite-step search framework built around improvable ideal directions and provide the theoretical foundation for the subsequent development of a globally optimal algorithm.

\section{Algorithm Design}

We begin by outlining the central idea of the proposed algorithm. By \cref{independent}, improvable ideal directions require no coordination and can therefore be treated as independent improvement units. In addition, \cref{finite ideal} shows that the realizability of each such direction can be decided in finitely many steps. Moreover, \cref{zhuai} establishes that queue discontinuity is the unique structural obstacle to the realization of improvable ideal directions, and \cref{Proposition3} gives a complete characterization of the local-optimum types. The key task in the subsequent algorithm design is thus to construct suitable repair operators for each type of local-optimum structure, so as to progressively eliminate these structural obstacles and restore the ability of the schedule to keep improving along improvable ideal directions.

Accordingly, we design a repair operator for each local optimum type and prove that every such operator is polynomial-time implementable and yields a strict decrease in the objective value. We further show that the iterative application of these operators progressively shrinks the set of jobs involved in structural violations, thereby ensuring finite termination.

We next introduce the basic definitions used throughout the algorithm. Taking the forward candidate set $\theta_f$ as an example, each candidate operation partitions the current job sequence into several structured subsequences. Consider the current sequence $\sigma_t$, and let jobs $i$ and $j$ satisfy $j>i$. Suppose that job $i$ is moved forward to the position immediately after job $j$. Then, except for job $i$, the remaining jobs preserve their original relative order and can be partitioned into the following three subsequences:
\[
P_O=\{\sigma_t(k)\mid k<i\},\quad
P_D=\{\sigma_t(k)\mid i<k\le j\},\quad
P_I=\{\sigma_t(k)\mid j<k\le n\},
\]
where $P_O$ denotes the subsequence preceding job $i$, $P_D$ denotes the subsequence between jobs $i$ and $j$ (excluding $i$ and including $j$), and $P_I$ denotes the subsequence following job $j$. Any of these subsequences may be empty.

Similarly, when a job is moved backward, the sequence is again partitioned into three parts, except that the relative positions of $P_I$ and $P_D$ are interchanged, so that the order becomes $P_O$, $P_I$, and $P_D$. Regardless of whether the move is forward or backward, its effect can be characterized in a unified manner through the interaction among these three subsequences.

\subsection{Machine idle time as the medium connecting $P_I$ and $P_D$}

By \cref{Theorem1}, the function $I_2$ characterizes the effect transmitted from $P_D$ to $P_I$ when a job is moved forward. Similarly, in the backward-move case, \cref{Theorem2} shows that the effect transmitted from $P_I$ to $P_D$ is characterized by $I_4$. However, in order to obtain an improved ordering within $P_I$ and $P_D$, the jobs inside these regions must in general be reordered, and such reordering introduces an additional effect on top of $I_2$ or $I_4$.

We first consider the forward-move case. Here, $I_2$ represents the effect transmitted from $P_D$ to $P_I$ before any additional algorithmic adjustment is applied. Introduce the notion of \emph{completion time} $C$, namely, the completion time of the last job in the queue. If $P_D$ and $P_I$ are each viewed as an aggregate block, then the completion time of $P_D$ is denoted by $C_{P_D}$.

By definition, $C_{P_D}$ can be written as
\[
C_{P_D}=C_{P_O}+\sum_{k\in P_D} p_k+\sum_{k\in P_D}\bigl(-\min(0,w_k)\bigr),
\]
where $C_{P_O}$ is the completion time of $P_O$, $\sum_{k\in P_D} p_k$ is the total processing time of the jobs in $P_D$, and the last term is the total machine idle time in $P_D$. Since both $C_{P_O}$ and the total processing time of $P_D$ are fixed, the variation in $C_{P_D}$ is determined solely by the total machine idle time in $P_D$, namely$\sum_{k\in P_D}\bigl(-\min(0,w_k)\bigr).$

Whenever the total machine idle time in $P_D$ changes, an additional adjustment is transmitted to $P_I$. This quantity is given by
\[
\Delta w_{minus}
=
\sum_{k\in P_D}\min(0,w_k^2)
-\sum_{k\in P_D}\min(0,w_k^3),
\]
where $w_k^2$ denotes the waiting time obtained after executing only the job move, whereas $w_k^3$ denotes the waiting time obtained after the subsequent reordering rule is applied.

Therefore, in the forward-move case, the total effect transmitted from $P_D$ to $P_I$ is given by
\[
I_2+\Delta w_{minus},
\]
and the waiting time of the forward-moved job $i$ is correspondingly updated as
\[
w_i^3=w_i^2+\Delta w_{minus}.
\]

In the backward-move case, the reordering adjustment takes place within $P_I$, and the resulting variation in machine idle time acts in the opposite direction on the starting time of $P_D$. Hence, the effect received by $P_D$ is expressed as
\[
I_4-\Delta w_{minus}.
\]
In this case, if the backward-moved job is denoted by $i$, then its waiting time remains unchanged.

Based on the above unified representation, the subsequent subsections introduce the concrete correction operators used in the algorithm, including  bottleneck breakthrough, and adjacent interchange. Together, these operators form the complete improvement mechanism of the proposed exact algorithm.

\subsection{Adjacent interchange rule for the increase flow}

For the processing region $P_I$ operating under increase flow, we introduce an \emph{adjacent interchange rule}. In the forward- and backward-adjustment cases, the initial value of the increase flow in $P_I$ is given by $\max(0,I_2)$ and
\[
\max\left(p_{\sigma_t(i)},\ \sum_{j=k}^{i} p_{\sigma_t(j)}-\sum_{j=k}^{i-1}\min(0,w_{\sigma_t(j)})-w_{\sigma_t(i)}\right),
\]
respectively. For notational convenience, we denote the initial increase flow by $f_I$, and let $\pi_{\sigma_t}$ represent the total amount of flow associated with sequence $\sigma_t$.

We first describe the operational procedure of the adjacent interchange rule. Let the current job sequence be $\sigma_t=\{\sigma_t(1),\dots,\sigma_t(n)\}$. Consider a finite subsequence $P_I=\{\sigma_t(i),\dots,\sigma_t(k)\},\; 1\le i<k\le n.$

\textbf{Step 1.}
Starting from position $i$, scan the sequence in the forward direction. For any adjacent pair of jobs $\sigma_t(j)$ and $\sigma_t(j+1)$, where $j\in[i,k-1]$, if
\begin{align*}
p_{\sigma_t(j)} > p_{\sigma_t(j+1)}
\quad \text{and} \quad
\Delta w_{P_I}+\Delta \pi_{P_I}<0,
\end{align*}
then interchange this adjacent pair.

Let $\sigma_t'$ denote the sequence after the interchange. The corresponding variations are defined by
\begin{align*}
\Delta w_{P_I}
&=
\sum_{h=j}^{n}\max\!\bigl(0,w_{\sigma_t'(h)}\bigr)
-
\sum_{h=j}^{n}\max\!\bigl(0,w_{\sigma_t(h)}\bigr),\\
\Delta \pi_{P_I}
&=
\pi_{\sigma_t'}-\pi_{\sigma_t},
\end{align*}
where $\pi_{\sigma_t}$ denotes the total flow in the current sequence, and $\pi_{\sigma_t'}$ denotes the total flow after the interchange.

\textbf{Step 2.}
If an interchange is performed, update $\sigma_t\leftarrow\sigma_t'$ and return to Step 1. If no adjacent pair satisfies the above conditions, terminate the procedure.

This rule provides a unified characterization, from the perspective of increase flow, of the trade-off between the waiting time and the total increase flow in both forward and backward adjustments. On the one hand, increase flow enlarges the waiting times of jobs in $P_I$; on the other hand, it also amplifies the benefit that can be generated by a job move (see \cref{Lemma1}). Therefore, by evaluating the sign of $\Delta w_{P_I}+\Delta \pi_{P_I}$, one can determine whether a given adjacent interchange is acceptable in the sense of reducing the total waiting time.

As shown in \cref{Lemma3}, any job capable of repairing a queue discontinuity must satisfy the nondecreasing release-time structure. Since any forward move or backward move can essentially be represented as a composition of adjacent interchanges, it is therefore necessary to analyze the structural properties of adjacent job pairs ordered nondecreasingly by release times.

\begin{lemma}[Conditions for an adjacent interchange to yield a no-aftereffect improvement]
\label{Lemma4}
Consider a job sequence $\sigma_t = \{\sigma_t(1),\dots,\sigma_t(n)\}$. For any adjacent job pair $(\sigma_t(j),\sigma_t(j+1))$ ordered nondecreasingly by release times, interchanging these two adjacent jobs reduces their total waiting time if and only if
\[
\Delta w_{positive}<0
\iff
\left(w_{\sigma_t(j+1)}>\frac12\bigl(p_{\sigma_t(j)}+p_{\sigma_t(j+1)}\bigr)\right)
\land
\left(p_{\sigma_t(j)}>p_{\sigma_t(j+1)}\right),
\]
where $\Delta w_{positive}
=
\max\bigl(0,w_{\sigma_t(j)}^2\bigr)
+\max\bigl(0,w_{\sigma_t(j+1)}^2\bigr)
-\max\bigl(0,w_{\sigma_t(j)}\bigr)
-\max\bigl(0,w_{\sigma_t(j+1)}\bigr).$

Furthermore, if the additional condition
\[
w_{\sigma_t(j+1)}^2\ge \min(0,w_{\sigma_t(j)})
\]
holds, then the effect of this interchange does not propagate to subsequent jobs. Hence, the waiting times of all downstream jobs remain unchanged, and the interchange is therefore a \emph{no-aftereffect improvement}.
\end{lemma}

\begin{proof}[Proof sketch]
For any adjacent job pair ordered nondecreasingly by release times, the conditions given in \cref{Proposition1} are automatically satisfied. By comparing the total waiting time of the two jobs and the variation in machine idle time before and after the interchange, one can derive both the necessary and sufficient condition for improvement and the condition for the absence of aftereffect. The derivation is straightforward and is therefore deferred to the Appendix~\ref{app:Lemma4}.
\end{proof}

\begin{lemma}[Global optimality of a local optimum under the adjacent interchange rule]
\label{Lemma5}
Under the adjacent interchange rule, consider a schedule whose initial sequence is ordered nondecreasingly by release times. If, for every adjacent pair satisfying $p_{\sigma_t(j)} > p_{\sigma_t(j+1)},$ one has $w_{\sigma_t(j+1)}^2>0,$ then any local optimum under this rule must be globally optimal.
\end{lemma}

\begin{proof}
If an adjacent pair satisfies $p_{\sigma_t(j)} > p_{\sigma_t(j+1)}$ and $w_{\sigma_t(j+1)}^2>0$, then the conditions of \cref{Lemma4} are satisfied. Hence, interchanging the adjacent pair $(\sigma_t(j),\sigma_t(j+1))$ constitutes a no-aftereffect improvement, that is, the interchange strictly reduces the total waiting time and does not alter the waiting times of jobs located after the interchange position.

Therefore, such an interchange does not change the discontinuity pattern or the downstream waiting-time propagation structure. Since the processing-time order of this pair changes from inverted to correctly ordered after the interchange, the same pair will not trigger another interchange. Meanwhile, the decision environment for all other untreated adjacent inverted pairs remains unchanged. By induction, repeated application of this rule eliminates all adjacent inverted pairs satisfying $p_{\sigma_t(j)} > p_{\sigma_t(j+1)}$, and the process terminates after finitely many steps.

Let the terminal sequence be $\sigma'$. Then there exists no adjacent pair in $\sigma'$ satisfying
\[
p_{\sigma'(j)} > p_{\sigma'(j+1)}.
\]
Hence, within each continuous subsequence, the jobs are ordered nondecreasingly by processing times.

On the other hand, because all accepted interchanges are no-aftereffect improvements, the discontinuity structure remains unchanged throughout the process. Therefore, the sequence can always be partitioned into several continuous subsequences, each with a fixed starting time. Consider any such continuous subsequence with starting processing time equal to a constant $T$. Since there is no machine idle time inside this block, the completion times of jobs in the block depend only on their relative order. Thus, the sequencing problem inside this block is equivalent to the classical single-machine scheduling problem without release times, namely $1||\sum C_j$, whose optimal order is the nondecreasing processing-time order (the SPT rule). Since minimizing total completion time is equivalent to minimizing total waiting time, the total waiting time of each continuous subsequence is minimized when the jobs are arranged in nondecreasing order of processing times.

In addition, because the discontinuity structure remains unchanged, the machine idle time throughout the process remains identical to that of the initial sequence. Since the initial sequence is ordered nondecreasingly by release times, it attains the minimum machine idle time. Hence, the machine idle time remains minimal throughout the entire process.

Combining the above facts, the terminal sequence $\sigma'$ preserves the discontinuity structure, and hence the minimum machine idle time, while also achieving the optimal order within each continuous subsequence. Therefore, the total waiting time of $\sigma'$ is globally minimal, which implies that $\sigma'$ is a global optimum.

Thus, under the adjacent interchange rule, any local optimum must be globally optimal.
\end{proof}

\begin{corollary}
\label{corollary1}
The contrapositive of \cref{Lemma5} implies the following: under the adjacent interchange rule, if a local optimum is not globally optimal, then there must exist at least one adjacent pair satisfying $p_{\sigma_t(j)} > p_{\sigma_t(j+1)}$ such that the queue becomes discontinuous after the interchange, that is,
 $w_{\sigma_t(j+1)}^2\le 0.$
\end{corollary}

However, the adjacent-interchange structure characterized by \cref{corollary1} alone is still insufficient to cover all cases in which a schedule is locally optimal but not globally optimal under the adjacent interchange neighborhood. Indeed, when every adjacent pair satisfies $p_{\sigma_t(j)} \le p_{\sigma_t(j+1)},$ the adjacent interchange rule is no longer triggered. Nevertheless, some jobs may still admit a forward move across multiple positions that further reduces the total waiting time, namely, an ideal forward move whose actual improvement value is negative.

\begin{proposition}
\label{Proposition4}
Therefore, after the adjacent interchange rule has been exhaustively applied, if a schedule is locally optimal in this neighborhood but is not globally optimal, then its nonoptimality must originate from one of the following two structures:
\begin{enumerate}
    \item an ideal forward move whose actual improvement value is negative, referred to as \emph{under-consumption};

    \item at least one adjacent pair satisfying $p_{\sigma_t(j)} > p_{\sigma_t(j+1)},$
    such that their interchange makes the queue discontinuous, referred to as \emph{over-consumption}.
\end{enumerate}

These two structures together provide a complete structural characterization of solutions that are locally optimal but not globally optimal under the adjacent interchange neighborhood. The former can be resolved by traversing ideal forward moves, whereas the latter corresponds to a backward-type local optimum. Hence, no new local-optimality type is introduced.
\end{proposition}

\begin{figure}[h]
    \centering
    \includegraphics[width=0.9\linewidth]{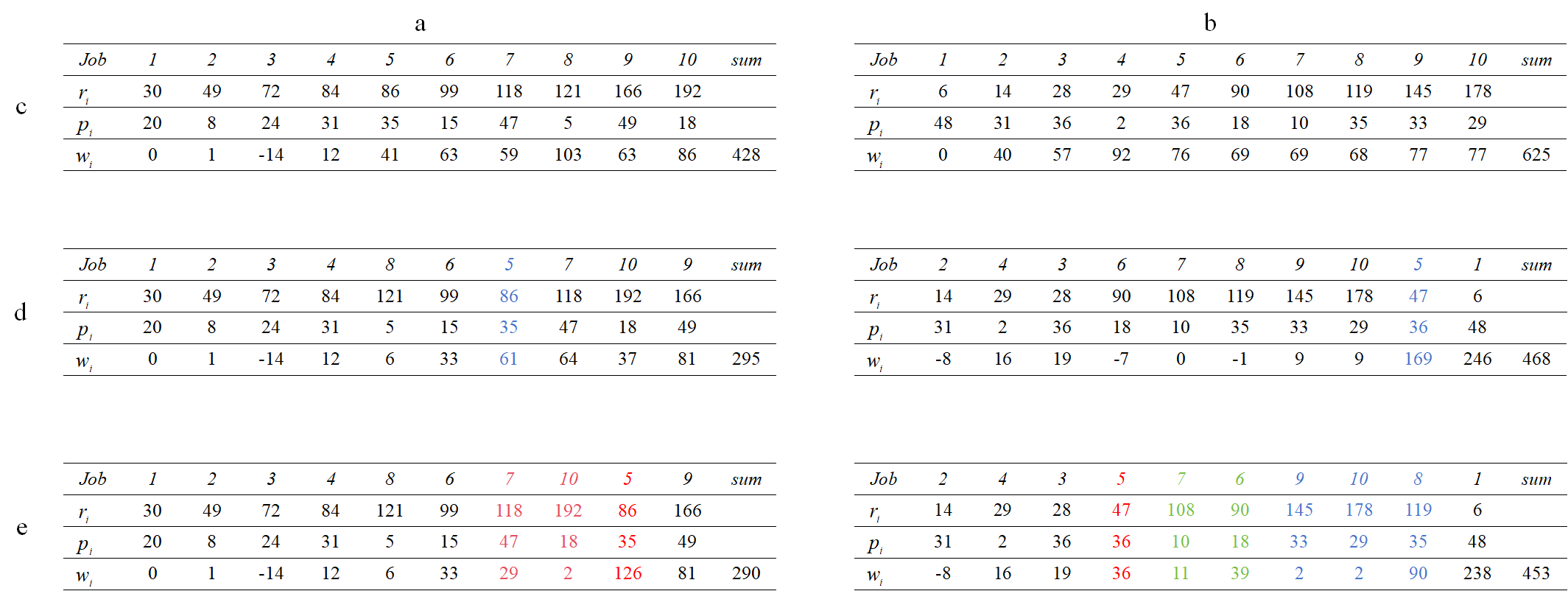}
    \caption{Two illustrative examples are provided. Panels (a) and (b) correspond to the two cases of under-consumption and over-consumption, respectively. Rows (c), (d), and (e) show, respectively, the sequence initially ordered nondecreasingly by release times, the result after applying the adjacent interchange rule, and the optimal sequence. In panel (a), since $p_5 < p_7$, the best result attainable by the adjacent interchange rule is the sequence shown in row (d), whereas the true optimum requires job 5 to be placed immediately after job 10. In panel (b), job 5 is located after job 10, whereas the optimal sequence requires job 5 to be scheduled before job 6 and job 8 to be scheduled after job 10.}
    \label{adjoint_exchange}
\end{figure}

\subsection{Bottleneck breakthrough rule for the decrease flow}

After handling the increase-flow region, we introduce the \emph{bottleneck breakthrough rule} to deal with decrease flow. It should be emphasized that this rule is applicable only to sequences that have already been processed by either the adjacent interchange rule or the solution-set traversal procedure. The initial values of decrease flow in the forward- and backward-move cases are $-(p_{\sigma_t(i)}-\min(0,w_{\sigma_t(i)})),\;\min(0,I_2),\; \text{and}\; I_4,$ which are uniformly denoted by $f_D$. Moreover, let $\pi_{\sigma_t}$ denote the total amount of the decrease flow of $\sigma_t$.

We next formalize the bottleneck breakthrough rule. Let the reference sequence be $\sigma_t^{o}=\{\sigma_t(1),\dots,\sigma_t(n)\}.$ Consider a finite subsequence $P_D=\{\sigma_t(i),\dots,\sigma_t(k)\},\; 1\le i\le k\le n,$ which represents the decrease-flow region currently under treatment.

The bottleneck breakthrough rule proceeds as follows.

\textbf{Step 1 (Termination condition 1).}
If $w_{\sigma_t(j)}+f_D\ge 0,\qquad \forall j\in[i,k],$ then the decrease flow does not induce any queue discontinuity, and the procedure terminates.

\textbf{Step 2 (Bottleneck identification).}
Otherwise, locate the first bottleneck position
\[
l^*=\min\left\{\,l\in[i,k-1]\mid w_{\sigma_t(l)}+f_D<0\,\right\}.
\]

\textbf{Step 3 (Adjustment step).}
Choose a job $\sigma_t(g)$ and move it to the position immediately before $\sigma_t(l^*)$, provided that the following conditions are satisfied:
\begin{align*}
&w_{\sigma_t(g)} > 0, \; w_{\sigma_t(g)}^2 > w_{\sigma_t(l^*)} \; (g>l^*),\\
&\Delta w_{P_D}+\Delta \pi_{P_D} < 0. 
\end{align*}

Let $\sigma_t'$ denote the resulting sequence after the adjustment. The corresponding variations are defined by
\begin{align*}
\Delta w_{P_D}
&=
\sum_{j=i}^{n}\max(0,w_{\sigma_t'(j)})
-
\sum_{j=i}^{n}\max(0,w_{\sigma_t^{o}(j)}), \\
\Delta \pi_{P_D}
&=
\pi_{\sigma_t'}-\pi_{\sigma_t^{o}}, 
\end{align*}
where $\pi_{\sigma_t^{o}}$ denotes the total flow before the adjustment, and $\pi_{\sigma_t'}$ denotes the total flow after the adjustment.

\textbf{Step 4 (Update rule).}
If multiple jobs satisfy conditions , select the one that minimizes
\[
\Delta w_{P_D}+\Delta \pi_{P_D},
\]
and update
\[
\sigma_t^{o}\leftarrow \sigma_t'.
\]
If no feasible job satisfies the above conditions, then reset the decrease-flow value by
\[
-f_D=\max(0,w_{\sigma_t(l^*)}).
\]

\textbf{Step 5 (Termination condition 2).}
Repeat Steps 1--4 until $l^*=k$.

The bottleneck breakthrough rule repairs queue discontinuities in the decrease-flow region according to a greedy criterion. At each iteration, it identifies the current bottleneck position and selects the job that yields the largest admissible improvement after balancing gain and cost. After each accepted adjustment, the bottleneck position is re-identified and the same procedure is repeated until no remaining discontinuity can be further improved. The structural properties underlying this rule are characterized in \cref{Lemma3}.

\subsection{Algorithm and optimality proof}

Based on the structural properties and candidate-screening conditions developed above, we are now in a position to present the pseudocode of the proposed algorithm. The algorithm iteratively selects promising moves from the forward and backward candidate sets, and repairs the queue discontinuities induced by such moves, thereby progressively decreasing the total waiting time.

We next establish that the algorithm terminates after finitely many iterations and that the solution it returns is globally optimal.

\begin{algorithm}[H]
\caption{Optimal Sort Algorithm}
\small
\label{algorithm1}
\begin{algorithmic}[1]
\State \textbf{Input:} Job set $N = \{1,2,\dots,n\}$, processing times $p_i$, release times $r_i$. Jobs are sorted in non-decreasing order of release times to obtain $\sigma_0$. Initialize $Improve = 1$, $Optimal_f = M$, and $t = 0$, where $M$ is a sufficiently large number.
\State \textbf{Output:} Optimal sequence $\sigma_t$

\While{$Improve = 1$}
    \State Set $Improve = 0$
    \State Compute the forward candidate solution set $\Theta_f$ of the current sequence $\sigma_t$, and calculate its total waiting time $w_o$. Traverse the sequence from the first job in the forward direction.
    \For{$f \in \Theta_f$}
        \For{$(job, position) \in f$}
            \State Initialize $f_D$ as $p_{job}-\min(0,w_{job})$
            \State Apply the bottleneck breakthrough rule to $P_D$
            \State Calculate the machine idle time $\Delta w_{minus}$ and update $I_2 \leftarrow I_2 + \Delta w_{minus}$
            \State Set $I_2$ as the initial value of $f_I$
            \If{$I_2 \geq 0$}
                \State Apply the adjacent interchange rule to $P_I$
            \Else
                \State Apply the bottleneck breakthrough rule to $P_I$
            \EndIf
            \State $\sigma_t^c = \text{Algorithm 2 }(\sigma_t^f, w_o, job)$
            \State Obtain the new sequence $\sigma_t^c$ and compute its total waiting time $w_c$
            \If{$Optimal_f > w_c$}
                \State Update: $\sigma_t \leftarrow \sigma_t^c$, $Improve \leftarrow 1$, $Optimal_f \leftarrow w_c$
            \EndIf
            \State $\sigma_t^{b} = \text{Algorithm 3 }(\sigma_t^c)$
            \State Obtain the sequence $\sigma_t^b$ and compute its total waiting time $w_b$
            \If{$Optimal_f > w_b$}
                \State Update: $\sigma_t \leftarrow \sigma_t^b$, $Improve \leftarrow 1$, $Optimal_f \leftarrow w_b$
            \EndIf
        \EndFor
    \EndFor
    \State $t \leftarrow t + 1$
\EndWhile
\end{algorithmic}
\end{algorithm}

\begin{algorithm}[H]
\small
\caption{The Forward Candidate Set Traversal}
\label{algorithm2}
\begin{algorithmic}[2]
\State \textbf{Input:} Job set $N = \{1,2,\dots,n\}$, processing times $p_i$, release times $r_i$, processing sequence $\sigma_t^f$, total waiting time $w_f$ of $\sigma_t^f$. Initialize $Improve = 1$, $Optimal_c = w_f$, $w_o$, and $job$.
\State \textbf{Output:} Sequence $\sigma_t^c$
\While{$Improve = 1$}
	\State Set $Improve = 0$
	\State Compute the forward candidate solution set $\Theta_f$ for $\sigma_t^f$, and traverse the sequence forward starting from the first job.
	\For{$f \in \Theta_f$}
		\For{$(job, position) \in f$}
			\State Initialize $f_D \leftarrow p_{job}-\min(0,w_{job})$
			\State Apply the bottleneck breakthrough rule to $P_D$
			\State Calculate the machine idle time $\Delta w_{minus}$ and update $I_2 \leftarrow I_2 + \Delta w_{minus}.$
			\State Compute $I_2$ and initialize $f_I \leftarrow I_2$
			\If {$I_2 \geq 0$}
				\State Apply the adjacent interchange rule to $P_I$
			\Else
				\State Apply the decrease flow law to $P_I$
			\EndIf
			\State Obtain sequence $\sigma_{t}^{'}$, apply the bottleneck breakthrough rule to $\sigma_{t}^{'}$ with $f_D=0$, skip $job$ and exclude it from queue gap compensation. Obtain sequence $\sigma_{t}^{''}$.
			\State Calculate the total waiting time $w_{tem}$ of $\sigma_{t}^{''}$
			\If{$Optimal_c > w_{tem}$}
				\State Update: $\sigma_{t}^{c} \leftarrow \sigma_{t}^{''}$, $Improve \leftarrow 1$, $Optimal_c \leftarrow w_{tem}$
			\EndIf
		\EndFor
	\EndFor
\EndWhile
\If{$Optimal_c > w_{o}$} 
	\State Update: $\sigma_{t}^{c} \leftarrow \sigma_{t}^{f}$, $Optimal_c \leftarrow w_{o}$
\EndIf
\end{algorithmic}
\end{algorithm}

\begin{algorithm}[H]
\small
\caption{The Backward Candidate Set Traversal}
\label{algorithm3}
\begin{algorithmic}[3]
\State \textbf{Input:} Job set $N = \{1,2,\dots,n\}$, processing times $p_i$, release times $r_i$, processing sequence $\sigma_t^c$. Initialize $Improve = 1$, $Optimal_b = M$, where $M$ is a sufficiently large constant.
\State \textbf{Output:} Sequence $\sigma_t^b$
\While{$Improve = 1$}
    \State Set $Improve = 0$ 
    \State Compute the backward candidate solution set $\Theta_b$ for $\sigma_t^c$, and traverse the sequence backward starting from the last job.
    \For{$b \in \Theta_b$}
        \For{$(job, position) \in b$}
            \State Initialize $f_I$ as $\max\left( p_{\sigma_t(|P_I|)}, \sum_{j=position}^{|P_I|} p_{\sigma(j)}-\sum_{j=position}^{|P_I|}\min(0,w_{\sigma(j)})-w_{\sigma_t(|P_I|)} \right)$ and apply the adjacent interchange rule to $P_I$
            \State Calculate the machine idle time $\Delta w_{minus}$ and update $I_4 \leftarrow I_4 - \Delta w_{minus}$ 
            \State Initialize $f_D \leftarrow I_4$
            \If {$I_4 \geq 0$}
				\State Apply the increase flow law to $P_D$
			\Else
				\State Apply the decrease flow law to $P_D$
			\EndIf
            \State Obtain sequence $\sigma_{t}^{'}$ and calculate its total waiting time $w_{tem}$ 
            \If{$Optimal_b > w_{tem}$} 
                \State Update: $\sigma_{t}^{b} \leftarrow \sigma_{t}^{'}$, $Improve \leftarrow 1$, $Optimal_b \leftarrow w_{tem}$
            \EndIf
        \EndFor
    \EndFor
\EndWhile
\end{algorithmic}
\end{algorithm}

In each iteration, the proposed algorithm explicitly enumerates all feasible forward and backward candidate moves, and selects an admissible move according to the maximum-improvement principle. The roles of the three algorithmic components are summarized below.

Given the current sequence $\sigma_t$ and its total waiting time $w_o$, Algorithm~\ref{algorithm1} explicitly applies each forward candidate move $\theta_{ij}^f$ in the forward candidate set of $\sigma_t$. In doing so, discontinuities that are originally hidden after the move are transformed into explicit queue discontinuities. Optimizing under $\theta_{ij}^f$ yields a sequence $\sigma_t^f$ with total waiting time $w_f$.

Algorithm~\ref{algorithm2} is then applied to $\sigma_t^f$. At each iteration, it follows a greedy rule and retains only the forward candidate with the largest improvement among the currently available forward moves, thereby producing a sequence $\sigma_t^c$ with total waiting time $w_c$. Two cases may occur.

\textbf{(1)} If no further forward move yields an improvement after Algorithm~\ref{algorithm1}, namely, if $w_c>w_o$, then no forward move can generate an actual improvement for the current sequence. In this case, the modification introduced by the forward move is retained so as to continue exploring backward moves. Hence, we set $\sigma_t^c=\sigma_t^f$ and $w_c=w_f$.

\textbf{(2)} If a further forward move still yields an improvement after Algorithm~\ref{algorithm1}, namely, if $w_c\le w_o$, then forward moves remain effective for the current sequence, and the improved sequence $\sigma_t^c$ is retained.

During the iterations of Algorithm~\ref{algorithm2}, the bottleneck breakthrough rule with $f_I=0$ is applied to the entire sequence in order to compensate for the omission of discontinuity information in the screening of the forward candidate set.

Algorithm~\ref{algorithm3} is subsequently applied to $\sigma_t^c$. Its main distinction from the previous two procedures lies in the traversal of the backward candidate set and the treatment of the $P_D$ region. Since the bottleneck breakthrough rule is essentially a backward move, and the backward candidate set may be viewed as its global counterpart, no additional treatment is imposed on the $P_D$ region. This asymmetry between forward and backward moves is deliberate: it prevents a forward move from being immediately canceled by a backward move, thereby avoiding entrapment in a local optimum.

Before proving global convergence, we first characterize the local structural patterns that may arise in solutions that are not globally optimal. We then show that each such structure can be improved by the combined neighborhood explored by the algorithm. It follows that the algorithm cannot terminate at a non-globally optimal solution, from which finite termination and global optimality are established.

\begin{lemma}[Completeness of the local-optimum classification]
\label{Lemma6}
Given a schedule $\sigma_t$, if $\sigma_t$ is not globally optimal under the greedy rule, then it must be of either forward type or backward type.
\end{lemma}

\begin{proof}
By \cref{Proposition3}, the non-global-optimal local structures arising in this problem admit only two basic types. Moreover, \cref{Proposition4} shows that the local structures induced by the adjacent swap rule are also subsumed by these same two categories. Hence, no additional local-optimum type arises, and the classification is complete.
\end{proof}

As shown by \cref{Lemma6}, any non-globally optimal solution must belong to at least one of the following two structural types: the forward-type local optimum or the backward-type local optimum. Therefore, to prove that the algorithm cannot terminate at a nonoptimal solution, it suffices to show that, for each of these two local structures, there always exists a feasible combined adjustment within the combined neighborhood enumerated by the algorithm that yields a strict improvement. Following this idea, we next discuss the breakability of forward-type and backward-type local optima, respectively.

\begin{lemma}[The combined neighborhood can escape forward-type local optimum]
\label{Lemma7}
Let $\sigma^t$ be the current schedule, and suppose that $\sigma^t$ is a forward-type local optimum. Then, within the combined neighborhood enumerated by Algorithm~\ref{algorithm2}, there exists a feasible combined adjustment that strictly decreases the total waiting time.
\end{lemma}

\begin{proof}
Since $\sigma^t$ is a forward-type local optimum, there exists an inducing job $\sigma^t(i)$ whose improvable ideal direction is forward, but this direction cannot yield a strict improvement when executed directly on the current schedule. At the same time, by the definition of a forward-type local optimum, after finitely many subsequent repair steps, the blocked forward ideal direction can still be converted into a strict improvement.

By \cref{zhuai}, the only structural obstacle preventing this ideal forward direction from being directly realized is a discontinuity arising in the propagation of the decrease flow. When $\sigma^t(i)$ is forcibly moved forward along its ideal direction, both pre-existing discontinuities and discontinuities newly created by the move become explicitly exposed in the corresponding decrease-flow region. Let $b$ denote the critical discontinuity exposed by this ideal forward direction.

Consider the local repair of the discontinuity at $b$. Under the greedy rule, the algorithm first selects the single job that provides the largest immediate compensation to the negative waiting time in the region $P_D$; denote this job by $\sigma^t(j)$. However, the fact that a single job provides the largest immediate compensation does not imply that it yields the globally best compensation structure. Indeed, the very meaning of a forward-type local optimum is that, compared with the single-job compensation currently preferred by the greedy rule, there exists a better combined repair set
\[
B=\{\sigma^t(h_1),\dots,\sigma^t(h_m)\},
\]
whose joint compensation for the same discontinuity $b$ is more effective, and hence is capable of releasing the blocked ideal forward improvement.

The key point is that this better combined repair does not lie outside the search range of the algorithm. By \cref{LemmaCoverB} in the appendix, every job in the set $B$ is contained in the forward candidate traversal range associated with the greedy compensation job $\sigma^t(j)$. Moreover, the algorithm does not enumerate forward candidates by arbitrary jumps, but expands them step by step along the forward direction. As the candidate position is progressively extended, the corresponding decrease-flow region and the explicit discontinuities within it are exposed successively. Therefore, the combined repair state induced by $B$ must be explicitly encountered at some stage of the forward candidate enumeration. Once this happens, the bottleneck breakthrough rule can reconstruct in $P_D$ a joint compensation structure that is superior to the current single-job compensation, thereby restoring the blocked ideal forward direction as a strict improving direction.

We next consider the increase-flow region. By \cref{increase flow}, if the increase-flow region $P_I$ itself cannot generate a net gain, then any subsequent adjustment within $P_I$ can only mitigate the additional loss caused by the increase flow, but cannot serve as an independent source of net improvement. Any improvement in this region must therefore rely on the coverage of the region by the decrease flow induced by other job moves. Accordingly, the algorithm minimizes this additional loss through adjacent swaps and forward candidate-set enumeration within the combined neighborhood.

Combining the above arguments, the forward candidate traversal, the bottleneck breakthrough rule, and the associated forward adjustments together form a finite combined repair path. On the one hand, this path exposes and repairs the critical discontinuity in $P_D$, thereby reconstructing a better compensation structure; on the other hand, it reduces the additional loss caused by the increase flow in $P_I$. Since, by definition of a forward-type local optimum, the blocked ideal forward direction can be converted into a strict improvement after repair, there must exist a feasible combined adjustment within the neighborhood enumerated by Algorithm~\ref{algorithm2} that produces a new schedule $\sigma'$ satisfying
\[
F(\sigma')<F(\sigma^t).
\]
Therefore, every forward-type local optimum can be broken by the combined neighborhood.
\end{proof}

\begin{lemma}[The combined neighborhood escapes backward-type local optima]
\label{Lemma8}
Let $\sigma^t$ be the current schedule, and suppose that $\sigma^t$ is a backward-type local optimum. Then, within the combined neighborhood induced by Algorithm~\ref{algorithm1} and Algorithm~\ref{algorithm3}, there exists a feasible combined adjustment that strictly decreases the total waiting time.
\end{lemma}

\begin{proof}
Since $\sigma^t$ is a backward-type local optimum, it follows from \cref{Proposition3} that there exists an inducing job $\sigma^t(i)$ whose improvable ideal direction is backward, but this ideal backward direction cannot yield a strict improvement when executed directly on the current schedule. At the same time, by the definition of a backward-type local optimum, after finitely many subsequent executable adjustments, the blocked ideal backward direction can still be converted into a strict improvement.

We first consider the decrease-flow region. By \cref{Lemma3}, repairing a discontinuity is essentially realized by moving some job backward to the corresponding discontinuity position, thereby restoring queue continuity. For any given discontinuity, the bottleneck breakthrough rule first enumerates its feasible repair jobs and then selects the currently best repair job according to the greedy criterion. Therefore, the bottleneck breakthrough rule is, in essence, a local greedy repair procedure for a single discontinuity. By contrast, the traversal of the backward candidate set is not restricted to one particular discontinuity, but performs a unified candidate comparison and greedy selection over all discontinuities that may be progressively exposed during the backward move. In other words, the bottleneck breakthrough rule may be viewed as the local implementation of the backward candidate traversal at a single discontinuity, whereas the backward candidate traversal constitutes its global extension over the entire backward-moving process.

Therefore, any critical discontinuity repair that restores the feasibility of the ideal backward direction must correspond to a discontinuity explicitly exposed at some backward candidate position; and all such candidate positions lie within the traversal range of the backward candidate set in Algorithm~\ref{algorithm3}. Hence, whenever there exists a path that restores the feasibility of the ideal backward direction through discontinuity repair, the backward candidate traversal in Algorithm~\ref{algorithm3} cannot miss the key repair steps required along that path. As a result, the improvement potential of the ideal backward direction in the decrease-flow region can be released through the joint action of the backward candidate traversal and the bottleneck breakthrough rule.

We next consider the increase-flow region. By \cref{Proposition4}, under the adjacent-swap rule, the local optima induced in the increase-flow region can only take two forms: the under-consumption type and the over-consumption type. An under-consumption local optimum essentially corresponds to an executable forward improvement whose potential has not yet been fully released; hence, it can be handled by the forward candidate traversal in Algorithm~\ref{algorithm1} together with the bottleneck breakthrough rule. An over-consumption local optimum, by contrast, arises because an excessive forward shift of some job introduces too much machine idle time, whereas the better structure should be realized by a combined forward adjustment of several jobs. Moreover, during the elimination of an over-consumption local optimum, the resulting structure may further induce new under-consumption local optima. However, by \cref{over under} in the appendix, regardless of whether such new under-consumption local optima are induced, all these subsequent local structures can be eliminated successively within the combined neighborhood induced by Algorithm~\ref{algorithm1} and Algorithm~\ref{algorithm3}. Consequently, the additional loss caused by the increase flow can be continuously compressed during the combined adjustment process. The technical derivation is given in Appendix~\ref{app:over under}.

Furthermore, by \cref{increase flow}, if the increase-flow region itself contains no independently feasible adjustment that strictly decreases the total waiting time, then any subsequent adjustment within it can at most reduce the additional loss caused by the increase flow, but cannot become an independent source of net improvement in the total waiting time. Therefore, in the combined repair process for a backward-type local optimum, the adjustments in the increase-flow region serve only to coordinate and compress the side effects; the actual restoration of strict improvement is still ensured by the fact that the critical discontinuity in the decrease-flow region of the ideal backward direction has been repaired.

Combining the above arguments, within the combined neighborhood induced by Algorithm~\ref{algorithm1} and Algorithm~\ref{algorithm3}, the critical discontinuities in the decrease-flow region can be repaired, while the additional loss in the increase-flow region can be further compressed. By the definition of a backward-type local optimum, it then follows that the ideal backward direction, which originally could not realize a strict improvement directly, can be transformed into a strict improvement under the above combined adjustments. Denote the resulting new schedule by $\tilde{\sigma}$. Then
\[
F(\tilde{\sigma})<F(\sigma^t).
\]
Therefore, whenever $\sigma^t$ is a backward-type local optimum, it can be broken by the combined neighborhood.
\end{proof}

The preceding two lemmas rule out the possibility that the algorithm may remain at either a forward-type local optimum or a backward-type local optimum. Combined with the complete classification in \cref{Lemma6}, it follows that the algorithm cannot terminate at any non-globally optimal solution. We thus obtain the following optimality criterion for termination.

\begin{theorem}
\label{Theorem3}
The algorithm cannot terminate at a non-globally optimal solution.
\end{theorem}

\begin{proof}
Suppose that the algorithm terminates at the sequence $\sigma^t$. By the definition of the algorithm, at this iteration all feasible forward moves, backward moves, and feasible adjustments induced by the bottleneck breakthrough rule and the adjacent interchange rule within the combined neighborhood have already been enumerated. Moreover, according to the maximum-improvement principle, the algorithm selects the feasible move that yields the largest decrease in the objective value. Since the algorithm has terminated, there exists no feasible combined adjustment within the currently enumerated combined neighborhood that can strictly decrease the total waiting time.

We prove by contradiction that $\sigma^t$ must then be globally optimal.

Assume that $\sigma^t$ is not globally optimal. Then there exists a feasible sequence $\sigma^*$ such that
\[
F(\sigma^*)<F(\sigma^t),
\]
where $F(\sigma)$ denotes the total waiting time of sequence $\sigma$.

By \cref{Lemma6}, any non-globally optimal solution must belong to at least one of the following two local structural types:
\begin{enumerate}
    \item a forward-type local optimum;
    \item a backward-type local optimum.
\end{enumerate}
Therefore, the terminal sequence $\sigma^t$ must also belong to at least one of these two types.

If $\sigma^t$ is a forward-type local optimum, then \cref{Lemma7} implies that there exists a feasible combined adjustment within the combined neighborhood enumerated by the algorithm that strictly decreases the total waiting time. This contradicts the fact that the algorithm has terminated at $\sigma^t$.

If $\sigma^t$ is a backward-type local optimum, then \cref{Lemma8} likewise implies that there exists a feasible combined adjustment within the combined neighborhood enumerated by the algorithm that strictly decreases the total waiting time. This again contradicts the fact that the algorithm has terminated at $\sigma^t$.

Both cases lead to a contradiction. Hence, the assumption is false, and $\sigma^t$ must be globally optimal. Therefore, the algorithm cannot terminate at a non-globally optimal solution.
\end{proof}

\cref{Theorem3} establishes the optimality guarantee of the terminal solution. The remaining question is whether the algorithm necessarily terminates after finitely many steps, and whether a sharper upper bound can be derived for the number of iterations of the outer \texttt{while} loop. In the following, we exploit the structural progression property of backward-type local optima to prove that the \texttt{while} loop in Algorithm~\ref{algorithm1} is executed at most $n$ times, thereby completing the proof of global convergence of the algorithm.

\begin{theorem}
\label{Theorem4}
The proposed algorithm terminates in finitely many steps and returns a global optimum.
\end{theorem}

\begin{proof}
By \cref{Theorem3}, once the algorithm terminates, it cannot stop at a non-global optimum. Hence, its output must be globally optimal. It therefore remains only to prove that the algorithm terminates in finitely many steps. More specifically, we show that the \texttt{while} loop in Algorithm~\ref{algorithm1} is executed at most $n$ times.

From the preceding optimality analysis and the algorithm design, once Algorithm~\ref{algorithm2} terminates, any remaining non-optimal local structure in the sequence can only be a backward-type local optimum. When backward moves and the adjacent-swap rule are applied to such a structure, two local configurations may arise: over-consumption and under-consumption. By \cref{over under}, the former can be resolved by Algorithm~\ref{algorithm3} itself, whereas the latter requires the joint action of the forward candidate traversal and the bottleneck breakthrough rule in Algorithm~\ref{algorithm1}, together with Algorithm~\ref{algorithm3}. This is precisely why the \texttt{while} loop is introduced in Algorithm~\ref{algorithm1}.

Let $b_t$ denote, at the beginning of the $t$-th iteration of the \texttt{while} loop, the position of the leftmost backward-type local optimum that has not yet been completely eliminated in the current sequence. Clearly,
\[
1\le b_t\le n.
\]
Furthermore, let $\mathcal{S}_t$ be the set of jobs involved in the under-consumption local optimum induced by the increase flow at position $b_t$, and define
\[
\Phi_t:=|\mathcal{S}_t|,
\qquad 0\le \Phi_t\le n.
\]

For a fixed position $b_t$, if the current local structure is an under-consumption local optimum induced by the increase flow, then one forward move together with the combined adjustment of Algorithm~\ref{algorithm3} strictly reduces the set of jobs involved in this local optimum type. In other words, after the adjustment, at least one job no longer participates in that local optimum structure, and therefore
\[
\Phi_{t+1}<\Phi_t.
\]
Thus, the participating job set of the leftmost under-consumption local optimum induced by the increase flow is irreversibly reduced. It follows that the local structure at position $b_t$ must be completely eliminated after finitely many such combined adjustments. Once the backward-type local optimum at this position is eliminated, the processing frontier cannot remain at the same position and must move strictly to the right.

Moreover, after Algorithm~\ref{algorithm3} eliminates the under-consumption local optimum induced by the increase flow in the current sequence, if the current sequence is still not globally optimal, then any newly generated non-optimal local structure can only appear in the $P_D$ region corresponding to the backward move. Since Algorithm~\ref{algorithm3} eliminates backward-type local optima, the local structure classification implies that this new structure can only be a forward-type local optimum. Because the region $P_D$ lies strictly to the right of the original increase-flow region $P_I$, the position of the leftmost unprocessed local optimum must also move strictly to the right.

Therefore, after each iteration of the \texttt{while} loop, the processing frontier advances strictly rightward. Since this frontier can move from the first job to at most the $n$-th job, the \texttt{while} loop in Algorithm~\ref{algorithm1} can be executed at most $n$ times. Hence, the entire algorithm terminates after finitely many steps.
\end{proof}

\textbf{Remark.}
\textbf{Remark.}
Both the analysis and the computational observations suggest that, for an instance with $n$ jobs, an under-consumption local optimum typically involves at least three jobs. As a consequence, the \texttt{while} loop is usually executed no more than $\left\lfloor \frac{n-1}{3} \right\rfloor + 2.$ The numerical results further indicate that this upper bound is conservative in practice.

In summary, the correctness proof of the proposed algorithm is now complete. \cref{Lemma6} establishes a complete classification of the local structures of non-globally optimal solutions; \cref{Lemma7} and \cref{Lemma8} show, respectively, that both forward-type and backward-type local optima can be broken through by the combined neighborhood; \cref{Theorem3} then implies that the algorithm cannot terminate at a non-globally optimal solution; and \cref{Theorem4} further proves that the algorithm necessarily reaches termination in finitely many steps. Therefore, the proposed algorithm solves the problem to global optimality in finitely many steps

\subsection{Complexity Analysis}

We now analyze the time complexity of the proposed algorithm for a set of $n$ jobs.

First, all jobs are sorted in nondecreasing order of release times. This preprocessing step can be completed in
\[
\mathcal{O}(n\log n)
\]
time.

Next, we construct the forward candidate set and the backward candidate set. In the worst case, for any position $i$, the number of forward candidate positions of job $\sigma_t(i)$ is at most $n-i$. Hence, the total number of forward candidates is at most
\[
\sum_{i=1}^{n-1}(n-i)=\frac{n(n-1)}{2}.
\]
Similarly, the total number of backward candidates is also at most $\frac{n(n-1)}{2}.$ Since the feasibility check for each candidate requires only constant time, the construction of both the forward and backward candidate sets has complexity $\mathcal{O}(n^2)$. Therefore, the total complexity of candidate-set construction is
\[
\mathcal{O}(n^2).
\]

For any individual candidate, its update procedure consists of two parts: the increase-flow part and the decrease-flow part. The increase-flow part is handled by the adjacent interchange rule, whereas the decrease-flow part is handled by the bottleneck breakthrough rule. Let $t$ denote the length of the relevant propagation region, where $t\le n$. In the worst case, the computational complexity of each of these two parts is
\[
\mathcal{O}(t^3)\subseteq \mathcal{O}(n^3).
\]
Hence, the evaluation and update of a single candidate requires
\[
\mathcal{O}(n^3)
\]
time.

Algorithms~\ref{algorithm2} and \ref{algorithm3} each perform one traversal of the forward candidate set and the backward candidate set, respectively. Their complexity is equal to the complexity of processing a single candidate multiplied by the size of the corresponding candidate set, namely,
\[
\mathcal{O}(n^3)\times \mathcal{O}(n^2)=\mathcal{O}(n^5).
\]
Since these two algorithms must exhaustively process the forward candidate set and the backward candidate set, and in the worst case each iteration resolves only the forward or backward candidate set of one job, as many as $n$ iterations may be required. Therefore, the overall complexity of Algorithms~\ref{algorithm2} and \ref{algorithm3} is
\[
\mathcal{O}(n^5)\times \mathcal{O}(n)=\mathcal{O}(n^6).
\]

In each iteration of Algorithm~\ref{algorithm1}, at most all forward candidates need to be traversed. Since Algorithms~\ref{algorithm2} and \ref{algorithm3} are nested within this outer loop, the complexity of one iteration of the main algorithm is at most
\[
\mathcal{O}(n^2)\times \mathcal{O}(n^6)=\mathcal{O}(n^8).
\]

Finally, consider the number of iterations of the main algorithm. In the worst case, the \texttt{while} loop is executed $n$ times. Then, by the preceding result, the total complexity of the main algorithm is
\[
\mathcal{O}(n)\times \mathcal{O}(n^8)=\mathcal{O}(n^9).
\]

In summary, the total time complexity of the proposed algorithm is
\[
\mathcal{O}(n^9).
\]

\section{Conclusion}
This paper studies the non-preemptive single-machine scheduling problem with heterogeneous release times and processing times, with the objective of minimizing the total waiting time. By modeling machine idle time as negative waiting time, we derive a unified waiting-time formula and reduce the original four-dimensional description of the problem to a two-dimensional structural representation. This reformulation enables a tractable structural analysis of the problem. Building on it, we introduce ideal directions induced by an optimal solution, establish a framework for analyzing the evolution from a non-optimal schedule to an optimal one, and show that improvable ideal directions admit no coordination relationship.

For each improvable ideal direction, we prove that queue discontinuity is the unique structural obstacle to its realization. We then give a complete characterization of the associated local optimum types, which reveals both how such local optima arise and how they can be eliminated. Based on these results, we develop an iterative structural-repair algorithm that repeatedly identifies and removes structural obstacles, thereby driving the current solution along ideal improvement paths toward a global optimum. We prove that every repair operator is polynomial-time implementable and yields a strict improvement in the objective value. Moreover, the number of jobs involved in structural violations decreases strictly during the iterative process, which ensures termination in a polynomial number of steps. The resulting overall worst-case time complexity is $\mathcal{O}(n^9)$.

Future work may focus on improving computational efficiency and extending the framework to broader scheduling settings.

\appendix
\section{Omitted Proofs in Section 2}

\subsection{ Proof of \cref{Lemma1}}\label{app:Lemma1}
\begin{proof}
The waiting time computation formula indicates that the queue leader is key to the calculation. When the waiting time of job $\sigma_t(j)$ drops by \(\Delta w_{\sigma_t(j)}^- = \Delta \) ($\Delta\geq0$), it is necessary to classify and discuss whether it is the queue leader. 

We label the neighboring jobs $\sigma_t(i)$ and $\sigma_t(i+1)$ as $j$ and $k$, respectively, for simplicity in the reasoning process, where $j < k$ indicating that job $j$ precedes job $k$ in processing. $w_j^2$ and $w_k^2$ are the waiting times of jobs $j$ and $k$, correspondingly, after the change. 

\textbf{Case 1: \( w_j > \Delta \).}

   In this case, job $j$ is neither the queue leader initially nor after the reduction. The queue leader remains unchanged, and the updated waiting time for job $j$ is:  
   \[   w_j^2 = w_j - \Delta > 0.   \]  
   
   The waiting time of job $k$ is updated as follows:  
   \[   w_k^2 = r_j + w_j^2 + p_j - r_k,   \]  
   since the original waiting time for job $k$ is:  
   \[   w_k = r_j + w_j + p_j - r_k,   \]  
   we can get the waiting time reduction for job $k$:  
   \begin{align*}
       \Delta w_k^- &= w_k - w_k^2 \\&=r_j + w_j + p_j - r_k-(r_j + w_j^2 + p_j - r_k)\\&= w_j - w_j^2 = \Delta.
   \end{align*}

\textbf{Case 2: \( 0 < w_j \leq \Delta \).} 

   In this scenario, we have $w_j^2 = w_j - \Delta \leq 0$, meaning that job $j$ becomes the new queue leader and starts processing immediately.  
   
   The new waiting time for job $k$ is:  
   \[
   w_k^2 = r_j^o + p_j - r_k,   \]        
   since its original waiting time is:  
   \[
   w_k = r_j + w_j + p_j - r_k,
   \]  
   the waiting time reduction of job $k$ is:  
   \begin{align*}
       \Delta w_k^- &= w_k - w_k^2 \\&=r_j + w_j + p_j - r_k-(r_j^o + p_j - r_k)\\&= w_j.
   \end{align*}
   (The superscript $o$ is only an indicator of the queue leader and does not change the values of the release time.)

\textbf{Case 3: \( w_j \leq 0 \).}

   In this case, job $j$ is already the queue leader initially, and a further reduction in its waiting time does not affect subsequent jobs. The updated waiting time for job $k$ is:  
   \[
   w_k^2 = r_j^o + p_j - r_k,
   \]  
     and the original waiting time for job $k$ is:  
   \[
   w_k = r_j^o + p_j - r_k,
   \]  
    the waiting time reduction for job $k$ can be get:  
   \begin{align*}
       \Delta w_k^- &= w_k - w_k^2\\&=r_j^o + p_j - r_k-(r_j^o + p_j - r_k)\\ &= 0.
   \end{align*}
   
\textbf{Conclusion:}  
We can conclude from these cases that: If $0<w_j \leq\Delta $, the decrease in waiting time of job $k$ is $w_j$; if $\Delta< w_j $, job $k$'s waiting time reduced $\Delta$. It always takes the smaller of the two.

\begin{align*}
    \Delta w_k^- &=\min(\Delta,w_j)\\&= \min(\Delta w_j^-, w_j).
\end{align*}

If $w_j \leq 0$, then:

\[
\Delta w_k^- = 0.
\]
We can express it concisely with a mathematical formula as follows: 
\begin{align*}
    \Delta w_k^- =\begin{cases}
        0,&\textit{if }w_{\sigma_t(k)} \leq 0,\\\min(\Delta w_{j}^-, w_{j}), & \textit{if } w_{\sigma_t(k)} > 0.
    \end{cases}
\end{align*}

This conclusion can be generalized to the waiting time changes of all subsequent jobs, finally obtaining a standard recursive formula, when $w_{\sigma_t(j)}$ is decreased:
\begin{align*}
    \Delta w_{\sigma_t(k+1)}^- =
  \begin{cases}
0 ,& \textit{if } w_{\sigma_t(k)} \leq 0,\\\min(\Delta w_{\sigma_t(k)}^-, w_{\sigma_t(k)}), & \textit{if } w_{\sigma_t(k)} > 0.
  \end{cases}
\end{align*}

where the initial condition is $w_{\sigma_t(j)} = \Delta$ and $k\in [j,n]$.  
\end{proof}
 
\subsection{ Proof of \cref{Lemma2}}\label{app:Lemma2}
\begin{proof}
    Similar to the demonstration in Lemma 1, the adjacent jobs $\sigma_t(i)$ and $\sigma_t(i+1)$ are denoted as $j$ and $k$, respectively, with $j < k$ meaning that job $j$ precedes job $k$ in processing. $w_j^2$ and $w_k^2$ are the revised waiting times for jobs $j$ and $k$, correspondingly, following the adjustment. 
   
\textbf{Case 1: \( w_j + \Delta \leq 0 \). } 
   Job $j$ is the queue leader initially. Even after its waiting time increases by $\Delta$ $w_j^2 = w_j + \Delta \leq 0,$ it remains the queue leader.
  
   Since job $j$ still leads the queue, the waiting time of job $k$ remains unchanged:
   \[w_k^2 = r_j^o + p_j - r_k,\]
   \[\quad w_k = r_j^o + p_j - r_k.\]
   Hence, the increase in waiting time for job $k$ is:
   \[\Delta w_k^+ = w_k^2 - w_k = 0.\]

    \textbf{Case 2:  }\( w_j \leq 0 \) and \( 0 < w_j + \Delta \).
   Job $j$ is the queue leader at the beginning. It is no longer the head when its waiting time is increased $\Delta$, and becomes a member of the continuous queue.
   \[ w_j^2 = w_j + \Delta > 0.\]
   The waiting time of job $k$ updates as:
   \[w_k^2 = r_j + w_j^2 + p_j - r_k,\]
   given the original waiting time:
   \[w_k = r_j^o + p_j - r_k,\]
   the waiting time increase for job $k$ is:
   \begin{align*}
       \Delta w_k^+ &= w_k^2 - w_k \\&= r_j + w_j^2 + p_j - r_k-(r_j^o + p_j - r_k) \\&=w_j^2\\&= w_j + \Delta.
   \end{align*}

\textbf{Case 3: }\( w_j > 0 \).
   Job $j$ is not the queue leader, increasing its waiting time does not change its queue status$ w_j^2 = w_j + \Delta > 0.$
     
   The updated waiting time of job $k$ is:
   \[
   w_k^2 = r_j + w_j^2 + p_j - r_k.
   \]
   Originally:
   \[
   w_k = r_j + w_j + p_j - r_k.
   \]
   Last, the waiting time increase for job $k$ is:
   \begin{align*}
       \Delta w_k^+ &= w_k^2 - w_k \\&=r_j + w_j^2 + p_j - r_k-(r_j + w_j + p_j - r_k)\\&= w_j^2 - w_j \\&= \Delta.
   \end{align*}

\textbf{Conclusion:}  
From these cases, we conclude: If $w_j \leq 0$ and $w_j+\Delta\leq0$, then we get $\Delta w_k^+=0$. If $w_j \leq 0$ and $w_j+\Delta>0$, it will be $\Delta w_k^+=w_j+\Delta$. $\Delta w_k^+$ is always the bigger one of $0$ and $w_j+\Delta$.
  \begin{align*}
      \Delta w_k^+ &= \max(0, w_j + \Delta) \\&= \max(0, w_j + \Delta w_j^+).
  \end{align*}
  
When $w_j > 0$, it seems that only one result exists:
  \[\Delta w_k^+ = \Delta = \Delta w_j^+.\]
  Finally, we get the concisely mathematical formula:
  \begin{align*}
      \Delta w_k^+=\begin{cases}
          \Delta w_j^+,&\textit{if } w_k>0,\\\max(0,\Delta w_j^+ +w_j),&\textit{if }w_k\leq0.
      \end{cases}
  \end{align*}
This conclusion can be further broadened to the waiting time changes of all subsequent jobs, finally, a general recursive formula is obtained, when $w_{\sigma_t(j)}$ is increased: 

\begin{equation*}
\Delta w_{\sigma_t(k+1)}^+ = \begin{cases}
\Delta w_{\sigma_t(k)}^+, & \textit{if } w_{\sigma_t(k)} > 0, \\
\max(0, \Delta w_{\sigma_t(k)}^+ + w_{\sigma_t(k)}), & \textit{if } w_{\sigma_t(k)} \leq 0.
\end{cases}
\end{equation*}

where the initial condition is $\Delta w_{\sigma_t(j)}^+ = \Delta$ and $k\in [j,n]$.  
\end{proof}

\subsection{ Proof of \cref{Proposition1}}\label{app:Proposition1}
\begin{proof}
The proposition’s first inequality and second inequality outline the relationships that waiting time and processing time are required to meet for the first-come, first-served rule and the first-come, last-served rule, respectively. 

For two adjacent jobs $\sigma_t(i)$ and $\sigma_t(i+1)$ in the queue $\sigma_t$.  When $r_{\sigma_{t} (i)}\leq r_{\sigma_{t} (i+1)}$, it reflects the first-come, first-served rule. In that case, the maximum waiting time for $\sigma_t(i+1)$ is the waiting time of  $\sigma_t(i)$ plus its processing time.  A greater waiting time could only occur if it arrives prior to $\sigma_t(i)$ but is serviced after $\sigma_t(i)$, which corresponds to a first-come, last-served rule.  

To simplify the derivation, let  $i = \sigma_t(i)$ and $i+1 = \sigma_t(i+1)$, with the service sequence always having job $i$ precede $i+1$. The release time of the previous queue leader and the cumulative processing time up to job $i$ are $r_0^o$  and $p_0^o$, respectively. As defined earlier, the waiting time for job $i$ is\[w_i=r_0^o + p_0^o - r_i.\]The waiting time calculation for job $i+1$ needs to account for changes of the queue leader: 
\begin{itemize}
    \item When $w_i>0, w_{i+1}=r_0^o + p_0^o + p_i - r_{i+1}.$
    \item When $w_i\leq 0,w_{i+1}=r_i^o+p_i - r_{i+1}.$
\end{itemize}
The two inequalities of the proposition that require proof are as follows: 
\begin{itemize}
    \item When $r_i\leq r_{i+1}, \max(0,w_{i})+p_{i}\geq\max(0,w_{i+1}).$
    \item When $r_i>r_{i+1}, \max(0,w_{i})+p_{i}<\max(0,w_{i+1}).$
\end{itemize}
As can be observed, the positivity or negativity of $w_i$ and $w_{i+1}$ is closely related to the proof of the inequalities. We will classify cases according to them.

Firstly, we prove the inequality for the first part under the first-come, first-served principle($r_{i+1} \geq r_i$): 
     
    \textbf{Case 1:} $w_i > 0, w_{i+1} > 0$, both job $i$ and $i+1$ are in the same continuous queue. We have\[\max(0, w_i) + p_i = w_i + p_i = r_0^o + p_0^o - r_i + p_i,\]
    \[\max(0, w_{i+1}) = w_{i+1} = r_0^o + p_0^o + p_i - r_{i+1},\]
    \[\because \text{the first-come, first-served } r_{i+1} \geq r_i, \therefore \max(0, w_i) + p_i\geq \max(0,w_(i+1)).\]
    
    The inequality holds.
    
    \textbf{Case 2:} $w_i > 0, w_{i+1} \leq 0$, job $i$ remains in the continuous queue, while job $i+1$ is the queue leader. At this point, we get   
    \begin{align*}
        &\max(0, w_i) + p_i \\&= w_i + p_i \\
            &= r_0^o + p_0^o - r_i + p_i \geq 0 = \max(0, w_{i+1})
    \end{align*}

    The inequality can be confirmed.  
    
    \textbf{Case 3:} $w_i \leq 0, w_{i+1} > 0$, job $i$ is the queue leader, while job $i+1$ belongs to the continuous queue. So we have 
    \[\max(0, w_i) + p_i = p_i,\]
    \[\max(0, w_{i+1}) = w_{i+1} = r_i^o + p_i - r_{i+1},\]
    \[\because r_{i+1}\geq r_i^o,\therefore p_i\geq r_i^o + p_i - r_{i+1},\]\[\therefore \max(0, w_i) + p_i\geq\max(0,w_{i+1}).\]  
    
    The inequality is shown to hold.  
    
    \textbf{Case 4:} $w_i \leq 0, w_{i+1} \leq 0$, both jobs are the queue leader.  At this time, we get
    \[\max(0, w_i) + p_i = p_i \geq 0 = \max(0, w_{i+1}).\]
    
    Obviously, the inequality continues to hold.

    As discussed above, the first inequality of the proposition is verified. We will explore the second portion of the proposition under the first-served, last-come rule. Under this rule, when the release time $r_i$ of job $i$ is later than that of job $i+1$, and the service order of job $i$ precedes that of job $i+1$, the waiting time $w_{i+1}$ of job $i+1$ is bound to be greater than 0. Therefore, only two cases need to be discussed.
    
    \textbf{Case 1:} $w_i > 0, w_{i+1} > 0$, both job $i$ and $i+1$ are in the same continuous queue. We have\[\max(0, w_i) + p_i = w_i + p_i = r_0^o + p_0^o - r_i + p_i,\]
    \[\max(0, w_{i+1}) = w_{i+1} = r_0^o + p_0^o + p_i - r_{i+1},\]
    \[\because \text{the first-come, last-served } r_{i+1} < r_i, \therefore \max(0, w_i) + p_i< \max(0,w_(i+1)).\]
    
    The inequality holds.
    
    \textbf{Case 2:} $w_i \leq 0, w_{i+1} > 0$, job $i$ is the queue leader, while job $i+1$ belongs to the continuous queue. So we have 
    \[\max(0, w_i) + p_i = p_i,\]
    \[\max(0, w_{i+1}) = w_{i+1} = r_i^o + p_i - r_{i+1},\]
    \[\because \text{the first-come, last-served }r_{i+1}< r_i^o,\therefore p_i< r_i^o + p_i - r_{i+1},\]\[\therefore \max(0, w_i) + p_i<\max(0,w_{i+1}).\]  
    
    The inequality is proved.  
    
    So far, we have validated the relationships between the waiting time and processing time of adjacent jobs required under the first-come, first-served, and first-come, last-served principles, respectively, and Proposition 1 is proven. \hfill $\square$
\end{proof}

\subsection{ Proof of \cref{Proposition2}}\label{app:Proposition2}
\begin{proof}
    Consider an arbitrary schedule $\sigma_t$ and any three consecutive jobs $\sigma_t(i)$, $\sigma_t(i+1)$, and $\sigma_t(i+2)$, with $i \in {1,\dots,n-2}$. Without loss of generality, denote these three jobs by $4$, $5$, and $6$, respectively, that is, let $\sigma_t(i)=4$, $\sigma_t(i+1)=5$, and $\sigma_t(i+2)=6$. This relabeling is adopted purely for notational simplicity and has no effect on the validity of the analysis. We focus on the situation where the waiting time of job $4$ is reduced by $\Delta w_4^-=\Delta$, with $\Delta \ge 0$.

\textbf{Situation 1: $w_4 \leq 0$}. 
 
  Job 4 is the queue leader at the beginning, and with its waiting time decreased $\Delta$, $w_4^2 = w_4 - \Delta \leq 0$, leading to the impact on the total waiting time.
  \begin{align*}
       \Delta w &= \max(0, w_4) - \max(0, w_4^2)  \\&= 0 - 0 = 0.
   \end{align*}
  We can get the values of  $ \Delta w_5^-$ according to the recurrence relation$ \Delta w_5^- = 0$, so the changes of the total waiting time equal with $ \Delta w_5^-$, $\Delta w=\Delta w_5^-$.
   
   \textbf{Situation 2:}$w_4 > 0$ and $w_4 - \Delta > 0$.   
   Job 4 is the member of the continuous queue initially and remains so after decreasing its waiting time by $\Delta$. The updated waiting time for job 4 is $w_4^2 = w_4 - \Delta > 0.$
   
   The change in the total waiting time is:
      \begin{align*}
          \Delta w &= \max(0, w_4) - \max(0, w_4^2)\\
      & = w_4 - (w_4 - \Delta) \\&= \Delta.
      \end{align*}
    $\Delta w_5^- = \min(\Delta, w_4) = \Delta$ according to the recurrence formulation. Obviously, $\Delta w=\Delta w_5^- $.

\textbf{Situation 3:} $w_4 > 0$ and $w_4 - \Delta \leq 0$.
  
  In this situation, job 4 transitions to become the new queue leader from the ones of the continuous queue, with $w_4^2 = w_4 - \Delta \leq 0$. $\Delta w_5^-$ is calculated as: \[\Delta w_5^- = \min(\Delta, w_4) = w_4.\]
   The total waiting time change is:
   \begin{align*}
       \Delta w &= \max(0, w_4) - \max(0, w_4^2) \\&= w_4 - 0  = w_4.
   \end{align*}
   It seems clear that $\Delta w=\Delta w_5^-$.

\textbf{Conclusion:} The reduction in job 4’s waiting time, lowering the total waiting time by $\Delta w$, precisely equals its backward transmission value $\Delta w_5^-$, thereby proving the proposition. This conclusion holds for all subsequent jobs, following the same proof process.

Similarly, When the waiting time of job 4 increases $\Delta w_4^+ = \Delta$ ($\Delta \geq 0$), using the same method of analysis can verify that the total waiting time increase, $\Delta w$, is the same as $\Delta w_5^+$.\hfill $\square$
\end{proof}

\subsection{ Proof of \cref{Theorem1}}\label{app:Theorem1}
\begin{proof}
We prove Theorem 1 in two steps. First, we establish the claimed propagation pattern by examining a one-position forward move. We then consider a two-position forward move and show that further forward moves only add processing-time and waiting-time information, without creating any new structural components. Therefore, the pattern admits a direct extension to arbitrarily long forward moves.
    
    Consider an arbitrary schedule $\sigma_t$ and any four consecutive jobs $\sigma_t(i)$, $\sigma_t(i+1)$, $\sigma_t(i+2)$, and $\sigma_t(i+3)$, where $i \in {1,\dots,n-3}$. We begin with the simplest case, in which job $\sigma_t(i)$ is moved to the position immediately after $\sigma_t(i+1)$, as illustrated in Figure~\ref{forward1}. For ease of exposition and without loss of generality, let $\sigma_t(i)=4$, $\sigma_t(i+1)=5$, $\sigma_t(i+2)=6$, and $\sigma_t(i+3)=7$. This relabeling is introduced solely for notational convenience and does not affect the subsequent derivation.
    
    Before starting the derivation, we need to define certain notations and provide proof of the impact of removing job 4 on the waiting time of job 5. 
    
    First, we prove that the removal of job 4 will decrease the waiting time of job 5 by $\Delta w_5=p_4-\min(0,w_4)$, and the proof processes as follows: 
    
    Suppose $r_0^o$ and $p_0^o$ are the release time of the previous queue leader and the cumulative processing time up to job 4, respectively.
    \begin{itemize}
        \item When $w_4>0$, the original waiting time of job 5 is $w_5=r_0^o+p_0^o+p_4-r_5,$ after we remove job 4, the waiting time of job 5 is updated as $w_5^2=r_0^o+p_0^o-r_5$, the differ between the two is $\Delta w_5=w_5-w_5^2=p_4$.
        \item  When $w_4\leq0$, the waiting time of job 5 starts with $w_5=r_0^o+p_0^o-w_4+p_4-r_5,$ after the removal of job 4, we get the new waiting time of job 5 $w_5^2=r_0^o+p_0^o-r_5$, the differ of the two is $\Delta w_5=w_5-w_5^2=p_4-w_4$.
    \end{itemize}
    
    We combine the two cases to reach the conclusion:$\Delta w_5=p_4-\min(0,w_4)$.
    
    Second, we record the effect of job $i$ on the total waiting time as $\Delta_i$. For example, if $w_5>0,w_5^2 = w_5 - \Delta w_5 \leq 0$, the total waiting time is reduced $\Delta_5 = -w_5$, instead of $\Delta w_5$. Likewise, when $w_5>0,w_5^2 = w_5 - \Delta w_5 > 0$, the total waiting time is decreased $\Delta_5 = -\Delta w_5$. So, the formula of $\Delta_5$ is that $\Delta_5=\max(0,w_5^2)-\max(0,w_5)$. Actually, $\Delta_i=-\Delta w_{i+1}$ as Proposition 2 suggests. We use $\Delta_5 = -\Delta w_6^-$ to calculate the impact of jobs who is not moved on the total waiting time, and $\Delta_5 = \max(0, w_5^2) - \max(0, w_5)$ to compute the effect of the moved job on the total waiting time.
    
    \begin{figure}
        \centering
        \includegraphics[width=0.8\linewidth]{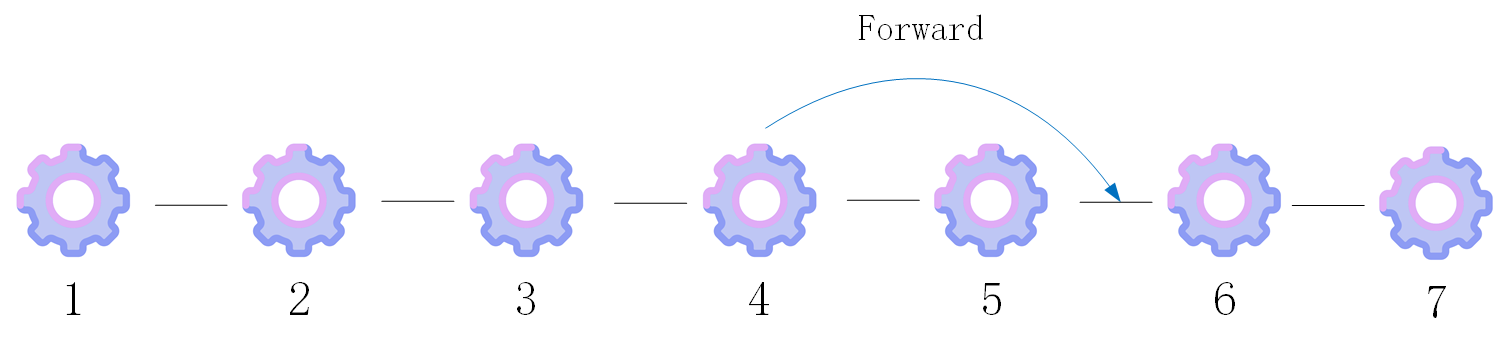}
        \caption{Job 4 moves to the position between jobs 5 and 6.}
        \label{forward1}
    \end{figure}
    
    Now, let's start the proven process. Beginning with the unique role of the queue leader, we list all possible situations and discuss the details of each situation.
    \begin{itemize}
        \item \textbf{Case 1:} $w_4\leq0,w_5>p_4-w_4$, under those conditions, job 4 is the queue leader, and once job 4 is moved forward, the waiting time of job 5 is positive, implying job 5 is one of the members of the continuous queue.
        \item \textbf{Case 2:} $w_4\leq0,0<w_5\leq p_4-w_4$, in this case, job 4 is the queue leader and job 5 is part of the continuous queue at the beginning. After job 4 moves forward, job 5 becomes the new queue leader.
        \item \textbf{Case 3:} $w_4\leq0,w_5\leq0$, both jobs 4 and 5 are the queue leader, after job 4 moves forward, job 5 still be the queue leader.
        \item \textbf{Case 4:} $w_4>0,w_5>p_4$, here, both jobs 4 and 5 belong to the continuous queue, and once job 4 moves, job 5's state is unchanged.
        \item \textbf{Case 5:} $w_4>0,0<w_5\leq p_4$, at this point, jobs 4 and 5 are part of the continuous queue. Job 5 is the new queue leader after job 4's relocation.
        \item \textbf{Case 6:} $w_4>0,w_5\leq0$, now, job 4 belongs to the continuous queue and job 5 is the queue leader. After job 4 moves forward, job 5 still be the queue leader.
    \end{itemize}

    We have listed the six cases above. Now, we will use mathematical formulas to deduce the impact of job 4’s movement on the queue and its mathematical laws. 
    
    \textbf{Case 1:} $w_4\leq0,w_5>p_4-w_4$
    
    Under those conditions, job 4 is the queue leader, and once job 4 is moved forward, the waiting time of job 5 is positive, implying job 5 is one of the members of the continuous queue. So we have:    \begin{align*}
        w_5^2 &= w_5 + \Delta w_5\\&=w_5-(p_4-w_4)>0,
    \end{align*} so job 5 remains in the continuous queue and the impact on the total waiting time is \begin{align*}
       \Delta_5&=-\Delta w_6^-=-(p_4-w_4).
    \end{align*}
        The new waiting time of job 4 is    \begin{align*}
    w_4^2 &= r_5 + w_5^2 + p_5 - r_4\\&=r_5+ w_5 - p_4 + w_4 + p_5   - r_4\\&= w_4 + p_5.
    \end{align*}
    (1) If $w_4 + p_5 \leq 0$, meaning job 4 still be the queue leader after its relocation:
    
    then\begin{align*}
        \Delta_4&=\max(0,w_4^2)-\max(0,w_4)\\&=0.
    \end{align*}   
   The waiting time change for job 6 is   \begin{align*}
       \Delta w_6 &= r_4^o + p_4 - r_6 - w_6\\&=r_4^o + p_4 - r_6 - (r_5 + w_5 + p_5 - r_6)\\&= r_4^o - r_5 + p_4 - w_5 - p_5\\&= -p_5.
   \end{align*}
   Thus, the total waiting time changes as follows:   \begin{align*}
       \Delta w &= \Delta_4 + \Delta_5 + \pi_5(-p_5)\\&=p_5 - (p_4 - w_4) - p_5 + \pi_6(-p_5).
   \end{align*}
    (2) If $w_4 + p_5 > 0$, job 4 remains in the continuous queue:
    
    The waiting time of job 4 is increased by $w_4+p_5$, so $\Delta_4=w_4+p_5$. The change of job 6's waiting time is    \begin{align*}
        \Delta w_6 &= r_4 + w_4^2 + p_4 - r_6 - w_6\\&= r_4 + w_4 + p_5 + p_4 - r_6 - (r_5 + w_5 + p_5 - r_6)\\&= r_4 + w_4 + p_4 - w_5 - r_5\\&=w_4.
    \end{align*}
    Finally, the total waiting time change is     \begin{align*}
        \Delta w &= \Delta_4 + \Delta_5 + \pi_5(w_4)\\&= w_4 + p_5 - (p_4 - w_4) + \pi_6(w_4)\\&= w_4 + p_5 - p_4 + w_4 + \pi_6(w_4).
    \end{align*}

    \textbf{Case 2:} $w_4 \leq 0 , 0 < w_5 \leq p_4 - w_4$
    
    In this case, job 4 is the queue leader and job 5 is part of the continuous queue at the beginning. After job 4 moves forward, job 5 becomes the new queue leader.
    
    At this situation,    \begin{align*}
        w_5^2 &= w_5 - (p_4 - w_4) \leq 0,\\\Delta_5&=-\Delta w_6^-=-w_5,
    \end{align*} job 5 becomes the new queue leader and decreases the total waiting time by $\Delta_5=-w_5$.
    The updated waiting time of job 4 is:    \begin{align*}
        w_4^2 &= r_5^o + p_5 - r_4\\&=r_5^o + p_5 - r_4 - p_4 + p_4\\&=p_4 + p_5 - w_5.
    \end{align*}
    (1) If \( p_4 + p_5 - w_5 \leq 0 \), it means job 4 still be the queue leader after the movement of it:
    
    For the reason that $w_4\leq0$ and $w_4^2\leq0$, so job 4 does not affect the total waiting time $\Delta_4=0$.
    The adjustment of job 6's waiting time is     \begin{align*}
        \Delta w_6 &= r_4^o + p_4 - r_6 - w_6\\&=r_4^o + p_4 - r_6 - (r_5 + w_5 + p_5 - r_6)\\&= w_5 - w_5 - p_5\\&= -p_5.
    \end{align*}
    Finally, the total waiting time change is    \begin{align*}
        \Delta w &= \Delta_4 + \Delta_5 + \pi_5(-p_5)\\&= p_5 - w_5 - p_5 + \pi_6(-p_5).
    \end{align*}
    (2) If $p_4 + p_5 - w_5 > 0$, job 4 is one of the members of the continuous queue:
    
    From the condition $w_4\leq0$ and $w_4^2=p_4 + p_5 - w_5$, we can get \begin{align*}
        \Delta_4&=\max(0,w_4^2)-\max(0,w_4)\\&=p_4 + p_5 - w_5.
    \end{align*}
    The waiting time adjustment for job 6 is    \begin{align*}
        \Delta w_6 &= r_4 + w_4^2 + p_4 - r_6 - w_6\\&=r_4 + p_4 + p_5 - w_5 + p_4 - r_6 - (r_5 + w_5 + p_5 - r_6)\\&=r_4 + p_4 - 2w_5 + p_4 - r_5\\&=p_4 - w_5.
    \end{align*}
    Consequently, the total waiting time changes as follows:    \begin{align*}
        \Delta w &= \Delta_4 + \Delta_5 + \pi_6(p_4 - w_5)\\&=p_4 + p_5 - w_5 - w_5 + \pi_6(p_4 - w_5)\\&= p_5 - w_5 + p_4 - w_5 + \pi_6(p_4 - w_5).
    \end{align*}
    
    \textbf{Case 3:}  $w_4 \leq 0 ,w_5 \leq 0$
    
    Both jobs 4 and 5 are the queue leader, after job 4 moves forward, job 5 still be the queue leader.
    
    At this stage, we can get the new waiting time of job 5 as follows:    \begin{align*}
        w_5^2 &= w_5 - (p_4 - w_4) \leq 0, \\ \Delta_5&=-\Delta w_6^-= 0,
    \end{align*} both the waiting time of job 5 is non-positive, job 5 still be the queue leader, and the impact on the total waiting time is $\Delta_5 = 0$.
    The updated waiting time of job 4 is:    \begin{align*}
        w_4^2 &= r_5^o + p_5 - r_4\\&=p_4 + p_5 - w_5.
    \end{align*}
    Since $w_5\leq0$, it follows that $w_4^2>0$. Besides, $w_4\leq0$, so we can get that \begin{align*}
        \Delta_4&=\max(0,w_4^2)-\max(0,w_4)\\&=p_4 + p_5 - w_5.
    \end{align*} 
    The new waiting time of job 6 is calculated as follows:    \begin{align*}
        \Delta w_6 &= r_4 + w_4^2 + p_4 - r_6 - w_6\\&=r_4 + p_4 + p_5 - w_5 + p_4 - r_6 - (r_5 + p_5 - r_6)\\&=r_4 + p_4 - w_5 + p_4 - r_5\\&=p_4.
    \end{align*}
    In the end, the total waiting time change is:    \begin{align*}
        \Delta w &= \Delta_4 + \Delta_5 + \pi_6(p_4)\\&=p_5 - w_5 + p_4 + \pi_6(p_4).
    \end{align*}
    
    \textbf{Case 4:} $w_4 > 0,w_5 > p_4$
    
    Here, both jobs 4 and 5 belong to the continuous queue, and once job 4 moves, job 5's state is unchanged.
    
    Under this condition,\begin{align*}
        w_5^2 &= w_5 - \Delta w_5 \\&= w_5 - p_4 > 0,
    \end{align*} job 5 is the members of the continuous queue, and the impact on the total waiting time is \begin{align*}
          \Delta_5&=-\Delta w_6^-= -p_4.
    \end{align*}
    The updated waiting time of job 4 is\begin{align*}
        w_4^2 &= r_5 + w_5^2 + p_5 - r_4\\&=r_5 - p_4 + w_5 + p_5 - r_4\\&=p_5 + w_4.
    \end{align*} 
    Because $w_4>0$, it is obvious that $w_4^2>0$, job 4 is one of the members of the continuous queue, and its contribution to the total waiting time is \begin{align*}
        \Delta_4&=\max(0,w_4^2)-\max(0,w_4)\\&=p_5
    \end{align*}
    For job 6, the waiting time change is \begin{align*}
        \Delta w_6 &= r_4 + w_4^2 + p_4 - r_6 - w_6\\&=r_4 + p_5 + w_4 + p_4 - r_6 - (r_5 + w_5 + p_5 - r_6)\\&= r_4 + w_4 + p_4 - r_5 - w_5\\&= 0.
    \end{align*}
    Consequently, the total waiting time change is \begin{align*}
        \Delta w &= \Delta_4 + \Delta_5 + \pi_6(0)\\&= p_5 - p_4 + \pi_6(0).
    \end{align*} 
    
    \textbf{Case 5:} $w_4 > 0, 0 < w_5 \leq p_4$
    
    At this point, jobs 4 and 5 are part of the continuous queue. Job 5 is the new queue leader after job 4's relocation.
    
    In this case,  \begin{align*}
        w_5^2&=w_5-\Delta w_5\\&=w_5 - p_4\leq0,\\\Delta_5&=-\Delta w_6^-=-w_5.
    \end{align*} Since job 5 becomes the new queue leader, the impact on the total waiting time is $\Delta_5=-w_5$.
    The new waiting time of job 4, after its relocation is \begin{align*}
        w_4^2 &= r_5^o + p_5 - r_4\\&=r_5^o + p_5 - r_4 - p_4 + p_4 + w_4 - w_4\\&=p_4 + p_5 - w_5 + w_4.
    \end{align*} Because $w_5\leq p_4,w_4>0$, it has that $w_4^2>0$ and job 4 in the continuous queue. The contribution of job 4 to the total waiting time is as follows:\begin{align*}
        \Delta_4 &=\max(0, w_4^2) -\max(0, w_4) \\&= p_4 + p_5 - w_5.
    \end{align*} 
    The adjustment of job 6 can be calculated as follows:\begin{align*}
        \Delta w_6 &= r_5^o + p_5 + p_4 - r_6 - w_6\\&=r_5^o + p_5 + p_4 - r_6 - (r_5^o + w_5 + p_5 - r_6)\\&=p_4 - w_5.
    \end{align*}
    Finally, the total waiting time change is \begin{align*}
        \Delta w &= \Delta_4 + \Delta_5 + \pi_6 (p_4 - w_5)\\&=p_5 - w_5 + p_4 - w_5 + \pi_6 (p_4 - w_5).
    \end{align*}
    
    \textbf{Case 6:}  $w_4 > 0 ,  w_5 \leq 0$
    
    Now, job 4 belongs to the continuous queue and job 5 is the queue leader. After job 4 moves forward, job 5 still be the queue leader.
    
    In this case,\begin{align*}
        w_5^2 &=w_5-\Delta w_5\\&= w_5 - p_4 \leq 0,\\\Delta_5&=-\Delta w_6^-= 0,
    \end{align*} job 5 still be the queue leader, and has no effect on the total waiting time.
    The new waiting time of job 4 is calculated as follows:\begin{align*}
        w_4^2 &= r_5^o + p_5 - r_4\\&=r_5^o + p_5 - r_4 + p_4 - p_4 + w_4 - w_4\\&=p_4 + p_5 + w_4 - w_5.
    \end{align*} Because $w_5 \leq 0,w_4>0$, $w_4^2>0$. The impact on the total waiting time is \begin{align*}
        \Delta_4 &= \max(0, w_4^2) - \max(0, w_4) \\&= w_4^2 - w_4 \\&= p_4 + p_5 - w_5.
    \end{align*}
    For job 6, the waiting time adjustment is \begin{align*}
        \Delta w_6 &= r_5^o + p_5 + p_4 - r_6 - w_6\\&=r_5^o + p_5 + p_4 - r_6 - (r_5^o + p_5 - r_6)\\&=p_4.
    \end{align*}
    Finally, the change of the total waiting time is \begin{align*}
        \Delta w &= p_4 - w_5 + p_5 + \pi_6 (p_4)\\&=p_5 - w_5 + p_4 + \pi_6 (p_4).
    \end{align*}
    
    So far, we have completed the discussion of all cases. In the next step, we will summarize the underlying laws. First, we gather the outcomes and observe their patterns. 
    
    \textbf{Case 1:} $w_4 \leq 0, w_5 > p_4 - w_4$
    \begin{align*}
        w_4^2 = w_4 + p_5.
    \end{align*}

    (1) If \( w_4 + p_5 \leq 0 \), then \[\Delta w = p_5 - (p_4 - w_4) - p_5 + \pi_6 (-p_5).\]
    
    (2) If \( w_4 + p_5 > 0 \), then \[\Delta w = w_4 + p_5 - p_4 + w_4 + \pi_6 (w_4).\]
    
    \textbf{Case 2:} $w_4 \leq 0, 0 < w_5 \leq p_4 - w_4$
    \begin{align*}
        w_4^2 = p_4 + p_5 - w_5.
    \end{align*}
    
    (1) If $p_4 + p_5 - w_5 \leq 0$, then:\begin{align*}
        \Delta w = p_5 - w_5 - p_5 + \pi_6 (-p_5).
    \end{align*}
    
    (2) If \( p_4 + p_5 - w_5 > 0 \), then:\begin{align*}
        \Delta w = p_5 - w_5 + p_4 - w_5 + \pi_6 (p_4 - w_5).
    \end{align*}
    
    \textbf{Case 3:} $w_4 \leq 0, w_5 \leq 0$ \begin{align*}
        w_4^2 &= p_4 + p_5 - w_5,\\\Delta w &= p_5 - w_5 + p_4 + \pi_6 (p_4).
    \end{align*}
    
    \textbf{Case 4:}$w_4 > 0, w_5 > p_4$ \begin{align*}
        w_4^2 &= p_5 + w_4,\\\Delta w &= p_5 - p_4 + \pi_6 (0).
    \end{align*}
    
    \textbf{Case 5:} $w_4 > 0, 0 < w_5 \leq p_4$\begin{align*}
        w_4^2 &= p_4 + p_5 - w_5 + w_4,\\\Delta w &= p_5 - w_5 + p_4 - w_5 + \pi_6 (p_4 - w_5).
    \end{align*}
    
    \textbf{Case 6:} $w_4 > 0, w_5 \leq 0$\begin{align*}
        w_4^2 &= p_4 + p_5 + w_4 - w_5,\\\Delta w &= p_5 - w_5 + p_4 + \pi_6 (p_4).
    \end{align*}

    After organizing the results, we can observe that the terms in $\pi(.)$ always have corresponding terms in the first half of the formula, and these terms are closely related to the conditions that segment scenarios. Besides, it is obvious that the remaining terms also connect to the splitting conditions. We call the terms in $\pi(.)$ as the latter terms $I_2$ and the remaining terms as the former terms $I_1$. Based on this observations, we summarize $I_1$ and $I_2$ separately as follows.

    We start by summarizing $I_2$: 
    
    We combine conditions (1) and (2) for $I_2$ as follows:
    
    \textbf{Case 1:} $w_4 \leq 0,w_5 > p_4 - w_4$
    
    (1) If $w_4 + p_5 \leq 0$:\begin{align*}
        I_2 = -p_5.
    \end{align*}
    
    (2) If $w_4 + p_5 > 0$:\[I_2 = w_4.\]
    Combining (1) and (2), it is obvious that we take the larger of $w_4,-p_5$ , yielding:\[I_2 = \max(w_4, -p_5).\]
    
    \textbf{Case 2:}$w_4 \leq 0,0 < w_5 \leq p_4 - w_4$
    
    (1) If $p_4 + p_5 - w_5 \leq 0$:\[I_2=-p_5.\]
    
    (2) If $p_4 + p_5 - w_5 > 0$:\[I_2 = p_4 - w_5.\]
    
    From the splitting conditions (1)(2), we can get that $I_2$ always take the larger of $p_4 - w_5,-p_5$, so:\[I_2 = \max(p_4 - w_5, -p_5).\]
    
    \textbf{Case 3:} $w_4 \leq 0,w_5 \leq 0$ \[I_2=p_4.\]
    
    \textbf{Case 4:} $w_4 > 0,w_5 > p_4$\[I_2 = 0.\]
    
    \textbf{Case 5:} $w_4 > 0,0 < w_5 \leq p_4$\[I_2 = p_4 - w_5.\]
    
    \textbf{Case 6:} $w_4 > 0,w_5 \leq 0$ \[I_2 = p_4.\]

    Integrate the case when $0<w_5$:
    
    \textbf{Cases 1, 2:} $w_4 \leq 0,0 < w_5$
    
    Both cases share the term  $-p_5$, with differing terms $p_4 - w_5$ and $w_4$. Under the given conditions, in Case 1, we have  $w_5 > p_4 - w_4$, and in Case 2, we have $0 < w_5 \leq p_4 - w_4$, it is clear that $I_2$ takes the largest of $p_4 - w_5$ and  $w_4$. Thus we get:\[I_2 = \max(p_4 - w_5, w_4, -p_5).\]
    
    \textbf{Case 3:} $w_4 \leq 0,w_5 \leq 0$ \[I_2=p_4.\]
    
    \textbf{Case 4, 5:} $w_4 > 0,0 < w_5$ 
    
    Based on the conditions of Case 4  $w_5 > p_4$ and Case 5  $0 < w_5 \leq p_4$, $I_2$ takes the larger of $p_4 - w_5$ and $0$:\[I_2 = \max(p_4 - w_5, 0).\]
    
    \textbf{Case 6:} $w_4 > 0,w_5 \leq 0$\[I_2 = p_4.\]

    Combine the conditions related to $w_4$.
    
    \textbf{Cases 1, 2, 4, 5:} $0<w_5$
    
    It compares $w_4$ with $0$, leading to: \[I_2 = \max(p_4 - w_5, \min(0,w_4), -p_5).\]
    
    \textbf{Cases 3, 6:} $w_5\leq0$ \[I_2=p_4.\]

    Based on the preceding analysis, the mathematical formula for $I_2$ can be expressed as follows:\begin{align*}
        I_2 = 
    \begin{cases} 
    \max(p_4 - w_5, \min(w_4, 0), -p_5), & \textit{}{if } w_5 > 0, \\ 
    p_4, & \textit{if } w_5 \leq 0.
    \end{cases}
    \end{align*}
    
    Similarly, the conditions for the previous term $I_1$ are merged as follows:
    First, summarizing the different cases of $I_1$, we have:
    
    \textbf{Case 1:} $w_4 \leq 0,w_5 > p_4 - w_4$ \[I_1 = p_5 - (p_4 - w_4).\]
    
    \textbf{Case 2:} $w_4 \leq 0,0 < w_5 \leq p_4 - w_4$\[I_1 = p_5 - w_5.\]
    
    \textbf{Case 3:} $w_4 \leq 0,w_5 \leq 0$ \[I_1 = p_5 - w_5.\]
    
    \textbf{Case 4:} $w_4 > 0,w_5 > p_4$ \[I_1 = p_5 - p_4.\]
    
    \textbf{Case 5:} $w_4 > 0,0 < w_5 \leq p_4$\[I_1 = p_5 - w_5.\]
    
    \textbf{Case 6:} $w_4 > 0,w_5 \leq 0$\[I_1 = p_5 - w_5.\]
    Now, let us integrate the conditions under $0<w_5$:
    
    \textbf{Cases 1, 2:} $w_4 \leq 0,w_5 > 0$
    
    In both cases, $I_1$ shares the common term $p_5$. The differing term in Case 1 is $w_5 > p_4 - w_4$, taking $p_4 - w_4$, while in Case 2, it is $0 < w_5 \leq p_4 - w_4$,taking $w_5$. So, $I_1$ always takes the minimum of $w_5$ and $p_4-w_4$:\[I_1 = p_5 - \min(p_4 - w_4, w_5)\]
    
    \textbf{Case 3:} $w_4 \leq 0,w_5 \leq 0$\[I_1 = p_5 - w_5.\]
    
    \textbf{Cases 4, 5:} $w_4 > 0,w_5 >0$ 
    
    Similarly, $I_1$ takes the smaller of $p_4$ and $w_5$ \[I_1 = p_5 - \min(p_4, w_5).\]
    
    \textbf{Case 6:} $w_4 > 0,w_5 \leq 0$ \[I_1 = p_5 - w_5.\]
    
    Finally, combining the conditions that $w_4\leq0$ and $w_4>0$:
    
    \textbf{Cases 1, 2, 4, 5:} $w_5>0$
    
    Then, we compare o with  $w_4$, taking the smaller of the two.  \[I_1 = p_5 -\min(p_4 -\min(0, w_4), w_5).\]
    
    \textbf{Cases 3, 6:} $w_5 \leq 0$
    
    The expression for $I_1$ is:\[I_1 = p_5 -w_5=p_5-\min(0, w_5).\]

    According to Lemma 1, $\min(p_4 -\min(0, w_4), w_5)$ and $0$ is actually equivalent to:
    \begin{align*}
        \Delta w_6^- =\begin{cases} 
        0, & \textit{if } w_5 \leq 0, \\\min(p_4 - \min(0, w_4), w_5), & \textit{if } w_5 > 0.
    \end{cases}
    \end{align*}
    
    And $\min(0,w_5)$ takes $0$ when $w_5\leq0$, else takes $w_5$. So, $I_1$ can be expressed as:\[I_1 = p_5 -\Delta w_6^--\min(0, w_5).\]

    In summary, we have obtained the formula for the change in total waiting time $\Delta w$ when job 4 is relocated between jobs 5 and 6:\begin{align*}
        \Delta w=I_1+I_2+\pi_{6} (I_2),
    \end{align*}
    where \begin{align*}
        I_1= p_5 -\Delta w_6^--\min(0, w_5),
    \end{align*} and initial value is $\Delta w_5^-=p_4-\min(0,w_4).$
    \begin{align*}
        I_2=\begin{cases}
            p_4,&\textit{if }w_5\leq0,\\\max\left(\min(w_4,0),-p_5,p_4-w_5\right),&\textit{if }w_5>0.
        \end{cases}
    \end{align*}

    However, we can compute the updated waiting times for all jobs in the queue except for the moved job 4, when we know the initial conditions $\Delta w_5^-=p_4-\min(0,w_4)$ and $I_2$, and using the recursive formulations of Lemma 1 and Lemma 2. To complete  the full information  of the queue, it is necessary to summarize the updated waiting time  $w_4^2$ of job 4. 

    To begin, we list the calculated outcomes for $w_4^2$:
    
    \textbf{Case 1:} $w_4 \leq 0, w_5 > p_4 - w_4$\[w_4^2 = w_4 + p_5.\]
    
    \textbf{Case 2:} $w_4 \leq 0, 0 < w_5 \leq p_4 - w_4$\[w_4^2 = p_4 + p_5 - w_5.\]
    
    \textbf{Case 3:} $w_4 \leq 0, w_5 \leq 0$ \[w_4^2 = p_4 + p_5 - w_5.\]
    
    \textbf{Case 4:} $w_4 > 0, w_5 > p_4$ \[w_4^2 = p_5 + w_4.\]
    
    \textbf{Case 5:} $w_4 > 0, 0 < w_5 \leq p_4$ \[w_4^2 = p_4 + p_5 - w_5 + w_4.\]
    
    \textbf{Case 6:} $w_4 > 0, w_5 \leq 0$\[w_4^2 = p_4 + p_5 - w_5 + w_4.\]
    Analogous to the condition combination process for $I_1$, it follows:
    
    \textbf{Cases 1, 2, 3:} $w_4 \leq 0$ \[w_4^2 = p_5 + p_4 - \max(0, \min(p_4 - w_4, w_5)) - \min(0, w_5).\]
    
    \textbf{Cases 4, 5, 6:} $w_4 > 0$ \[w_4^2 = p_5 + p_4 - \max(0, \min(p_4, w_5)) - \min(0, w_5) + w_4.\]
    Finally, integrating the conditions for $w_4$, we have:
    
    \textbf{Cases 1, 2, 3, 4, 5, 6:}
    \begin{align*}
        w_4^2 &= p_5 + p_4 - \max(0, \min(p_4-\min(0,w_4), w_5)) - \min(0, w_5) + \max(0,w_4)\\&=p_5 + p_4 - \Delta w_6^-- \min(0, w_5) + \max(0,w_4).
    \end{align*}
    We summarize all the information about job 4 relocating to the position between jobs 5 and 6. The equation for the change in total waiting time is as follows:\begin{align*}
        \Delta w=I_1+I_2+\pi_{6}(I_2),
    \end{align*} 
    where \[I_1=p_5-\Delta w_6^--\min(0,w_5),\] and initial value is $\Delta w_5^-=p_4-\min(0,w_4).$
    \begin{align*}
        I_2=\begin{cases}
            p_4,&\textit{if }w_5\leq0,\\\max(\min(w_4,0),-p_5,p_4-w_5),&\textit{if }w_5>0.
        \end{cases}
    \end{align*}
    The updated waiting time of job 4 is:\begin{align*}
        w_4^2=p_4+p_5-\Delta w_6^--\min(0,w_5)+\max(0,w_4).
    \end{align*} where initial value is $\Delta w_5^-=p_4-\min(0,w_4).$
    
    We have achieved a breakthrough from “0” to “1”, by obtaining the formula above. But does this formula still hold when job 4 moves to a more forward position? Or, in other words, what is the law for extending this formula? Thus, if we want to generalize it to an infinite scope, it is necessary to examine and broaden this formula from “1” to “2”. We examine whether the law holds by shifting job 4 to the location between jobs 6 and 7 in the next step, as illustrated in Figure~\ref{forward2}.
    \begin{figure}
        \centering
        \includegraphics[width=0.8\linewidth]{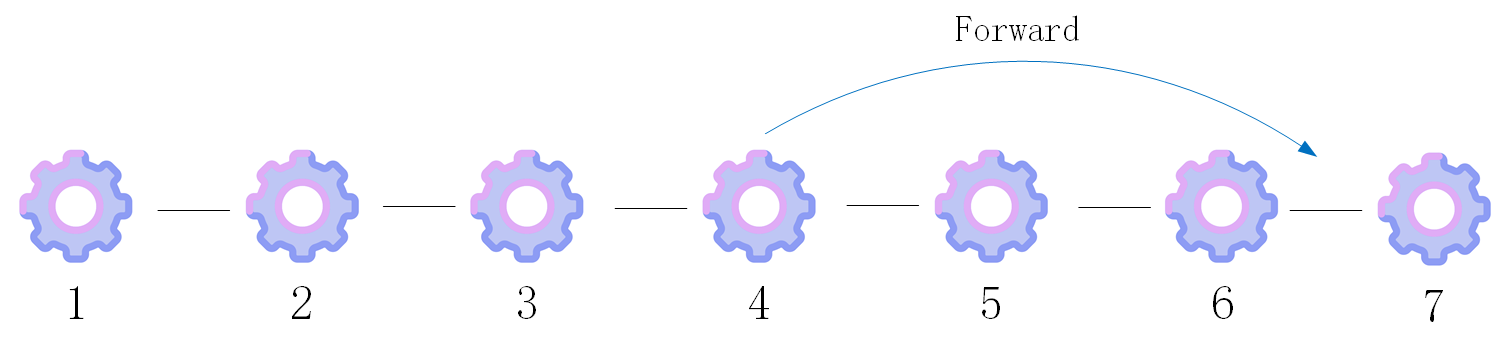}
        \caption{Fig 2. Job 4 moves to the position between jobs 6 and 7.}
        \label{forward2}
    \end{figure}
    
    Similar to the case of job 4 moving to the position between jobs 5 and 6, we also start from the special role of the queue leader and discuss it by cases:\begin{itemize}
        \item \textbf{Case 1:} $w_4 \leq 0, w_5 > p_4 - w_4, w_6 > p_4 - w_4$, initially, job 4 is the queue leader, and jobs 5 and 6 are in the contiguous queue. After moving job 4, jobs 5 and 6 still be in the contiguous queue.
        \item \textbf{Case 2:} $w_4 \leq 0, w_5 > p_4 - w_4, 0 < w_6 \leq p_4 - w_4$, in this case, job 4 is the queue leader, and jobs 5 and 6 are in the contiguous queue. Once job 4 is moved, job 5 remains in the contiguous queue, while job 6 becomes the new queue leader.
        \item \textbf{Case 3:} $w_4 \leq 0, w_5 > p_4 - w_4, w_6 \leq 0$, jobs 4 and 6 are the queue leader, and job 5 is one of the members of the contiguous queue, in the beginning. After job 4 is moved, jobs 5 and 6 are not changed.
        \item \textbf{Case 4:} $w_4 \leq 0, 0 < w_5 \leq p_4 - w_4, w_6 > w_5 $, originally, job 4 is the queue leader, and jobs 5 and 6 are in the continuous queue. Following job 4's movement, job 5 becomes the new queue leader, while job 6 remains in the contiguous queue. 
        \item \textbf{Case 5:} $ w_4 \leq 0, 0 < w_5 \leq p_4 - w_4, w_6 \leq 0$, at this point, jobs 4 and 6 are the queue leader, and job 5 is in the continuous queue. Once job 4 is relocated, both jobs 5 and 6 are the new queue leader.
        \item \textbf{Case 6:} $ w_4 \leq 0, 0 < w_5 \leq p_4 - w_4, w_6 \leq 0$, in this case, jobs 4 and 6 are the queue leader, and job 5 is part of the continuous queue. Once job 4 is relocated, both jobs 5 and 6 become the new queue leader.
        \item \textbf{Case 7:} $w_4 \leq 0, w_5 \leq 0, w_6 > 0$, in this stage, jobs 4 and 5 are the queue leader, and job 6 is in the contiguous queue. After moving job 4, jobs 5 and 6 are unchanged.
        \item \textbf{Case 8:} $w_4 \leq 0, w_5 \leq 0, w_6 \leq 0$, under this condition, jobs 4,5, and 6 are the queue leader. Following job 4's movement, jobs 5 and 6 still be the queue leader.
        \item \textbf{Case 9:} $w_4 > 0, w_5 > p_4, w_6 > p_4$, initially, jobs 4,5, and 6 are part of the continuous queue. When job 4 is relocated, jobs 5 and 6 remain unchanged.
        \item \textbf{Case 10:} $w_4 > 0, w_5 > p_4, 0 < w_6 \leq p_4$, in this case, jobs 4,5, and 6 are in the continuous queue. Once job 4 is moved, job 5 is unchanged, while job 6 shifts to the new queue leader. 
        \item \textbf{Case 11:} $w_4 > 0, w_5 > p_4, w_6 \leq 0$, here, jobs 4 and 5 are part of the continuous queue, and job 6 is the queue leader. Once job 4 is moved, jobs 5 and 6 remain their role.
        \item \textbf{Case 12:} $w_4 > 0, 0 < w_5 \leq p_4, w_6 > w_5$, initially, jobs 4, 5, and 6 are in the continuous queue. Once job 4's relocation, job 5 becomes the new queue leader, while job 6 is unchanged.
        \item \textbf{Case 13:} $w_4 > 0, 0 < w_5 \leq p_4, 0 < w_6 \leq w_5$, under this condition, jobs 4, 5, and 6 are one of the members of the continuous queue. Following job 4's movement, jobs 5 and 6 become the new queue leader.
        \item \textbf{Case 14:} $w_4 > 0, 0 < w_5 \leq p_4, w_6 \leq 0$, at this point, jobs 4 and 5 are in the continuous queue, while job 6 is the queue leader. After job 4's relocation, jobs 5 and 6 are the queue leader.
        \item \textbf{Case 15:} $w_4 > 0, w_5 \leq 0, w_6 > 0$, at this stage, jobs 4 and 6 are part of the continuous queue, and job 5 is the queue leader. After moving job 4, jobs 5 and 6 are unchanged.
        \item \textbf{Case 16:} $w_4 > 0, w_5 \leq 0, w_6 \leq 0$, job 4 is in the contiguous queue, while jobs 5 and 6 are queue heads, in the beginning. Following job 4's movement, jobs 5 and 6 remain the queue leader.
    \end{itemize}
    
    \textbf{Case 1:} $w_4 \leq 0, w_5 > p_4 - w_4, w_6 > p_4 - w_4$
    
    Initially, job 4 is the queue leader, and jobs 5 and 6 are in the contiguous queue. After moving job 4, jobs 5 and 6 still be in the contiguous queue. In this case, we can get from the conditions that:\begin{align*}
       \Delta_5 &= -\Delta w_6^- = -(p_4 - w_4),\\\Delta_6 &=- \Delta w_7^- = -(p_4 - w_4),
    \end{align*}

    The new waiting time of job 4 is:
    \begin{align*}
        w_4^2 &= r_6 + w_6^2 + p_6 - r_4\\&=r_6 - r_5 + r_5 + w_4 - p_4 + w_6 + p_6 - r_4\\&=p_5 + w_5 - w_6 + w_4 - p_4 + w_6 + p_6 + p_4 - w_5\\&=p_5 + p_6 + w_4.
    \end{align*}
    
    (1) if $p_5 + p_6 + w_4 \leq 0$, job 4 is the queue leader, we can obtain follows:\begin{align*}
        \Delta_4&=\max(0,w_4^2)-\max(0,w_4)=0,\\\Delta w_7 &= r_4^o + p_4 - r_7 - w_7\\&=r_4^o + p_4 - r_7 - (r_4^o + p_4 + p_5 + p_6 - r_7)\\&=-(p_5 + p_6).
    \end{align*}
    
    The change of total waiting time is:\begin{align*}
        \Delta w &= \Delta_4 + \Delta_5 + \Delta_6 + \pi_7 (\Delta w_7)\\&=-2(p_4 - w_4) + p_5 + p_6 - (p_5 + p_6) + \pi_7 (-(p_5 + p_6)).
    \end{align*} 
    
    (2) if $p_5 + p_6 + w_4 > 0$, job 4 is in the continuous queue:\begin{align*}
        \Delta_4&=\max(0,w_4^2)-\max(0,w_4)\\&=p_5 + p_6 + w_4.
    \end{align*}
    \begin{align*}
        \Delta w_7 &= r_4 + w_4^2 + p_4 - r_7 - w_7\\&=r_4 + p_5 + p_6 + w_4 + p_4 - r_7 - (r_4 + p_4 + p_5 + p_6 - r_7)\\&=w_4.
    \end{align*}
    
    The adjustment of total waiting time is as follows:\begin{align*}
        \Delta w &= \Delta_4 + \Delta_5 + \Delta_6 + \pi_7 (\Delta w_7)\\&=p_5 + p_6 - 2(p_4 - w_4) + w_4 + \pi_7 (w_4).
    \end{align*}
    
    \textbf{Case 2:} $w_4 \leq 0 , w_5 > p_4 - w_4 ,0 < w_6 \leq p_4 - w_4$
    
    In this case, job 4 is the queue leader, and jobs 5 and 6 are in the contiguous queue. Once job 4 is moved, job 5 remains in the contiguous queue, while job 6 becomes the new queue leader. 
    Under this conditions:\begin{align*}
        \Delta_5 &= -\Delta w_6^- = -(p_4 - w_4),\\\Delta_6 &=- \Delta w_7^- = -w_6,
    \end{align*}
    
    The updated waiting time of job 4 is:\begin{align*}
        w_4^2 &= r_6^o + p_6 - r_4 \\&= r_6^o - r_4 - p_4 - p_5 + p_4 + p_5 + p_6 \\&= p_4 + p_5 + p_6 - w_6.
    \end{align*}
    
    (1) If $p_4 + p_5 + p_6 - w_6\leq0$, job 4 still be the queue leader after its relocation, it follows:\begin{align*}
        \Delta_4 = \max(0, w_4^2) - \max(0, w_4) = 0,
    \end{align*}
    \begin{align*}
        \Delta w_7 &= r_4^o + p_4 - r_7 - w_7 \\&= r_4^o + p_4 - r_7 - (r_4^o + p_4 + p_5 + p_6 - r_7) \\&= -(p_5 + p_6).
    \end{align*}
    
    At this point, the impact on the total waiting time is formulated as follows:\begin{align*}
        \Delta w &= \Delta_4 + \Delta_5 + \Delta_6 + \pi_7 (\Delta w_7) \\&= p_5 + p_6 + (w_4 - p_4) - w_6 - (p_5 + p_6) + \pi_7 (-(p_5 + p_6)).
    \end{align*}
    
    (2) If $p_4 + p_5 + p_6 - w_6>0$, job 4 is one of the members of the continuous queue, we can obtain that:
    \begin{align*}
        \Delta_4 &= \max(0, w_4^2) - \max(0, w_4) \\&= p_4 + p_5 + p_6 - w_6,
    \end{align*}
    \begin{align*}
        \Delta w_7 &= r_4 + w_4^2 + p_4 - r_7 - w_7 \\&= r_4 + p_4 + p_5 + p_6 - w_6 + p_4 - r_7 - (r_4^o + p_4 + p_6 + p_6 - r_7) \\&= p_4 - w_6.
    \end{align*}
    Now, the impact on the total waiting time is formulated as follows:
    \begin{align*}
        \Delta w &= \Delta_4 + \Delta_5 + \Delta_6 + \pi_7 (\Delta w_7) \\&= p_5 + p_6 - w_6 - (p_4 - w_4) + p_4 - w_6 + \pi_7 (p_4 - w_6).
    \end{align*}
    
    \textbf{Case 3:} $w_4 \leq 0, w_5 > p_4 - w_4, w_6 \leq 0$
    
    Jobs 4 and 6 are the queue leader, and job 5 is one of the members of the contiguous queue, in the beginning. After job 4 is moved, jobs 5 and 6 not changed.
    This moment \begin{align*}
        \Delta_5 &= -\Delta w_6^- = -(p_4 - w_4), \\\Delta_6 &= -\Delta w_7^- = 0.
    \end{align*}
    
    The new waiting time of job 4 is\begin{align*}
        w_4^2 &= r_6^o + p_6 - r_4\\&=r_6^o + p_6 - r_4 - p_4 - p_5 + p_4 + p_5\\&=p_4 + p_5 + p_6 - w_6.
    \end{align*}
    $\because w_6\leq0,\therefore w_4^2>0$, so we have:\begin{align*}
        \Delta_4 &= \max(0, w_4^2) - \max(0, w_4) \\&= p_4 + p_5 + p_6 - w_6.
    \end{align*}
    \begin{align*}
        \Delta w_7 &= r_4 + w_4^2 + p_4 - r_7 - w_7\\&=r_4 + p_4 + (p_4 + p_5 + p_6 - w_6) - r_7 - (r_4 + p_4 + p_5 + p_6 - w_6 - r_7)\\&= p_4.
    \end{align*}
    In the end, the total waiting time calculated as follows:\begin{align*}
        \Delta w &= \Delta_4 + \Delta_5 + \Delta_6 + \pi_7 (\Delta w_7)\\&=p_5 + p_6 - (p_4 - w_4) - w_6 + p_4 + \pi_7 (p_4).
    \end{align*}
    
    \textbf{Case 4:} $w_4 \leq 0, 0 < w_5 \leq p_4 - w_4, w_6 > w_5$
    
    Originally, job 4 is the queue leader, and jobs 5 and 6 are in the continuous queue. Following job 4's movement, job 5 becomes the new queue leader, while job 6 remains in the contiguous queue.     
    Now, we can get that:\begin{align*}
        \Delta_5 &= -\Delta w_6^- = -w_5, \\\Delta_6 &= -\Delta w_7^- = -w_5.
    \end{align*}
    
    The $w_4^2$ is calculated:\begin{align*}
        w_4^2 &= r_5^o + p_5 + p_6 - r_4\\&=p_4 + p_5 + p_6 - w_5.
    \end{align*}
    
    (1) When $p_4 + p_5 + p_6 - w_5 \leq 0$, job 4 is the queue leader, following job's waiting time is:\begin{align*}
        \Delta_4 &= \max(0, w_4^2) - \max(0, w_4) = 0,\\
        \Delta w_7 &= r_4^o + p_4 - r_7 - w_7\\&=r_4^o + p_4 - r_7 - (r_4^o + p_4 + p_5 + p_6 - r_7)\\&=-(p_5 + p_6).
    \end{align*}
    
    Thus, the total waiting time change can be expressed:\begin{align*}
        \Delta w &= \Delta_4 + \Delta_5 + \Delta_6 + \pi_7 (\Delta w_7)\\&=p_5 + p_6 - w_5 - w_5 - (p_5 + p_6) + \pi_7 (-(p_5 + p_6)).
    \end{align*}
    
    (2) When $p_4 + p_5 + p_6 - w_5 > 0$, job 4 becomes part of the continuous queue.
    \begin{align*}
        \Delta_4 &= \max(0, w_4^2) - \max(0, w_4) \\&= p_4 + p_5 + p_6 - w_5.\\
        \Delta w_7 &= r_4 + w_4^2 + p_4 - r_7 - w_7\\&=r_4 + p_4 + p_5 + p_6 - w_5 + p_4 - r_7 - w_7\\&=p_4 - w_5.
    \end{align*}
    
    Finally, the change of the total waiting time is:\begin{align*}
        \Delta w &= \Delta_4 + \Delta_5 + \Delta_6 + \pi_7 (\Delta w_7)\\&=p_5 + p_6 - 2w_5 + (p_4 - w_5) + \pi_7 (p_4 - w_5).
    \end{align*}
    
    \textbf{Case 5:} $w_4 \leq 0, 0 < w_5 \leq p_4 - w_4, 0 < w_6 \leq w_5$
    
    At this point, jobs 4 and 6 are the queue leader, and job 5 is in the continuous queue. Once job 4 is relocated, both jobs 5 and 6 are the new queue leader.
    In this stage, the impact of jobs 5 and 6 on the total waiting time is\begin{align*}
        \Delta_5 &= -\Delta w_6^- = -w_5, \\\Delta_6 &= -\Delta w_7^- = -w_6.
    \end{align*}
    
    The new waiting time of job 4 is\begin{align*}
        w_4^2 &= r_6^o + p_6 - r_4\\&= p_6 + p_4 + p_5 - w_6.
    \end{align*}
    
    (1) If $p_6 + p_4 + p_5 - w_6 \leq 0 $, job 4 is the queue leader.
    \begin{align*}
        \Delta_4 &= \max(0, w_4^2) - \max(0, w_4) = 0\\
        \Delta w_7 &= r_4^o + p_4 - r_7 - w_7\\&=r_4^o + p_4 - r_7 - (r_4^o + p_4 + p_5 + p_6 - r_7)\\&=-(p_5 + p_6).
    \end{align*}
    
    In summary, the adjustment of total waiting time is:\begin{align*}
        \Delta w &= \Delta_4 + \Delta_5 + \Delta_6 + \pi_7 (\Delta w_7)\\&=p_5 + p_6 - w_5 - w_6 - (p_5 + p_6) + \pi_7 (-(p_5 + p_6)).
    \end{align*}
    
    (2) If $p_6 + p_4 + p_5 - w_6 > 0$, then job 4 is within the continuous queue.
    \begin{align*}
        \Delta_4 &= \max(0, w_4^2) - \max(0, w_4) \\&= p_6 + p_4 + p_5 - w_6.\\
        \Delta w_7 &= r_4 + w_4^2 + p_4 - r_7 - w_7\\&=r_4 + p_4 + p_5 + p_6 - w_6 + p_4 - r_7 - (r_4 + p_4 + p_5 + p_6 - r_7)\\&=p_4 - w_6.
    \end{align*}
    
    Finally, the change of total waiting time is:\begin{align*}
        \Delta w &= \Delta_4 + \Delta_5 + \Delta_6 + \pi_7 (\Delta w_7)\\&=p_5 + p_6 - w_5 - w_6 + (p_4 - w_6) + \pi_7 (p_4 - w_6).
    \end{align*}
    
    \textbf{Case 6:} $w_4 \leq 0, 0 < w_5 \leq p_4 - w_4, w_6 \leq 0$
    
    In this case, jobs 4 and 6 are the queue leader, and job 5 is part of the continuous queue. Once job 4 is relocated, both jobs 5 and 6 become the new queue leader.    
    Under this condition, it can be obtained as follows:\begin{align*}
        \Delta_5 &= -\Delta w_6^- = -w_5, \\\Delta_6 &= -\Delta w_7^- = 0.
    \end{align*}
    The new waiting time of job 4, following its movement is:\begin{align*}
        w_4^2 &= r_6^o + p_6 - r_4\\&=p_4 + p_5 + p_6 - w_6.
    \end{align*}
    $\because  w_6 \leq 0\therefore w_4^2>0.$
    \begin{align*}
        \Delta_4 &= \max(0, w_4^2) - \max(0, w_4) \\&= p_4 + p_5 + p_6 - w_6.\\
        \Delta w_7 &= r_4 + w_4^2 + p_4 - r_7 - w_7\\&=r_4 + p_4 - r_7 + p_4 + p_5 + p_6 - w_6 - w_7\\&=p_4.
    \end{align*}
    The total waiting time change is:\begin{align*}
        \Delta w &= \Delta_4 + \Delta_5 + \Delta_6 + \pi_7 (\Delta w_7)\\&=p_5 + p_6 - w_5 - w_6 + p_4 + \pi_7 (p_4).
    \end{align*}
    
    \textbf{Case 7:} $w_4 \leq 0, w_5 \leq 0, w_6 > 0$
    
    In this stage, jobs 4 and 5 are the queue leader, and job 6 is in the contiguous queue. After moving job 4, jobs 5 and 6 are unchanged.    
    In this situation, we can get that:\begin{align*}
        \Delta_5 &=- \Delta w_6^- = 0, \\\Delta_6 &= -\Delta w_7^- = 0.
    \end{align*}
    
    The updated waiting time is:\begin{align*}
        w_4^2 &= r_5^o + p_5 + p_6 - r_4\\&=p_4 + p_5 + p_6 - w_5.
    \end{align*}
    
    Since $w_5\leq0$, we have $w_4^2>0$.\begin{align*}
        \Delta_4 &= \max(0, w_4^2) - \max(0, w_4) \\&= p_4 + p_5 + p_6 - w_5.\\
        \Delta w_7 &= r_4^o + p_4 + p_5 + p_6 - r_7 - w_7\\&=p_4.
    \end{align*}
    
    Thus, the adjustment of total waiting time is:\begin{align*}
        \Delta w &= \Delta_4 + \Delta_5 + \Delta_6 + \pi_7 (\Delta w_7)\\&=p_5 + p_6 - w_5 + p_4 + \pi_7 (p_4).
    \end{align*}
   
    \textbf{Case 8:} $w_4 \leq 0, w_5 \leq 0, w_6 \leq 0 $
    
    Under this condition, jobs 4,5, and 6 are the queue leader. Following job 4's movement, jobs 5 and 6 still be the queue leader. Now,\begin{align*}
        \Delta_5 &= -\Delta w_6^- = 0,\\\Delta_6 &= -\Delta w_7^- = 0.
    \end{align*}
    
    The updated waiting time is:\begin{align*}
        w_4^2 &= r_6^o + p_6 - r_4\\&=r_6^o + p_6 - r_5 + r_5 - r_4\\&=p_4 + p_5 + p_6 - w_5 - w_6.
    \end{align*}
    $\because w_5 \leq 0, w_6 \leq 0,\therefore w_4^2 > 0$.\begin{align*}
        \Delta_4 &= \max(0, w_4^2) - \max(0, w_4) \\&= p_4 + p_5 + p_6 - w_5 - w_6.\\
        \Delta w_7 &= r_4^o + p_4 + p_5 + p_6 - r_7 - w_7\\&=p_4.
    \end{align*}
    
    Finally, the total waiting time change is:\begin{align*}
        \Delta w &= \Delta_4 + \Delta_5 + \Delta_6 + \pi_7 (\Delta w_7)\\&=p_5 + p_6 - w_5 - w_6 + p_4 + \pi_7 (p_4).
    \end{align*}
    
    \textbf{Case 9:} $w_4 > 0, w_5 > p_4, w_6 > p_4$
    
    Initially, jobs 4,5, and 6 are part of the continuous queue. When job 4 is relocated, jobs 5 and 6 remain unchanged.    
    In this case, we have:\begin{align*}
        \Delta_5 &= -\Delta w_6^- = -p_4, \\\Delta_6 &=- \Delta w_7^- = -p_4.
    \end{align*}
    
    The updated waiting time is:\begin{align*}
        w_4^2 &= r_6 + w_6^2 + p_6 - r_4\\&=r_6 + w_6 - p_4 + p_6 - r_4 - w_4 + w_4 - p_5 + p_5\\&=p_5 + p_6 + w_4.
    \end{align*}
    
    $\because w_4>0,\therefore w_4^2>0$.\begin{align*}
        \Delta_4 &= \max(0, w_4^2) - \max(0, w_4) \\&= p_5 + p_6.\\
        \Delta w_7 &= r_4 + w_4^2 + p_4 - r_7 - w_7\\&=r_4 + p_5 + p_6 + w_4 + p_4 - r_7 - w_7\\&=0.
    \end{align*}
    
    The change of total waiting time is:\begin{align*}
        \Delta w &= \Delta_4 + \Delta_5 + \Delta_6 + \pi_7 (\Delta w_7)\\&=p_5 + p_6 - 2p_4 + \pi_7 (0).
    \end{align*}
    
    \textbf{Case 10:} $w_4 > 0, w_5 > p_4, 0 < w_6 \leq p_4$
    
    In this case, jobs 4,5, and 6 are in the continuous queue. Once job 4 is moved, job 5 is unchanged, while job 6 shifts to the new queue leader.     
    At this point, it follows:
    \begin{align*}
        \Delta_5 &= -\Delta w_6^- = -p_4, \\\Delta_6 &= -\Delta w_7^- = -w_6.
    \end{align*}
    
    The updated waiting time is:\begin{align*}
        w_4^2 &= r_6^o + p_6 - r_4\\&=r_6^o + p_6 - r_4 - w_4 + w_4 - p_4 + p_4 - p_5 + p_5\\&=p_4 + p_5 + p_6 + w_4 - w_6.
    \end{align*}
    
    Since $w_6\leq p_4,w_4>0$, we have $w_4^2>0$, so:\begin{align*}
        \Delta_4 &= \max(0, w_4^2) - \max(0, w_4)\\&=w_4^2 - w_4 = p_4 + p_5 + p_6 - w_6.\\
        \Delta w_7 &= r_4 + w_4^2 + p_4 - r_7 - w_7\\&=r_4 + p_4 + p_5 + p_6 - w_6 + w_4 + p_4 - r_7 - w_7\\&=p_4 - w_6.
    \end{align*}
    
    Finally, the total waiting time change is:\begin{align*}
        \Delta w &= \Delta_4 + \Delta_5 + \Delta_6 + \pi_7 (\Delta w_7)\\&=p_5 + p_6 - (p_4 + w_6) + (p_4 - w_6) + \pi_7 (p_4 - w_6).
    \end{align*}
    
    \textbf{Case 11:} $w_4 > 0, w_5 > p_4, w_6 \leq 0$
    
    Here, jobs 4 and 5 are part of the continuous queue, and job 6 is the queue leader. Once job 4 is moved, jobs 5 and 6 remain their role.    
    In this case, we get:\begin{align*}
        \Delta_5 &= -\Delta w_6^- = -p_4, \\\Delta_6 &= -\Delta w_7^- = 0.
    \end{align*}
    
    The new waiting time of job 4 is calculated as follows:\begin{align*}
        w_4^2 &= r_6^o + p_6 - r_4\\&=r_6^o + p_6 - r_4 - w_4 + w_4 - p_4 + p_4 - p_5 + p_5\\&=p_4 + p_5 + p_6 + w_4 - w_6.
    \end{align*}
    
    Under conditions $w_4 > 0, w_6 \leq 0$, we can get $w_4^2 > 0$, so:\begin{align*}
        \Delta_4 &= \max(0, w_4^2) - \max(0, w_4)\\&=p_4 + p_5 + p_6 - w_6.\\
        \Delta w_7 &= r_4 + w_4^2 + p_4 - r_7 - w_7\\&=r_4 + p_4 + p_5 + p_6 + w_4 + p_4 - r_7 - w_7\\&=p_4.
    \end{align*}
    
    In the end, we have the formula of the total waiting time change is as follows:\begin{align*}
        \Delta w &= \Delta_4 + \Delta_5 + \Delta_6 + \pi_7 (\Delta w_7)\\&=p_5 + p_6 - (p_4 + w_6) + p_4 + \pi_7 (p_4).
    \end{align*}
    
    \textbf{Case 12:} $w_4 > 0, 0 < w_5 \leq p_4, w_6 > w_5$
    
    Initially, jobs 4, 5, and 6 are in the continuous queue. Once job 4's relocation, job 5 becomes the new queue leader, while job 6 is unchanged.    
    This moment, we have:\begin{align*}
        \Delta_5 &= -\Delta w_6^- = -w_5, \\\Delta_6 &= -\Delta w_7^- = -w_5.
    \end{align*}
    
    The new waiting of job 4 is:\begin{align*}
        w_4^2 &= r_5^o + p_5 + p_6 - r_4\\&=p_4 + p_5 + p_6 + w_4 - w_5.
    \end{align*}
    $\because  p_4 > w_5,\therefore  w_4^2 > 0$, so:\begin{align*}
        \Delta_4 &= \max(0, w_4^2) - \max(0, w_4)\\&=p_4 + p_5 + p_6 - w_5.\\
        \Delta w_7 &= r_4 + w_4^2 + p_4 - r_7 - w_7\\&=r_4 + p_4 + p_5 + p_6 + w_4 - w_5 + p_4 - r_7 - (r_4 + w_4 + p_4 + p_5 + p_6 - r_7)\\&=p_4 - w_5.
    \end{align*}
    
    In summary, the change of total waiting time is:\begin{align*}
        \Delta w &= \Delta_4 + \Delta_5 + \Delta_6 + \pi_7 (\Delta w_7)\\&=p_5 + p_6 - 2w_5 + (p_4 - w_5) + \pi_7 (p_4 - w_5).
    \end{align*}
    
    \textbf{Case 13:} $w_4 > 0, 0 < w_5 \leq p_4, 0 < w_6 \leq w_5$
    
    Under this condition, jobs 4, 5, and 6 are one of the members of the continuous queue. Following job 4's movement, jobs 5 and 6 become the new queue leader.    
    It follows:\begin{align*}
        \Delta_5 &= -\Delta w_6^- = -w_5, \\\Delta_6 &= -\Delta w_7^- = -w_6.
    \end{align*}
    
    The new waiting time of job 4 is:\begin{align*}
        w_4^2 &= r_6^o + p_6 - r_4\\&= p_4 + p_5 + p_6 + w_4 - w_6.
    \end{align*}
    $\because p_4 \geq w_5 \geq w_6,w_4>0\therefore w_4^2 > 0$. So:\begin{align*}
        \Delta_4 &= \max(0, w_4^2) - \max(0, w_4)\\&=p_4 + p_5 + p_6 - w_6.\\
        \Delta w_7 &= r_6^o + p_4 + p_6 - w_7\\&=r_6^o + p_4 + p_6 - (r_6 + w_6 + p_6 - r_7)\\&= p_4 - w_6.
    \end{align*}
    
    Consequently, the change of total waiting time is:\begin{align*}
        \Delta w &= \Delta_4 + \Delta_5 + \Delta_6 + \pi_7 (\Delta w_7)\\&=p_5 + p_6 - w_5 - w_6 + (p_4 - w_6) + \pi_7 (p_4 - w_6).
    \end{align*}
    
    \textbf{Case 14:} $w_4 > 0, 0 < w_5 \leq p_4, w_6 \leq 0$
    
    At this point, jobs 4 and 5 are in the continuous queue, while job 6 is the queue leader. After job 4's relocation, jobs 5 and 6 are the queue leader.    
    In this situation, it follows:\begin{align*}
        \Delta_5 &= -\Delta w_6^- = -w_5, \\\Delta_6 &= -\Delta w_7^- = 0.
    \end{align*}
    
    The new waiting time of job 4 is:\begin{align*}
        w_4^2 &= r_6^o + p_6 - r_4\\&=r_6^o + p_6 - r_4 - w_4 + w_4 - p_5 + p_5 - p_4 + p_4\\&=p_4 + p_5 + p_6 + w_4 - w_6.
    \end{align*}
    
    When we have those conditions: $w_6 \leq 0,w_4 > 0$, it follows that $w_4^2 > 0$:\begin{align*}
        \Delta_4 &= \max(0, w_4^2) - \max(0, w_4)\\&=p_4 + p_5 + p_6 - w_6.\\
        \Delta w_7 &= r_6^o + p_6 + p_4 - w_7\\&=r_6^o + p_6 + p_4 - (r_6^o + p_6 - r_7)\\&=p_4.
    \end{align*}
    
    Finally, we have:\begin{align*}
        \Delta w &= \Delta_4 + \Delta_5 + \Delta_6 + \pi_7 (\Delta w_7)\\&=p_5 + p_6 - w_5 - w_6 + p_4 + \pi_7 (p_4).
    \end{align*}
    
    \textbf{Case 15:} $w_4 > 0, w_5 \leq 0, w_6 > 0$
    
    At this stage, jobs 4 and 6 are part of the continuous queue, and job 5 is the queue leader. After moving job 4, jobs 5 and 6 are unchanged.    
    It can be deducted that:\begin{align*}
        \Delta_5 &= -\Delta w_6^- = 0,\\\Delta_6 &= -\Delta w_7^- = 0.
    \end{align*}
    
    The new waiting time of job 4 is:\begin{align*}
        w_4^2 &= r_5^o + p_5 + p_6 - r_4\\&=r_5^o + p_5 + p_6 - r_4 - p_4 + p_4 - w_4 + w_4\\&=p_4 + p_5 + p_6 + w_4 - w_5.
    \end{align*}
    
    Since $w_4>0,w_5\leq0$, we have $w_4^2>0$, it follows that:\begin{align*}
        \Delta_4 &= \max(0, w_4^2) - \max(0, w_4)\\&=p_4 + p_5 + p_6 - w_5.\\
        \Delta w_7 &= r_5^o + p_5 + p_6 + p_4 - r_7 - w_7\\&=p_4.
    \end{align*}
    
    Thus, the total waiting time change is calculated as follows:
    \begin{align*}
        \Delta w &= \Delta_4 + \Delta_5 + \Delta_6 + \pi_7 (\Delta w_7)\\&=p_5 + p_6 - w_5 + p_4 + \pi_7 (p_4).
    \end{align*}
    
    \textbf{Case 16:} $w_4 > 0, w_5 \leq 0, w_6 \leq 0$
    
    Job 4 is in the contiguous queue, while jobs 5 and 6 are queue heads, in the beginning. Following job 4's movement, jobs 5 and 6 remain the queue leader. In this case,  it follows that:\begin{align*}
        \Delta_5 &= -\Delta w_6^- = 0, \\\Delta_6 &= -\Delta w_7^- = 0.
    \end{align*}
    
    The updated time of job 4 is:\begin{align*}
        w_4^2 &= r_6^o + p_6 - r_4\\&=r_6 + p_6 - r_4 - w_4 - p_4 - p_5 + w_5 + p_4 + p_5 - w_5\\&=p_4 + p_5 + p_6 + w_4 - w_6 - w_5.
    \end{align*}    
    $\because w_4 > 0, w_5 \leq 0, w_6 \leq 0,$ we have $w_4^2 > 0$, it follows that:\begin{align*}
        \Delta_4 &= \max(0, w_4^2) - \max(0, w_4)\\&=p_4 + p_5 + p_6 - w_5 - w_6.\\
        \Delta w_7 &= r_6 + p_6 + p_4 - r_7 - w_7\\&=p_4.
    \end{align*}
    
    Finally, we have the change of total waiting time is:\begin{align*}
        \Delta w &= \Delta_4 + \Delta_5 + \Delta_6 + \pi_7 (\Delta w_7)\\&=p_5 + p_6 - w_5 - w_6 + p_4 + \pi_7 (p_4).
    \end{align*}

    In the same way, we consolidate and integrate $w_4$, $I_2$, and $I_1$ following the approach used for the job 4 forward movement only one position.
    
    First, let us integrate the conditions for $I_2$: 
    
     \textbf{Cases 1, 2, 4, 5, 9, 10, 12, and 13} $w_5 > 0, w_6 > 0$
     
     \[I_2=\max(\min(0,w_4),-(p_5+p_6),p_4-w_5,p_4-w_6).\]
     
     \textbf{Case 3}: $w_4\leq0,w_5>p_4-w_4,w_6\leq0$
     
     \[I_2=p_4.\]
     
     \textbf{Case 6}: $w_4\leq0,0<w_5\leq p_4-w_4,w_6\leq0$\[I_2=p_4.\]
     
     \textbf{Case 7}: $w_4\leq0,w_5\leq0,w_6>0$\[I_2=p_4.\]
     
     \textbf{Case 8}: $w_4\leq0,w_5\leq0,w_6\leq0$\[I_2=p_4.\]
     
     \textbf{Case 11}: $w_4>0,w_5>p_4,w_6\leq0$\[I_2=p_4.\]
     
     \textbf{Case 14}: $w_4>0,0<w_5\leq p_4,w_6\leq0$\[I_2=p_4.\]
     
     \textbf{Case 15}: $w_4>0,w_5\leq0,w_6>0$\[I_2=p_4.\]
     
     \textbf{Case 16}: $w_4>0,w_5\leq0,w_6\leq0$\[I_2=p_4.\]
     
     From observation, it is obvious that when either $w_5\leq 0$ or $w_6 \leq 0$, we have that:  \[I_2 = p_4.\] If both them are greater than 0: \begin{align*}
         I_2 &= \max(\min(w_4, 0), -(p_5 + p_6), p_4 - w_5, p_4 - w_6) \\&=\max\left(\min(w_4, 0), -(\sum_{k=5}^6 p_k), \max_{k\in{5,6}} (p_4 - w_k)\right).
     \end{align*} 

     Therefore, $I_2$ can be expressed mathematically as:\begin{align*}
     I_2=
         \begin{cases}
             p_4,&\textit{otherwise},\\\max\left(\min(w_4, 0), -(\sum_{k=5}^6 p_k), \max_{k\in{5,6}} (p_4 - w_k)\right),&\textit{if } \boldsymbol{w_k}\succ0,k\in {\{5,6\}}.
         \end{cases}
     \end{align*}
     
     Likewise, we proceed to combine the conditions for the $I_1$:
         
    \textbf{Cases 1, 2, 4, 5, 9, 10, 12, 13:} $w_5 > 0, w_6 > 0$\[I_1 = p_5 + p_6 - \min(p_4 - \min(w_4, 0), w_5) - \min(p_4 - \min(w_4, 0), w_5, w_6).\]    
    
    \textbf{Cases 3, 6, 11, 14:} $w_5 > 0, w_6 \leq 0$\[I_1 = p_5 + p_6 - \min(p_4 - \min(w_4, 0), w_5) - \min(w_6, 0).\]
    
    \textbf{Cases 7, 15:} $w_5 \leq 0, w_6 > 0$\[I_1 = p_5 + p_6 - \min(w_5, 0).\]
    
    \textbf{Cases 8, 16:} $w_5 \leq 0, w_6 \leq 0$\[I_1 = p_5 + p_6 - \min(w_5, 0) - \min(w_6, 0).\]
    
    If we omit $\min(w_5, 0)$ and $\min(w_6, 0)$, it can observe:\begin{itemize}
        \item When $w_5 > 0$ and $w_6 > 0$, we have\begin{align*}
            \Delta w_6^-&= \min(p_4 - \min(w_4, 0), w_5),\\  \Delta w_7^-&=\min(p_4 - \min(w_4, 0), w_5, w_6).
        \end{align*}
        \item When $w_5 > 0$ and $w_6 \leq 0$, leading to:\begin{align*}
            \Delta w_6^-&= \min(p_4 - \min(w_4, 0), w_5),\\  \Delta w_7^-&=0.
        \end{align*}
        \item When $w_5 \leq 0$, the term is $\Delta w_6^-=0$ and $\Delta w_7^-=0$.
    \end{itemize}
    
    Obviously, it is conforms to the law of Lemma 1. So, $I_1$ can be expressed as:\begin{align*}
        I_1 &= p_5 + p_6 - (\Delta w_6^- + \min(w_5,0)) - (\Delta w_7^- + \min(w_6,0))\\&=\sum_{k=5}^{6} \left( p_k - \min(w_k, 0) \right) - \sum_{k=6}^{7} \Delta w_k^-.
    \end{align*}

    In the next step, we combine the conditions for $w_4^2$: 
         
    \textbf{Cases 1, 2, 4, 5, 9, 10, 12, 13:}$w_5>0,w_6>0$\[w_4^2=p_4+p_5+p_6-\min(w_5,w_6,p_4-\min(0,w_4 ))+\max(0,w_4).\]
    
    \textbf{Cases 3, 6, 7, 8, 11, 14, 15, 16:} $w_6\leq0$\[w_4^2=p_4+p_5+p_6-\min(0,w_5)-\min(0,w_6)+\max(0,w_4).\]

    According to the regularity of Lemma 1, it is obviously that $\Delta w_7^-=\min(w_5, w_6, p_4 - \min(0, w_4))$ when $w_5 > 0$ and $w_6 > 0$. When either $w_5$ or $w_6$ is less than or equal to 0, $\Delta w_7^-=0$. So, we have:\[w_4^2=p_4+p_5+p_6-\Delta w_7^--\min(0,w_5)-\min(0,w_6)+\max(0,w_4).\]

    \textbf{Conclusion:} Now, that we have obtained all the information for job 4 moving to a position between jobs 6 and 7, the calculation formula for the change of total waiting time is:\[\Delta w=I_1+I_2+\pi_{7}(I_2).\] where \begin{align*}
        I_1=\sum_{k=5}^6 (p_k-\min(0,w_k ))-\sum_{k=6}^7 \Delta w_k^-,
    \end{align*} 
    
    \begin{align*}
        I_2=\begin{cases}
            p_4,&\textit{otherwise,}\\\max\left(\min(w_4,0),-(\sum_{k=5}^6 p_k),\max_{k\in\{5,6\}}(p_4-w_k)\right),&\textit{if } \boldsymbol{w_k}\succ0,k\in\{5,6\}.
        \end{cases}
    \end{align*} 
    
    The new waiting time of job 4 is:\begin{align*}
        w_4^2=\sum_{k=4}^6 p_k -\sum_{k=5}^6 \min(0,w_5) -\Delta w_7^-+\max(0,w_4).
    \end{align*} 

    And the initial value is $\Delta w_5^-=p_4-\min(0,w_4)$.

    So far, we have established the law of the job moves forward. It can be described as: Given the current job sequence $\sigma_t$, if job $\sigma_t(i)$ is moved forward to the position immediately after job $\sigma_t(k)$ (where $k \in [i+1,n]$), the change in total waiting time $\Delta w$ follows this formula:\begin{align*}
        \Delta w = I_1 + I_2 + \pi_{\sigma_t(k+1)}(I_2),
    \end{align*}
    where\begin{align*}
         I_1 = \sum_{j=i+1}^{k} \left(p_{\sigma_t(j)}- \min( 0,w_{\sigma_t(j)})\right)- \sum_{j=i+2}^{k+1} \Delta w_{\sigma_t(j)}^-,
    \end{align*} 

    \begin{align*}
        I_2 = \begin{cases} 
            p_{\sigma(i)}, 
                & \text{\textnormal{otherwise}}, \\\max \Bigl( \begin{aligned}[t]\min(0,w_{\sigma(i)} ),\quad -{\sum_{j=i+1}^{k}} p_{\sigma(j)}, \quad\max_{\mathclap{\substack{j \in [i+1,k]}}} \ (p_{\sigma(i)} -w_{\sigma(j)}) \Bigr),\end{aligned} &\text{\textnormal{if} } \mathbf{w}_{\sigma(j)} \succ \mathbf{0}, \forall j \in [i+1,k].
        \end{cases}
    \end{align*}

    Furthermore, after the relocation of job $\sigma_t(i)$, its new waiting time $w_{\sigma_t(i)}^2$ is calculated as:
    \begin{align*}
        w_{\sigma_t(i)}^2 = \sum_{j=i}^{k} p_{\sigma_t(j)} -\sum_{j=i+1}^{k} \min(0, w_{\sigma_t(j)})- \Delta w_{\sigma_t(k+1)}^-+\max(0, w_{\sigma_t(i)}).
    \end{align*} 
    
    And $\Delta w_{\sigma_t(j)}^-$ follows  the recursive formula given in Lemma 1, with the initial value:\begin{align*}\Delta w_{\sigma_t(i+1)}^- = p_{\sigma_t(i)} - \min( 0,w_{\sigma_t(i)}).\end{align*}
\end{proof}
\subsection{ Proof of \cref{Theorem2}}\label{app:Theorem2}
\begin{proof}
   The proof of \cref{Theorem2} is analogous to that of \cref{Theorem1}. Consider an arbitrary schedule $\sigma_t$ and any four consecutive jobs $\sigma_t(i)$, $\sigma_t(i+1)$, $\sigma_t(i+2)$, and $\sigma_t(i+3)$, with $i \in {1,\dots,n-3}$. We start from the simplest case, in which job $\sigma_t(i+2)$ is moved to the position immediately preceding $\sigma_t(i+1)$, as shown in Figure~\ref{forward1}. Without loss of generality, denote these four jobs by $2$, $3$, $4$, and $5$, respectively; that is, let $\sigma_t(i)=2$, $\sigma_t(i+1)=3$, $\sigma_t(i+2)=4$, and $\sigma_t(i+3)=5$. This relabeling is adopted purely for notational simplicity and has no effect on the validity of the subsequent analysis.
    
    Inspired by the proof of Theorem 1, we adjust the formula toward the form of $\Delta w_j^+$. The notation is the same as the proof of Theorem 1. In the same way, $\Delta_5 = \Delta w_6^+$ is used to calculate the impact of non-moving jobs on the total waiting time, while $\Delta_5 = \max(0, w_5^2) - \max(0, w_5)$ calculates the effect of the moving job. Beginning with the special role of the queue leader, we classify the scenarios.
    
    \begin{figure}
        \centering
        \includegraphics[width=0.8\linewidth]{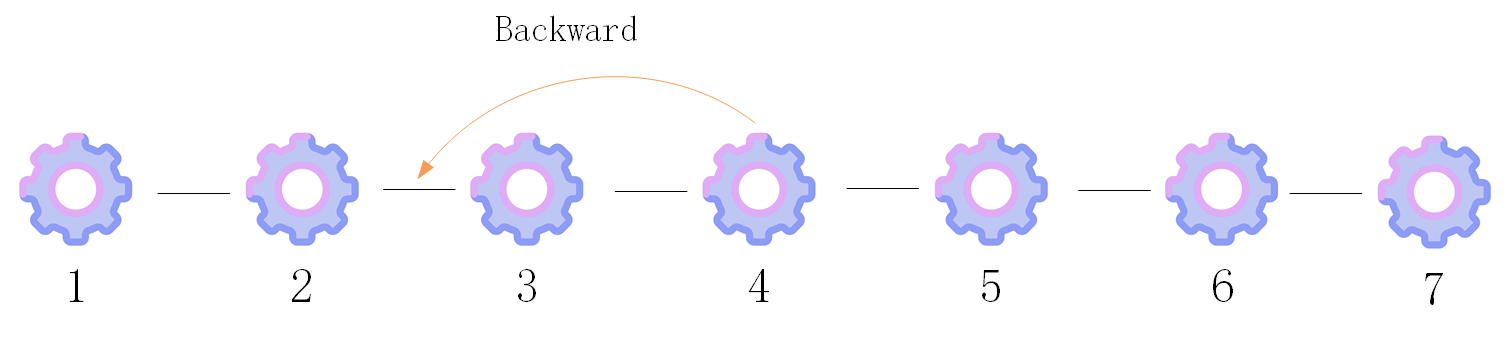}
        \caption{Fig 3. Job 4 is moved to the position between jobs 3 and 4.}
        \label{backward1}
    \end{figure}

    \textbf{Case 1:} $w_4\leq0,w_3\leq0$
    
    Initially, both jobs 3 and 4 are the queue leader.    
    In this case, we can calculate the change of job 4:\begin{align*}
        \Delta w_4&=r_0^o+p_0-r_4-w_4\\&=r_0^o+p_0-r_4-(r_0^o+p_0+p_3-r_4 )\\&=w_3-p_3.
    \end{align*}
    The updated waiting time of job 4 is:\begin{align*}
        w_4^2&=w_4+\Delta w_4\\&=w_3+w_4-p_3.
    \end{align*}Since $w_4\leq0,w_3<0$, we can get $w_4^2\leq0$. Then, job 4 still be the queue leader after its relocation.
    The impact of job 4 on the total waiting time is:\[\Delta_4=\max(0,w_4^2 )-\max(0,w_4 )=0.\]
    The updated information of job 3 is as follows:\begin{align*}
        \Delta w_3&=r_4^o+p_4-r_3-w_3\\&=p_3+p_4-w_4.\\
        w_3^2&=w_3+\Delta w_3\\&=p_3+p_4+w_3-w_4>0.\\
        \Delta_3&=\Delta w_4^+=p_3+p_4-w_4.
    \end{align*}
    Because $w_3^2>0$, job 3 is part of the continuous queue. The change of job 5 is:\begin{align*}
        \Delta w_5&=r_3+w_3^2+p_3-r_5-w_5\\&=r_3+p_3+p_4+w_3-w_4+p_3-r_5-(r_3+w_3+p_3+p_4-r_5-w_4 )\\&=p_3.
    \end{align*}
    The change of total waiting time is calculated as below:\begin{align*}
        \Delta w&=\Delta_3+\Delta_4+\pi_5 (\Delta w_5 )\\&=p_4-w_4+p_3+\pi_5 (p_3).
    \end{align*}
    
    \textbf{Case 2:} $w_4\leq0,w_3>0$
    
    Job 4 is the queue leader, while job 3 is part of the continuous queue.    
    Here, the change waiting time for job 4, and its new waiting time are, respectively:\begin{align*}
        \Delta w_4&=r_0^o+p_0-r_4-w_4\\&=r_0^o+p_0-r_4-(r_0^o+p_0+p_3-r_4 )\\&=-p_3.\\
        w_4^2&=w_4+\Delta w_4\\&=w_4-p_3\leq0.
    \end{align*}    
    For $w_4^2\leq0$, job 4 is the queue leader even its movement. The impact of job 4 on the total waiting time is:\begin{align*}
        \Delta_4=\max(0,w_4^2 )-\max(0,w_4 )=0.
    \end{align*}
    The new information about job 3 is:\begin{align*}
        \Delta w_3&=r_4^o+p_4-r_3-w_3\\&=p_3-w_4+p_4.\\
        w_3^2&=w_3+\Delta w_3\\&=p_3+p_4+w_3-w_4>0.\\
        \Delta_3&=\Delta w_4^+=p_3+p_4-w_4.
    \end{align*}
    For $w_3^2>0$, job 3 is in the continuous queue. The change of job 5's waiting time is as follows:\begin{align*}
        \Delta w_5&=r_3+w_3^2+p_3-r_5-w_5\\&=r_3+p_3+w_3-w_4+p_3-r_5+p_4-(r_3+w_3+p_3+p_4-r_5-w_4 )\\&=p_3.
    \end{align*}
    The impact on the total waiting time is:\begin{align*}
        \Delta w&=\Delta_3+\Delta_4+\pi_5 (\Delta w_5 )\\&=p_4-w_4+p_3+\pi_5 (p_3).
    \end{align*}
    
    \textbf{Case 3:} $w_4>0,w_3\leq0$
    
    Job 4 is one of the members of the continuous queue, job 3 is the queue leader.    
    Under this condition, the change of job 4's waiting time, and its new waiting time are, respectively:
    \begin{align*}
        \Delta w_4&=r_0^o+p_0-r_4-w_4\\&=r_0^o+p_0-r_4-(r_0^o+p_0-w_3+p_3-r_4 )\\&=w_3-p_3.\\
        w_4^2&=w_4+\Delta w_4\\&=w_3+w_4-p_3.
    \end{align*}
    
    Comparing 0 with $w_4^2$:
    
    (1) If $w_3+w_4-p_3\leq0$, job 4 is the queue leader, and its impact on the total waiting time is:\begin{align*}
        \Delta_4=\max(0,w_4^2 )-\max(0,w_4 )=-w_4.
    \end{align*}
    
    The information about job 3, after job 4's relocation, is listed:\begin{align*}
        \Delta w_3&=r_4^o+p_4-r_3-w_3\\&=r_4^o+p_4-r_3-w_3+p_3-p_3\\&=p_3+p_4-w_3-w_4.\\
        w_3^2&=\Delta w_3+w_3\\&=p_3+p_4-w_4.
    \end{align*}
    
    It is necessary to compare 0 with $w_3^2$: 
    
    1) When $p_3+p_4-w_4\leq0$, job 3 is the queue leader.
    The impact of job 3 on the total waiting time is:\[\Delta_3=\Delta w_4^+=0.\]
    
    The adjustment of job 5 is as follows:\begin{align*}
        \Delta w_5&=r_3^o+p_3-r_5-w_5\\&=r_3^o+p_3-r_5-(r_3^o+p_3+p_4-r_5 )\\&=-p_4.
    \end{align*}
    
    The change of total waiting time is:\begin{align*}
        \Delta w&=\Delta_3+\Delta_4+\pi_5 (\Delta w_5 )\\&=p_4-w_4-p_4+\pi_5 (-p_4).
    \end{align*}
    
    2) When $p_3+p_4-w_4>0$, it means job 3 is part of the continuous queue.
    The impact of job 3' waiting time change on the total waiting time is:\[\Delta_3=\Delta w_4^+=p_3+p_4-w_4.\]
    
    The change of job 5's waiting time is:\begin{align*}
        \Delta w_5&=r_3+w_3^2+p_3-r_5-w_5\\&=r_3^o+p_3+p_4-w_4+p_3-r_5-(r_3^o+p_3+p_4-r_5 )\\&=p_3-w_4.
    \end{align*}
    
    The change of total waiting time is:\begin{align*}
        \Delta w &= \Delta_3 + \Delta_4 + \pi_5 (\Delta w_5)\\&=p_4 - w_4 + p_3 - w_4 + \pi_5 (p_3 - w_4).
    \end{align*}
    
    (2) If $w_3+w_4-p_3>0$, job 4 is one of the members of the continuous queue, its impact on the total waiting time is:\[\Delta_4=w_3-p_3.\]
    
    The new information of job 3 is listed:\begin{align*}
        \Delta w_3&=r_4^o+w_4^2+p_4-r_3-w_3\\&=r_4+w_4+w_3-p_3+p_4-r_3-w_3\\&=r_3^o+p_3-w_4+w_4+w_3-p_3+p_4-r_3-w_3\\&=p_4.
        \\w_3^2&=\Delta w_3+w_3\\&=w_3+p_4.
    \end{align*}
    
    It needs to compare 0 with $w_3^2$ as well.
    
    1) When $w_3+p_4>0$, job 3 is in the continuous queue, its impact on the total waiting time is:\[\Delta_3=\Delta w_4^+=w_3+p_4.\]
    
    The change of job 5 is:\begin{align*}
        \Delta w_5&=r_3+w_3^2+p_3-r_5-w_5\\&=r_3+w_3+p_4+p_3-r_5-(r_3^o+p_3+p_4-r_5 )\\&=w_3.
    \end{align*}
    
    Then, the adjustment of total waiting time is:\begin{align*}
        \Delta w &= \Delta_3 + \Delta_4 + \pi_5 (\Delta w_5)\\&=p_4-p_3+w_3+w_3+\pi_5 (w_3).
    \end{align*}
    
    2) When $w_3+p_4\leq0$, job 3 is the queue leader, it changes the total waiting time by:\[\Delta_3=\Delta w_4^+=0.\]
    
    The change of job 5 is:\begin{align*}
        \Delta w_5&=r_3^o+p_3-r_5-w_5\\&=r_3^o+p_3-r_5-(r_3+p_3+p_4-r_5 )\\&=-p_4.
    \end{align*}
    
    The change of total waiting time is as follows:\begin{align*}
         \Delta w &= \Delta_3 + \Delta_4 + \pi_5 (\Delta w_5)\\&=w_3-p_3+p_4-p_4+\pi_5 (-p_4).
    \end{align*}
    
    \textbf{Case 4:} $w_4>0,w_3>0$, both jobs 3 and 4 are the members of the continuous queue.
    
    At this point, the change of job 4's waiting time is:\begin{align*}
        \Delta w_4&=r_0^o+p_0-r_4-w_4\\&=r_0^o+p_0-r_4-(r_0^O+p_0+p_3-r_4 )\\&=-p_3.
    \end{align*}
    
    The updated waiting time of job 4 is:\begin{align*}
        w_4^2&=w_4+\Delta w_4\\&=w_4-p_3.
    \end{align*}
    
    It needs to know whether job 4 is the queue leader, therefore comparing 0 with $w_4^2$.
    
    (1) If $w_4-p_3>0$, job 4 is part of the continuous queue, while the change of total waiting time caused by job 4 is:\[\Delta_4=\max(0,w_4^2 )-\max(0,w_4 )=-p_3.\]
    
    The information about job 3, once job 4 is moved, is:\begin{align*}
        \Delta w_3&=r_4+w_4^2+p_4-r_3-w_3\\&=r_4+w_4-p_3+p_4-r_3-w_3\\&=p_4.\\
        w_3^2&=\Delta w_3+w_3\\&=w_3+p_4>0.
    \end{align*}
    
    For $w_3^2>0$, job 3 is in the continuous queue, and its impact on the total waiting time is:\[\Delta_3=\Delta w_4^+=p_4.\]
    
    The change of job 5's waiting time is:\begin{align*}
        \Delta w_5&=r_3+w_3^2+p_3-r_5-w_5\\&=r_3+w_3+p_4+p_3-r_5-(r_3+w_3+p_3+p_4-r_5 )\\&=0.
    \end{align*}
    
    The total waiting time is changed as follows:\begin{align*}
        \Delta w &= \Delta_3 + \Delta_4 + \pi_5 (\Delta w_5)\\&=p_4-p_3+\pi_5 (0).
    \end{align*}
    
    (2) If $w_4-p_3\leq0$, job 4 becomes the queue leader, and the total waiting change caused by job 4 is:\[\Delta_4=\max(0,w_4^2 )-\max(0,w_4 )=-w_4.\]
    
    The new information of job 3 is as follows:\begin{align*}
        \Delta w_3&=r_4^o+p_4-r_3-w_3\\&=r_3+w_3+p_3-w_4+p_4-r_3-w_3\\&=p_3+p_4-w_4.\\
        w_3^2&=w_3+\Delta w_3\\&=p_3+p_4-w_4+w_3>0.
    \end{align*}
    
    Because $w_3^2>0$, job 3 is one of the members of the continuous queue, and its impact on the total waiting time is:\[\Delta_3=\Delta w_4^+=p_3+p_4-w_4.\]
    
    Then, the change of total waiting time is:\begin{align*}
        \Delta w &= \Delta_3 + \Delta_4 + \pi_5 (\Delta w_5)\\&=p_4-w_4+p_3-w_4+\pi_5 (p_3-w_4).
    \end{align*}
    
    Following the way of Theorem 1’s proof, we partition the formula according to the items in $\pi(.)$, denoting the terms in $\pi(.)$ as $I_4$ and the remaining terms as $I_3$.
    
    To begin, we summarize $I_4$ to get the expression:\[I_4=\max(-p_4,\min(w_3,0),p_3-\max(w_4,0)).\]

     $I_3$ doesn’t fully depend on the splitting condition and needs to move toward the form of $\Delta w_j^+$. Here, we provide the detailed process.
     
     \textbf{Case 1:} $w_4\leq0,w_3\leq0$\[I_3=p_4-w_4.\]
     
     \textbf{Case 2:} $w_4\leq0,w_3>0$\[I_3=p_4-w_4.\]
     
     \textbf{Case 3:} $w_4>0,w_3\leq0$
     
     (1) If $w_3+w_4-p_3\leq0$:\[I_3=p_4-w_4.\]
     
     (2) If $w_3+w_4-p_3>0$:\[I_3=p_4-p_3+w_3.\]
     
     \textbf{Case 4:} $w_4>0,w_3>0$
     
     (1) If $w_4-p_3>0$:\[I_3=p_4-p_3.\]
     (2) If $w_4-p_3\leq0$:\[I_3=p_4-w_4.\]

     Combining the conditions for (1) and (2):
     
     \textbf{Case 1:} $w_4\leq0,w_3\leq0$\[I_3=p_4-w_4.\]
     
     \textbf{Case 2:} $w_4\leq0,w_3>0$\[I_3=p_4-w_4.\]
     
     \textbf{Case 3:} $w_4>0,w_3\leq0$
     
     In this case, it compares $-w_4$ with $w_3-p_3$, and $I_3$ takes the greater of them:\[I_3=p_4+\max(-p_3+w_3,-w_4).\] 
    
     \textbf{Case 4:} $w_4>0,w_3>0$
     
     In the same way, $I_3$ gets the bigger one of $-p_3$ and $-w_4$: \[I_3=p_4+\max(-p_3,-w_4).\]

     As mentioned at the beginning, the formula's form should be similar to $\Delta w_j^+$. We perform an equivalent reformulation of $I_3$ for each case:
     
     \textbf{Case 1:} $w_4\leq0,w_3\leq0$\begin{align*}
         I_3&=p_4-w_4\\&=p_4-w_4+w_3-p_3-w_3+p_3\\&=\max(p_4,p_3+p_4-w_3-w_4 )+w_3-p_3.
     \end{align*}
     For $w_4\leq0,w_3\leq0$, so $p_3+p_4-w_3-w_4\geq p_4.$ Thus, $\max(p_4,p_3+p_4-w_3-w_4 )$ always taking $p_3+p_4-w_3-w_4$, above equation holds.
     
     \textbf{Case 2:} $w_4\leq0,w_3>0$\begin{align*}
         I_3&=p_4-w_4\\&=p_4-w_4-p_3+p_3\\&=\max(p_4,p_4+p_3-w_4 )-p_3.
     \end{align*}
     Similarly, $w_4\leq0$, we can get $p_4+p_3-w_4>p_4$. So $\max(p_4,p_4+p_3-w_4 )$ always is $p_4+p_3-w_4$, the equation holds.
     
     \textbf{Case 3:} $w_4>0,w_3\leq0$ \begin{align*}
         I_3&=p_4+\max(-w_4,w_3-p_3 )\\&=\max(p_4,p_3+p_4-w_3-w_4 )+w_3-p_3.
     \end{align*}
     It only involves adjusting the terms in $\max$ by a constant, without changing their relative magnitudes, thus it will not change the result of the $\max$ operation as well.
     
     \textbf{Case 4:} $w_4>0,w_3>0$\begin{align*}
         I_3&=p_4+\max(-p_3,-w_4 )\\&=\max(p_4,p_3+p_4-w_4 )-p_3.
     \end{align*}
     Similarly to the Case 3, the equation holds.

     Performing combination of similar terms leads to:
     
     \textbf{Case 1, 3:} $w_3\leq0$\begin{align*}
         I_3=\max(p_4,p_3+p_4-w_3-w_4 )+w_3-p_3.
     \end{align*}
    
     \textbf{Case 2, 4:} $w_3>0$\begin{align*}
         I_3=\max(p_4,p_4+p_3-w_4 )-p_3.
     \end{align*}

     Integrating the conditions for $w_3$, we have:\begin{align*}
         I_3&=\max(p_4,p_3+p_4-\min(0,w_3 )-w_4 )+\min(0,w_3 )-p_3\\&=\Delta w_3^++\min(0,w_3)-p_3.
     \end{align*} 

     Now, we have the information of $I_3$ and $I_4$, and only $w_4^2$ is missing to know everything. $w_4^2$ is summarized as:\[w_4^2=\min(0,w_3)+w_4-p_3.\]

     In the end, we get the formula of the change of total waiting time, when job 4 moves to the position between job 2 and 3:\begin{align*}
         \Delta w=I_3+I_4+\pi_5 (I_4 ).
     \end{align*} where \begin{align*}
         I_3&=\Delta w_3^++\min(0,w_3)-p_3,\\I_4&=\max(-p_4,\min(w_3,0),p_3-\max(w_4,0)).
     \end{align*} 
     The initial value is $\Delta w_3^+=\max(p_4,p_3+p_4-\min(0,w_3 )-w_4 ).$
     
     The updated waiting time of job 4 can be calculated as:\[w_4^2=\min(0,w_3 )-p_3+w_4.\]

     We have achieved the breakthrough from "0" to "1", we next verify the law from "1" to "2" by relocating job 4 to the position between jobs 1 and 2. There are eight cases according to the queue leader.
     \begin{figure}
         \centering
         \includegraphics[width=0.8\linewidth]{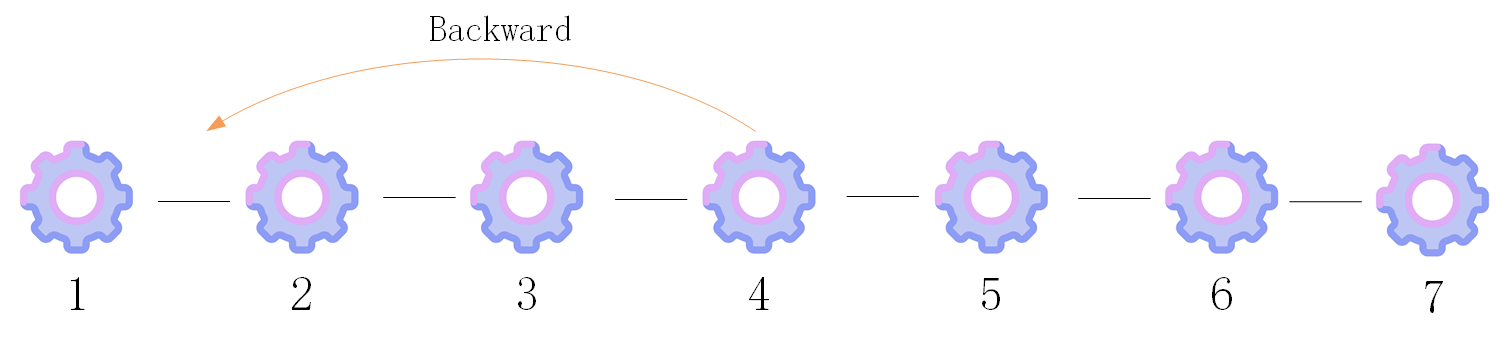}
         \caption{Fig 4. Job 4 moves to the position between job 1 and 2.}
         \label{backward2}
     \end{figure}
     
     \textbf{Case 1:} $w_4\leq0,w_2\leq0,w_3\leq0$
     
     When job 4 is not relocated, jobs 2, 3, and 4 are each the queue leader.     
     In this case, the waiting time's change of job 4 and its new waiting time are respectively:\begin{align*}
         \Delta w_4&=r_0^o+p_0-r_4-w_4\\&=r_0^o+p_0-r_4-(r_0^o+p_0-w_2+p_2-w_3+p_3-r_4 )\\&=w_2+w_3-p_2-p_3.\\
         w_4^2&=w_4+\Delta w_4\\&=w_2+w_3+w_4-p_2-p_3>0.
     \end{align*}
     
     For $w_4^2>0$, job 4 is part of the continuous queue. Once job 4 is moved, the impact of job 4 on the total waiting time is:\begin{align*}
     \Delta_4 &= \max(0, w_4^2) - \max(0, w_4)\\&=w_2 + w_3 + w_4 - p_2 - p_3.
     \end{align*}     
     
     The updated information about job 2 is:\begin{align*}
         \Delta w_2 &= r_4^o + p_4 - r_2 - w_2\\&= r_2^o + p_2 - w_3 + p_3 - w_4 + p_4 - r_2 - w_2\\&=p_2 + p_3 + p_4 - w_2 - w_3 - w_4.\\
         w_2^2 &= w_2 + \Delta w_2\\&=p_2 + p_3 + p_4 - w_3 - w_4 > 0.
     \end{align*}
     
     In this point, job 2 is one of the members of the continuous queue. The affect of job 2 on the total waiting time is:\[\Delta_2=\Delta w_3^+=p_2+p_3+p_4-w_3-w_4.\]
     
     The new information of job 3 is:\begin{align*}
         \Delta w_3& = r_4^o + p_4 + p_2 - r_3 - w_3\\&=r_4^o + p_4 + p_2 - (r_3 + p_3 - p_3) - w_3\\&=p_2 + p_3 + p_4 - w_4 - w_3.\\
         w_3^2 &= w_3 + \Delta w_3\\&=p_2 + p_3 + p_4 - w_4 > 0.
     \end{align*}
     
     Because $w_3^2>0$, job 3 is in the continuous queue, the impact of job 3 on the total waiting time is:\[\Delta_3=\Delta w_4^+=p_2+p_3+p_4-w_4.\]
     
     The change of job 5's waiting time is:\begin{align*}
         \Delta w_5 &= r_4^o + p_4 + p_3 + p_2 - r_5 - w_5\\&=r_4^o + p_2 + p_3 + p_4 - r_5 - (r_4^o + p_4 - r_5)\\&=p_2 + p_3.
     \end{align*}
     
     The adjustment of total waiting time is:\begin{align*}
         \Delta w &= \Delta_2 + \Delta_3 + \Delta_4 + \pi_5 (\Delta w_5)\\&=p_2 + p_3 + p_4 - w_3 - w_4 + p_2 + p_3 + p_4 - w_4 + 0 + \pi_5 (p_2 + p_3)\\&=p_2 + p_3 + p_4 - w_3 - w_4 + p_4 - w_4 + p_2 + p_3 + \pi_5 (p_2 + p_3).
     \end{align*}
     
     \textbf{Case 2:} $w_4\leq0,w_2\leq0,w_3>0$
     
     Before job 4's relocation, both jobs 2 and 4 are the queue leader, and job 3 in the continuous queue.
     
     The change of information about job 4 is:\begin{align*}
         \Delta w_4&=r_0^o+p_0-r_4-w_4\\&=r_0^o+p_0-r_4-(r_0^o+p_0-w_2+p_2+p_3-r_4 )\\&=w_2-p_2-p_3,\\
         w_4^2&=w_4+\Delta w_4\\&=w_4+w_2-p_2-p_3<0.
     \end{align*}
     
     Job 4 is the queue leader now, the impact of job 4 on the total waiting time is:\[\Delta_4=\max(0,w_4^2 )-\max(0,w_4 )=0.\]
     
     The updated waiting time of job 2 and its new waiting time are:\begin{align*}
         \Delta w_2&=r_4^o=p_4-r_2-w_2.\\w_2^2&=w_2+\Delta w_2\\&=p_2+p_3+p_4-w_4>0.
     \end{align*}
     
     For $w_2^2>0$, job 2 is in the continuous queue, its impact on the total waiting time is:\begin{align*}
         \Delta_2=\Delta w_3^+=p_2+p_3+p_4-w_4.
     \end{align*}
     
     The new information about job 3 is:\begin{align*}
         \Delta w_3&=r_2+w_2^2+p_2-r_3-w_3\\&=r_2+p_2+p_3+p_4-w_4+p_2-r_3-w_3\\&=p_2+p_3+p_4-w_4.\\
         w_3^2&=w_3+\Delta w_3\\&=p_2+p_3+p_4-w_4+w_3.
     \end{align*}
     
     For $w_4\leq0,w_3>0$, so $w_3^2>0$ and job 3 is part of the continuous queue. The impact of job 3 on the total waiting time is:\[\Delta_3=\Delta w_4^+=p_2+p_3+p_4-w_4.\]
     
     The change of job 5's waiting time is:\begin{align*}
         \Delta w_5&=r_3+w_3^2+p_3-r_5-w_5\\&=p_2+p_3+p_4-w_4+w_4+w_5-p_4-w_5\\&=p_2+p_3.
     \end{align*}
     
     The adjustment of total waiting time is:\begin{align*}
         \Delta w&=\Delta_2+\Delta_3+\Delta_4+\pi_5 (\Delta w_5 )\\&=p_2+p_3+p_4-w_4+p_2+p_3+p_4-w_4+\pi_5 (p_2+p_3 )\\&=p_2+p_3+p_4-w_4+p_4-w_4+p_2+p_3+\pi_5 (p_2+p_3 ).
     \end{align*}
     
     \textbf{Case 3:} $w_4\leq0,w_2>0,w_3\leq0$
     
     Before job 4 is moved, jobs 3 and 4 are the queue leader, and job 2 is part of the continuous queue.     
     In this moment, the new waiting time of job 4 and its change of waiting time are respectively:\begin{align*}
         \Delta w_4&=r_0^o+p_0-r_4-w_4\\&=w_3-p_2-p_3,\\
         w_4^2&=w_4+\Delta w_4\\&=w_3+w_4-p_2-p_3<0.
     \end{align*}
     
     So, job 4 is the queue leader, and the change in total waiting time caused by job 4 is:\[\Delta_4=\max(0,w_4^2 )-\max(0,w_4 )=0.\]
     
     The updated information of job 2 is:\begin{align*}
         \Delta w_2&=r_4^o+p_4-r_3-w_2\\&=r_4^o+p_4-(r_2+w_2+p_2-p_2-w_3+w_3+p_3-p_3 )\\&=p_4+p_3+p_2-w_3-w_4.\\
         w_2^2&=w_2+\Delta w_2\\&=p_2+p_3+p_4+w_2-w_3-w_4>0.
     \end{align*}
     
     For job 2 is in the continuous queue, its influence on the total waiting time is:\[\Delta_2=\Delta w_3^+=p_2+p_3+p_4-w_3-w_4.\]
     
     The new information for job 3 is as follows:\begin{align*}
         \Delta w_3&=r_4^o+p_2+p_4-r_3-w_3\\&=r_4^o+p_2+p_4-r_3-w_3+p_3-p_3\\&=p_2+p_3+p_4-w_4-w_3.\\
         w_3^2&=w_3+\Delta w_3\\&=p_2+p_3+p_4-w_4>0.
     \end{align*}
     
     Job 3 is in the continuous queue for $w_3^2>0$, the impact on the total waiting time is:\[\Delta_3=\Delta w_4^+=p_2+p_3+p_4-w_4.\]
     
     The change of job 5's waiting time is:\begin{align*}
         \Delta w_5&=r_4^o+p_4+p_3+p_2-r_5-w_5\\&=p_2+p_3.
     \end{align*}
     
     Finally, the change of total waiting time is:\begin{align*}
         \Delta w &= \Delta_2 + \Delta_3 + \Delta_4 + \pi_5 (\Delta w_5)\\&=p_2 + p_3 + p_4 - w_3 - w_4 + p_2 + p_3 + p_4 - w_4 + \pi_5 (p_2 + p_3)\\&=p_2 + p_3 + p_4 - w_3 - w_4 + p_4 - w_4 + p_2 + p_3 + \pi_5 (p_2 + p_3).
     \end{align*}
     
     \textbf{Case 4:} $w_4\leq0,w_2>0,w_3>0$
     
     Now, job 4 is the queue leader, and both job 2 and 3 are in the continuous queue.
     
     The updated information of job 4 is given by:\begin{align*}
         \Delta w_4&=r_0^o+p_0-r_4-w_4\\&=r_0^o+p_0-r_4-(r_0^o+p_0+p_2+p_3-r_4 )\\&=-p_2+p_3.\\
         w_4^2&=w_4+\Delta w_4\\&=w_4-p_2-p_3<0.
     \end{align*}
     
     Job 4 is the queue leader, its influence on the total waiting time is:\[\Delta_4=\max(0,w_4^2 )-\max(0,w_4 )=0.\]
     
     The new information about job 2 is:\begin{align*}
         \Delta w_2&=r_4^o+p_4-r_2-w_2\\&=r_4^o+p_4-(w_4-w_2-p_2-p_3+r_4 )-w_2\\&=p_2+p_3+p_4-w_4.\\
         w_2^2&=w_2+\Delta w_2\\&=p_2+p_3+p_4-w_4+w_2>0.
     \end{align*}
     
     For job 2 is in the continuous queue, the change in total waiting time due to job 2 is:\begin{align*}
         \Delta_2=\Delta w_3^+=p_2+p_3+p_4-w_4.
     \end{align*}
     
     The updated information of job 3 is:\begin{align*}
         \Delta w_3&=r_4^o+p_4+p_2-r_3-w_3\\&=r_4^o+p_4+p_2-r_3-w_3-p_3+p_3\\&=p_2+p_3+p_4-w_4.\\
         w_3^2&=w_3+\Delta w_3\\&=p_2+p_3+p_4-w_4+w_3>0.
     \end{align*}
     
     The change of job 5's waiting time is:\begin{align*}
         \Delta w_5&=r_4^o+p_2+p_3+p_4-r_5-w_5\\&=p_2+p_3+w_5-w_5\\&=p_2+p_3.
     \end{align*}
     
     The total waiting time's change is:\begin{align*}
         \Delta w &= \Delta_2 + \Delta_3 + \Delta_4 + \pi_5 (\Delta w_5)\\&=p_2 + p_3 + p_4 - w_4 + p_2 + p_3 + p_4 - w_4 + \pi_5 (p_2 + p_3)\\&=p_2 + p_3 + p_4 - w_4 + p_4 - w_4 + p_2 + p_3 + \pi_5 (p_2 + p_3).
     \end{align*}
     
     \textbf{Case 5:} $w_4>0,w_2\leq0,w_3\leq0$
     
     Before job 4's relocation, both jobs 2 and 3 are the queue leader, and job 4 is in  the continuous queue.     
     The updated information about job 4 is:\begin{align*}
         \Delta w_4&=r_0^o+p_0-r_4-w_4\\&=r_0^o+p_0-r_4-(r_0^o+p_0-r_4-w_2+p_2-w_3+p_3 )\\&=w_2+w_3-p_2-p_3.\\
         w_4^2&=w_4+\Delta w_4\\&=w_4+w_3+w_2-p_2-p_3.
     \end{align*}
     
     (1) If $w_4+w_3+w_2-p_2-p_3>0$, job 4 is part of continuous queue.\begin{align*}
         \Delta_4&=\max(0,w_4^2 )-\max(0,w_4 )\\&=w_3+w_2-p_2-p_3.
     \end{align*}
     
     The new information of job 2 is:\begin{align*}
         \Delta w_2&=r_4+w_4^2+p_4-r_2-w_2\\&=r_4+w_2+w_3+w_4-p_2-p_3+p_4-r_2-w_2\\&=p_4+w_4-w_4\\&=p_4.\\
         w_2^2&=w_2+\Delta w_2\\&=w_2+p_4.
     \end{align*}
     
     1) When $w_2+p_4\leq0$, job 2 is the queue leader, the influence of job 2 on the total waiting time is:\[\Delta_2=\Delta w_3^+=0.\]
     The new information of job 3 is:\begin{align*}
         \Delta w_3&=r_2^o+p_2-r_3-w_3\\&=w_3-w_3=0.\\
         w_3^2&=w_3+\Delta w_3\\&=w_3\leq 0.
     \end{align*}
     
     Then, job 3 is the queue leader, the effect of job 3 on the total waiting time is:\[\Delta_3=\Delta w_4^+=0.\]
     
     The change of job 5's waiting time is:\begin{align*}
         \Delta w_5&=r_3^o+p_3-r_5-w_5\\&=r_3^o+p_3+p_4-p_4-r_5-w_5\\&=-p_4.
     \end{align*}
     
     The adjustment of total waiting time is:\begin{align*}
         \Delta w &= \Delta_2 + \Delta_3 + \Delta_4 + \pi_5 (\Delta w_5)\\&=w_3 + w_2 - p_2 - p_3 + \pi_5 (-p_4)\\&=p_4 + w_3 + w_2 - p_2 - p_3 - p_4 + \pi_5 (-p_4).
     \end{align*}
     
     2) When $w_2+p_4>0$, job 2 belongs to the continuous queue, its influence on the total waiting time is:\begin{align*}
         \Delta_2=\Delta w_3^+=w_2+p_4.
     \end{align*}
     
     The updated information about job 3 is:\begin{align*}
         \Delta w_3&=r_2+w_2^2+p_2-r_3-w_3\\&=w_2^2\\&=w_2+p_4.
         \\w_3^2&=w_3+\Delta w_3\\&=w_2+w_3+p_4.
     \end{align*}
     
    \romannumeral 1\text{ }) When $w_2+w_3+p_4\leq0$, job 3 is the queue leader, and its affect on the total waiting time can be listed:\[\Delta_3=\Delta w_4^+=0.\]
    The change of job 5 is:\begin{align*}
        \Delta w_5&=r_3^o+p_3-r_5-w_5\\&=r_3^o+p_3+p_4-p_4-r_5-w_5\\&=-p_4.
    \end{align*}
    
    The change of the total waiting time is:\begin{align*}
        \Delta w &= \Delta_2 + \Delta_3 + \Delta_4 + \pi_5 (\Delta w_5)\\&=p_4 + p_4 + w_2 + w_2 + w_3 - p_2 - p_3 - p_4 + \pi_5 (-p_4).
    \end{align*}
    
    \romannumeral 2\text{ }) When $w_2+w_3+p_4>0$, job 3 is one of the members of the continuous queue, and the change in total waiting time caused by job 3 is:\[\Delta_3=\Delta w_4^+=w_2+w_3+p_4.\]
    
    The change of job 5's waiting time is:\begin{align*}
        \Delta w_5&=r_3+w_3^2+p_3-r_5-w_5\\&=w_3^2-p_4\\&=w_2+w_3.
    \end{align*}
    
    Finally, the change of total waiting time is as follows:\begin{align*}
        \Delta w &= \Delta_2 + \Delta_3 + \Delta_4 + \pi_5 (\Delta w_5)\\&=p_4 + p_4 + w_2 - p_2 + w_3 - p_3 + w_2 + w_2 + w_3 + \pi_5 (w_2 + w_3).
    \end{align*}
    
    (2) If $w_4+w_3+w_2-p_2-p_3\leq0$, job 4 is the queue leader, the impact of job 4 on the total waiting time is:\[\Delta_4=\max(0,w_4^2 )-\max(0,w_4 )=-w_4.\]
    
    The updated information about job 2 is:\begin{align*}
        \Delta w_2&=r_4^o+p_4-r_2-w_2\\&=r_4^o+p_4-r_2-p_2+p_2-w_3+w_3-p_3+p_3-w_2\\&=p_2+p_3+p_4-w_2-w_3-w_4.\\
        w_2^2&=w_2+\Delta w_2\\&=p_2+p_3+p_4-w_3-w_4.
    \end{align*}
    
    1) When $p_2+p_3+p_4-w_3-w_4>0$, job 2 is part of continuous queue, its influence on the total waiting time is:\begin{align*}
        \Delta_2=\Delta w_3^+=p_2+p_3+p_4-w_3-w_4.
    \end{align*}
    
    The information of job 3, after job 4's relocation, is as follows:\begin{align*}
        \Delta w_3&=r_2+w_2^2+p_2-r_3-w_3\\&=p_2+p_3+p_4-w_3-w_4.\\
        w_3^2&=w_3+\Delta w_3\\&=p_2+p_3+p_4-w_4.
    \end{align*}
    
    \romannumeral 1 \text{ }) When $p_2+p_3+p_4-w_4\leq0$, job 3 is the queue leader, and its impact on the total waiting time is:\begin{align*}
        \Delta_3=\Delta w_4^+=0.
    \end{align*}
    
    The change of job 5's waiting time is:\begin{align*}
        \Delta w_5&=r_3^o+p_3-r_5-w_5\\&=r_3^o+p_3+p_4-p_4-r_5-w_5\\&=-p_4.
    \end{align*}
    
    In the end, the change of total waiting time is:\begin{align*}
        \Delta w &= \Delta_2 + \Delta_3 + \Delta_4 + \pi_5 (\Delta w_5)\\&=p_2 + p_3 + p_4 - w_3 - w_4 - w_4 + \pi_5 (-p_4)\\&= p_3 + p_2 + p_4 - w_3 - w_4 + p_4 - w_4 - p_4 + \pi_5 (-p_4).
    \end{align*} 
    
    \romannumeral 2 \text{ }) When $p_2+p_3+p_4-w_4>0$, job 3 belongs to the continuous queue, its impact on the total waiting time is:\[\Delta_3=\Delta w_4^+=p_2+p_3+p_4-w_4.\]
    
    The change of job 5's waiting time is:\begin{align*}
        \Delta w_5&=r_3+w_3^2+p_3-r_5-w_5\\&=r_3+w_3^2+p_3+p_4-p_4-r_5-w_5\\&=w_3^2-p_4\\&=p_2+p_3-w_4.
    \end{align*}
    
    The change of total waiting time is:\begin{align*}
        \Delta w &= \Delta_2 + \Delta_3 + \Delta_4 + \pi_5 (\Delta w_5)\\& = p_2 + p_3 + p_4 - w_3 - w_4 + p_4 - w_4 + p_2 + p_3 - w_4 + \pi_5 (p_2 + p_3 - w_4).
    \end{align*}

    2) When $p_2+p_3+p_4-w_3-w_4\leq0$, job 2 is the queue leader, and its impact on the total waiting time is:\begin{align*}
        \Delta_2=\Delta w_3^+=0.
    \end{align*}
    
    The new information about job 3, after job 4's moved, is as follows:\begin{align*}
        \Delta w_3&=r_2^o+p_2-r_3-w_3=0,\\w_3^2&=w_3+\Delta w_3=w_3.
    \end{align*}
    
    For $w_3\leq0$, job 3 still be the queue leader, and its impact on the total waiting time is:\begin{align*}
        \Delta_3=\Delta w_4^+=0.
    \end{align*}
    
    The change of total waiting time is:\begin{align*}
        \Delta w_5&=r_3^o+p_3-r_5-w_5\\&=r_3^o+p_3+p_4-p_4-r_5-w_5\\&=-p_4.
    \end{align*}
    
    The adjustment of total waiting time is as follows:\begin{align*}
        \Delta w &= \Delta_2 + \Delta_3 + \Delta_4 + \pi_5 (\Delta w_5)\\&= -w_4 + \pi_5 (-p_4)\\&= p_4 - w_4 - p_4 + \pi_5 (-p_4).
    \end{align*}
    
     \textbf{Case 6:} $w_4>0,w_2\leq0,w_3>0$
     
     Before job 4's movement, job 2 is the queue leader, both job 3 and 4 are in the continuous queue.
     
     The change of job 4's waiting time and its new waiting time are, respectively:\begin{align*}
         \Delta w_4&=r_0^o+p_0-r_4-w_4\\&=r_0^o+p_0-r_4-(r_0^o+p_0-w_2+p_2+p_3-r_4 )\\&=w_2-p_2-p_3.\\
         w_4^2&=w_4+\Delta w_4\\&=w_2+w_4-p_2-p_3.
     \end{align*}
     
     (1) If $w_2+w_4-p_2-p_3>0$, job 4 is part of the continuous queue.\begin{align*}
         \Delta_4&=\max(0,w_4^2 )-\max(0,w_4 )\\&=w_2-p_2-p_3.
     \end{align*}
     
     The new information of job 2 is:\begin{align*}
         \Delta w_2&=r_4+w_4^2+p_4-r_2-w_2\\&=w_2+w_4+p_4-p_2-p_3+r_4-r_2-p_2+p_2-p_3+p_3-w_2\\&=w_4+p_4-p_2-p_3+p_2+p_3-w_4\\&=p_4.\\
         w_2^2&=w_2+\Delta w_2\\&=w_2+p_4.
     \end{align*}
     
     1) When $w_2+p_4>0$, job 2 belongs to the continuous queue, and its influence on the total waiting time is:\begin{align*}
         \Delta_2=\Delta w_3^+=w_2+p_4.
     \end{align*}
     
     The updated information of job 3 is:\begin{align*}
         \Delta w_3&=r_2+w_2^2+p_2-r_3-w_3\\&=r_2+w_2+p_4+p_2-r_3-w_3\\&=w_2+p_4.\\
         w_3^2&=w_3+\Delta w_3\\&=w_3+w_2+p_4>0.
     \end{align*}
     
     For $w_3^2>0$, job 3 is one of the members of the continuous queue, and its impact on the total waiting time is:\begin{align*}
         \Delta_3=\Delta w_4^+=w_2+p_4.
     \end{align*}
     
     The change of job 5's waiting time is:\begin{align*}
         \Delta w_5&=r_3+w_3^2+p_3-r_5-w_5\\&=r_3+w_2+p_4+w_3+p_3-r_5-(r_3+w_3+p_3+p_4-r_5 )\\&=w_2.
     \end{align*}
     
     At last, the total waiting time's change is:\begin{align*}
         \Delta w &= \Delta_2 + \Delta_3 + \Delta_4 + \pi_5 (\Delta w_5)\\&= w_2 - p_2 - p_3 + w_2 + p_4 + p_4 + w_2 + \pi_5 (w_2).
     \end{align*}
     
     2) When $w_2+p_4\leq0$, job 2 is the queue leader, and its affect on the total waiting time is:\[\Delta_2=\Delta w_3^+=0.\]
     
     The new information about job 3 is:\begin{align*}
         \Delta w_3&=r_2^o+p_2-r_3-w_3\\&=0.\\
         w_3^2&=w_3+\Delta w_3\\&=w_3.
     \end{align*}
     
     The total waiting time's change is:\begin{align*}
         \Delta w &= \Delta_2 + \Delta_3 + \Delta_4 + \pi_5 (\Delta w_5)\\&= w_2 - p_2 - p_3 + p_4 - p_4 + \pi_5 (-p_4).
     \end{align*}
     
     (2) If $w_2+w_4-p_2-p_3\leq0$, job 4 is the queue leader, we can calculate that:\[\Delta_4=\max(0,w_4^2 )-\max(0,w_4 )=0.\]
     The updated information of job 2 is:\begin{align*}
         \Delta w_2&=r_4^o+p_4-r_2-w_2\\&=r_4^o+p_4-(w_4-p_2-p_3+r_4 )-w_2\\&=p_2+p_3+p_4-w_4-w_2.\\
         w_2^2&=w_2+\Delta w_2\\&=p_2+p_3+p_4-w_4.
     \end{align*}
     
     1) When $p_2+p_3+p_4-w_4>0$, job 2 is part of continuous queue, and its impact on the total waiting time is:\begin{align*}
         \Delta_2=\Delta w_3^+=p_2+p_3+p_4-w_4.
     \end{align*}
     
     The new information of job 3 is:\begin{align*}
         \Delta w_3&=r_2+p_2+w_2^2-r_3-w_3\\&=p_2+p_2+p_3+p_4-w_4+r_2-r_3-w_3\\&=p_2+p_3+p_4-w_4.\\
         w_3^2&=w_3+\Delta w_3\\&=p_2+p_3+p_4-w_4+w_3>0.
     \end{align*}
     
     Now, job 3 belongs to the continuous queue, and its influence on the total waiting time is:\[\Delta_3=\Delta w_4^+=p_2+p_3+p_4-w_4.\]
     
     The change of job 5's waiting time is:\begin{align*}
         \Delta w_5&=r_3+w_3^2+p_3-r_5-w_5\\&=r_3+p_2+p_3+p_4-w_4+w_3+p_3-r_5-(r_3+w_3+p_3+p_4-r_5 )\\&=p_2+p_3-w_4.
     \end{align*}
     
     Finally, the change of total waiting time is:\begin{align*}
         \Delta w &= \Delta_2 + \Delta_3 + \Delta_4 + \pi_5 (\Delta w_5)\\&= p_2 + p_3 + p_4 - w_4 + p_2 + p_3 + p_4 - w_4 - w_4 + \pi_5 (p_2 + p_3 - w_4)\\&= p_4 - w_4 + p_2 + p_3 + p_4 - w_4 + p_2 + p_3 - w_4 + \pi_5 (p_2 + p_3 - w_4).
     \end{align*}
     
     2) When $p_2+p_3+p_4-w_4\leq0$, job 2 is the queue leader, and its effect on the total waiting time is:\begin{align*}
         \Delta_2=\Delta w_3^+=0.
     \end{align*}
     
     The updated information about job 3 is:\begin{align*}
         \Delta w_3&=r_2^o+p_2-r_3-w_3=0,\\w_3^2&=w_3>0
     \end{align*}
     
     For $w_3^2>0$, the impact of job 3 on the total waiting time is:\[\Delta_3=\Delta w_4^+=0.\]
     
     The change of job 5's waiting time is:\begin{align*}
         \Delta w_5&=r_3+w_3+p_3-r_3-w_5\\&=r_3+w_3+p_3-r_5-(r_3+w_3+p_3+p_4-r_5 )\\&=-p_4.
     \end{align*}
     
     The adjustment of total waiting time is:\begin{align*}
         \Delta w &= \Delta_2 + \Delta_3 + \Delta_4 + \pi_5 (\Delta w_5)\\&= -w_4 + \pi_5 (-p_4)\\&=p_4 - w_4 - p_4 + \pi_5 (-p_4).
    \end{align*}
    
     \textbf{Case 7:} $w_4>0,w_2>0,w_3\leq0$
     
     Now, job 3 is the queue leader, both jobs 2 and 4 are in the continuous queue.  
     The change of job 4's waiting time and its new waiting time are respectively:\begin{align*}
         \Delta w_4&=r_0^o+p_0-r_4-w_4\\&=r_0^o+p_0-r_4-(r_0^o+p_0+p_2-w_3+p_3-r_4 )\\&=w_3-p_2-p_3.\\
         w_4^2&=w_4+\Delta w_4\\&=w_4+w_3-p_2-p_3.
     \end{align*}
     
     (1) When $w_4+w_3-p_2-p_3>0$, job 4 is part of continuous queue, and the adjustment of total waiting time due to job 4 is:\[\Delta_4=\max(0,w_4^2 )-\max(0,w_4 )\\=w_3-p_2-p_3.\]
     The updated information about job 2 is:\begin{align*}
         \Delta w_2&=r_4+w_4^2+p_4-r_2-w_2\\&=w_4+w_3-p_2-p_3+r_4+p_4-r_2-w_2\\&=w_4+w_3-p_2-p_3+r_2+p_2+w_2-w_3+p_3-w_4+p_4-r_2-w_2\\&=p_4.\\
         w_2^2&=w_2+\Delta w_2\\&=p_4+w_2>0.
     \end{align*}
     
     For $w_2^2>0$, job 2 belongs to the continuous queue, the influence of job 2 on the total waiting time is:\[\Delta_2=\Delta w_3^+=p_4.\]
     
     The updated information of job 3 is:\begin{align*}
         \Delta w_3&=r_2+w_2^2+p_2-r_3-w_3\\&=r_2+w_2+p_4+p_2-r_3-w_3\\&=p_4.\\
         w_3^2&=w_3+\Delta w_3\\&=p_4+w_3.
     \end{align*}
     
     1) If $p_4+w_3\leq0$, at this point, job 3 is the queue leader, and its impact on the total waiting time is:\[\Delta_3=\Delta w_4^+=0.\]
     
     The change of waiting time of job 5 is:\begin{align*}
         \Delta w_5=r_3^o+p_3-r_5-w_5=-p_4.
     \end{align*}
     
     Finally, the change of total waiting time is:\begin{align*}
         \Delta w &= \Delta_2 + \Delta_3 + \Delta_4 + \pi_5 (\Delta w_5)\\&=w_3 - p_2 - p_3 + p_4 + p_4 - p_4 + \pi_5 (-p_4)    \end{align*}
         
    2) If $p_4+w_3>0$, job 3 is in the continuous queue, the impact of job 3 on the total waiting time is:\[\Delta_3=\Delta w_4^+=w_3+p_4.\]
    
    The change of job 5's waiting time is:\begin{align*}
        w_5&=r_3+w_3^2+p_3-r_5-w_5\\&=p_4+w_3-p_4\\&=w_3.
    \end{align*}
    
    The change of total waiting time is as follows:\begin{align*}
        \Delta w &= \Delta_2 + \Delta_3 + \Delta_4 + \pi_5 (\Delta w_5)\\&=p_4 + w_3 - p_2 - p_3 + p_4 + w_3 + \pi_5 (w_3).
    \end{align*}
    
    (2) When $w_4+w_3-p_2-p_3\leq \leq0$, job 4 is the queue leader, and its impact on the total waiting time is:\[\Delta_4=\max(0,w_4^2 )-\max(0,w_4 )=-w_4.\]
    
    The updated waiting time is :\begin{align*}
        \Delta w_2&=r_4^o+p_4-r_2-w_2\\&=(r_2+p_2+w_2-w_3+p_3-w_4 )+p_4-r_2-w_2\\&=p_2+p_3+p_4-w_3-w_4.\\
        w_2^2&=w_2+\Delta w_2\\&=p_2+p_3+p_4-w_3-w_4+w_2
    \end{align*}
    
    Because $w_2>0$ and $p_2+p_3>w_3+w_4$,we can get that $w_2^2>0$. Thus, job 2 is in the continuous queue,and the change of total waiting time caused by job 2 is:\begin{align*}
        \Delta_2=\Delta w_3^+=p_2+p_3+p_4-w_3-w_4.
    \end{align*}
    
    The new information of job 3, once job 4 is moved, is as follows:\begin{align*}
        \Delta w_3&=r_2+w_2^2+p_2-r_3-w_3\\&=r_2+p_2+p_3+p_4-w_3-w_4+w_2+p_2-r_3-(r_2+w_3+p_2-r_3 )\\&=p_2+p_3+p_4-w_3-w_4.\\
        w_3^2&=w_3+\Delta w_3\\&=p_2+p_3+p_4-w_4.
    \end{align*}
    
    1) If $p_2+p_3+p_4-w_4>0$, job 3 is one of the members of the continuous queue, and its influence on the total waiting time is:\begin{align*}
        \Delta_3=\Delta w_4^+=p_2+p_3+p_4-w_4-w_3.
    \end{align*}
    
    The change of job 5's waiting time is:\begin{align*}
        \Delta w_5&=r_3+w_3^2+p_3-r_5-w_5\\&=r_3+p_2+p_3+p_4-w_4+p_3-r_5-(r_3^o+p_3+p_4-r_5 )\\&=p_2+p_3-w_4.
    \end{align*}
    
    Finally, the adjustment of total waiting time is as follows:\begin{align*}
        \Delta w &= \Delta_2 + \Delta_3 + \Delta_4 + \pi_5 (\Delta w_5)\\&=2p_4 - 2w_4 + p_2 + p_3 - w_3  + p_2 + p_3 - w_4 + \pi_5 (p_2 + p_3 - w_4).
    \end{align*}
    
    2) If $p_2+p_3+p_4-w_4\leq0$, job 3 is the queue leader, and its impact on the total waiting time is:\[\Delta_3=\Delta w_4^+=0.\]
    
    The change of job 5's waiting time is:\begin{align*}
        \Delta w_5&=r_3^o+p_3-r_5-w_5\\&=r_3^o+p_3-r_5-(r_3^o+p_3+p_4-r_5 )\\&=-p_4.
    \end{align*}
    
    Finally, the adjustment of total waiting time caused by job 4's relocation is calculated as:\begin{align*}
        \Delta w &= \Delta_2 + \Delta_3 + \Delta_4 + \pi_5 (\Delta w_5)\\&=p_4 - w_4 + p_2 + p_3 + p_4 - w_3 - w_4 - p_4 + \pi_5 (-p_4).
    \end{align*}
        
     \textbf{Case 8:} $w_4>0,w_2>0,w_3>0$
     
     When not moved, job 2, 3 and 4 are belong to the continuous queue.     
     The updated information about job 4, after its relocation, we can get:\begin{align*}
         \Delta w_4&=r_0^o+p_0-r_4-w_4\\&=r_0^o+p_0-r_4-(r_0^o+p_0+p_2+p_3-r_4 )\\&=-p_2-p_3.\\w_4^2&=w_4+\Delta w_4\\&=w_4-p_2-p_3.
     \end{align*}
     
     (1) If $w_4-p_2-p_3>0$, it means that job 4 is in the continuous queue, and the affect of job 4 on the total waiting time is:\[\Delta_4=\max(0,w_4^2 )-\max(0,w_4 )=-p_2-p_3.\]
     
     The new information of job 2 is:\begin{align*}
         \Delta w_2&=r_4+w_4^2+p_4-r_2-w_2\\&=r_4+w_4-p_2-p_3+p_4-r_2-w_2\\&=p_4.\\
         w_2^2&=w_2+\Delta w_2\\&=w_2+p_4.
     \end{align*}
     
     For $w_2>0$, we can get $w_4^2>0$. Job 2 in the continuous queue, the impact of job 2 on the total waiting time is:\[\Delta_2=\Delta w_3^+=p_4.\]
     
     The change of job 3's waiting time and its new waiting time are, respectively:\begin{align*}
         \Delta w_3&=r_2+w_2^2+p_2-r_3-w_3\\&=r_2+w_2+p_4+p_2-r_3-w_3\\&=p_4.\\
         w_3^2&=w_3+\Delta w_3\\&=w_3+p_4.
     \end{align*}
     
     For $w_3>0$, so $w_3^2>0$, job 3 is in the continuous queue, the affect of job 3 on the total waiting time is:\[\Delta_3=\Delta w_4^+=p_4.\]
     
     The change of job 5's waiting time is:\begin{align*}
         \Delta w_5&=r_3+w_3^2+p_3-r_5-w_5\\&=r_3+w_3+p_4+p_3-r_5-w_5\\&=0.
     \end{align*}
     
     The change of total waiting time is:\begin{align*}
         \Delta w &= \Delta_2 + \Delta_3 + \Delta_4 + \pi_5 (\Delta w_5)\\&=p_4 + p_4 - p_2 - p_3 + \pi_5 (0).         
     \end{align*}
     
     (2) If $w_4-p_2-p_3\leq0$, job 4 is the queue leader, and its influence on the total waiting time is:\[\Delta_4=\max(0,w_4^2 )-\max(0,w_4 )=-w_4.\]
     
     The updated information of job 2 is:\begin{align*}
         \Delta w_2&=r_4+p_4-r_2-w_2\\&=p_2+p_3+p_4-w_4.\\
         w_2^2&=w_2+\Delta w_2\\&=p_2+p_3+p_4-w_4+w_2.
     \end{align*}
     
     For $w_4\leq p_2-p_3$ and $w_2>0$, so $w_2^2>0$. At this point, job 2 is in the continuous queue, and its impact on the total waiting time is:\[\Delta_2=\Delta w_3^+=p_2+p_3+p_4-w_4.\]
     
     The new information of job 3 is:\begin{align*}
         \Delta w_3&=r_2+w_2^2+p_2-r_3-w_3\\&=r_2+p_2+p_3+p_4-w_4+w_2+p_2-r_3-w_3\\&=p_2+p_3+p_4-w_4.\\
         w_3^2&=w_3+\Delta w_3\\&=p_2+p_3+p_4-w_4+w_3>0.
     \end{align*}
     
     Similar to job 2, job 3 belongs to the continuous queue, and the total waiting time change due to job 3 is:\[\Delta_3=\Delta w_4^+=p_2+p_3+p_4-w_4.\]
     
     The change of job 5's waiting time is:\begin{align*}
         \Delta w_5&=r_3+w_3^2+p_3-r_5-w_5\\&=r_3+p_2+p_3+p_4-w_4+w_3+p_3-r_5-w_5\\&=p_2+p_3-w_4.
     \end{align*}
     
     Finally, the change of total waiting time is:\begin{align*}
         \Delta w &= \Delta_2 + \Delta_3 + \Delta_4 + \pi_5 (\Delta w_5)\\&= p_4 - w_4 + p_4 + p_2 + p_3 - w_4 + p_2 + p_3 - w_4 + \pi_5 (p_2 + p_3 - w_4).
     \end{align*}

     Now, let us integrate the conditions for $I_4$. $I_4$ can be expressed as:\[I_4=\max(-p_4,\min(0,w_2 )+\min(0,w_3),p_2+p_3-\max(0,w_4))\]

     In the next step, we get the formula of $I_3$ is:\begin{align*}
         I_3=&\min(0,w_2 )+\min(0,w_3 )-p_2-p_3\\&+\max(0,\min(0,w_2 )+p_4,p_2+p_3+p_4-\min(0,w_3 )-w_4 )\\&+\max(p_4,p_2+p_3+p_4-\min(0,w_2 )-\min(0,w_3 )-w_4).
     \end{align*}
     
     Currently, only the  information for $w_4^2$ is miss, we summarize the $w_4^2$ as:
     \begin{align*}
         w_4^2=\min(0,w_2)+\min(0,w_3)-p_2-p_3 +w_4.
     \end{align*}

     \textbf{Conclusion:} By summarizing all the above information, we get the law that governing the change of queue when job 4 is relocated the position between jobs 1 and 2. The change of total waiting follows this formula:\[\Delta w=I_3+I_4+\pi_5(I_4).\]
     where\begin{align*}
         I_3=&\min(0,w_2 )+\min(0,w_3 )-p_2-p_3\\&+\max(0,\min(0,w_2 )+p_4,p_2+p_3+p_4-\min(0,w_3 )-w_4 )\\&+\max(p_4,p_2+p_3+p_4-\min(0,w_2 )-\min(0,w_3 )-w_4 ).
     \end{align*}
     \begin{align*}
         I_4=\max(-p_4,\min(0,w_2 )+\min(0,w_3),p_2+p_3-\max(0,w_4)).
     \end{align*}
     It is obviously that  $\max(0,\min(0,w_2 )+p_4,p_2+p_3+p_4-\min(0,w_3 )-w_4 )$ and $\max(p_4,p_2+p_3+p_4-\min(0,w_2 )-\min(0,w_3 )-w_4 )$ are in line with recursive formula of Lemma 2. So, the equivalent formula of $I_3$ can be written as:\begin{align*}
         I_3&=\min(0,w_2 )+\min(0,w_3 )-p_2-p_3+\sum_{j=2}^3 \Delta w_j^+ \\&=\sum_{j=2}^3 (\min(0,w_j )-p_j) +\sum_{j=2}^3 \Delta w_j^+.
     \end{align*} where the initial value is $\Delta w_2^+=\max(p_4,p_2+p_3+p_4-\min(0,w_2 )-\min(0,w_3 )-w_4 )$.
      The new waiting time of $w_4^2$ can be calculated as:\[w_4^2=\min(0,w_2)+\min(0,w_3)-p_2-p_3 +w_4.\]

      This law can be applied to the backward movement of job $\sigma_{t} (i)$ in any sequence $\sigma_{t}$ to the position before $\sigma_{t}(k)$(where $k\in [1,i-1]$). The change of total waiting time can be expressed as:\[\Delta w=I_3+I_4+\pi_{\sigma_t (k+1)} (I_4).\] where \begin{align*}
          I_3&=\sum_{j=k}^{i-1} (\min(0,w_{\sigma_t (j)} )-p_{\sigma_t (j)})+\sum_{j=k}^{i-1} \Delta w_{\sigma_t (j)}^+.\\
          I_4&=\max \left(-p_{\sigma_t (i)},\sum_{j=k}^{i-1} \min(0,w_{\sigma_t (j)} ) ,\sum_{j=k}^{i-1} p_{\sigma_t (j)} -\max(0,w_{\sigma_t (i)}) \right).
      \end{align*} 
      where the initial value of the recursive formula is:\begin{align*}
          \Delta w_{\sigma_t (k)}^+=\max \left(p_{\sigma_t (i)},\sum_{j=k}^i p_{\sigma_t (j)} -\sum_{j=k}^{i-1} \min(0,w_{\sigma_t (j)} ) -w_{\sigma_t (i)} \right).
      \end{align*}

    The job $\sigma_{t} (i)$'s new waiting time $w_{\sigma_{t} (i)}^2$ is calculated as:\[w_{\sigma_t (i)}^2=\sum_{j=k}^{i-1}(\min(0,w_{\sigma_t (j)} ) -p_{\sigma_t (j)}) +w_{\sigma_t (i)}.\]
\end{proof}

\section{Omitted Proofs in Section 3}
\subsection{ Proof of \cref{forward condition}}\label{app:forward condition}
\begin{proof}
\cref{Theorem1} provides the formula for the variation in the total waiting time caused by a forward move of a job. Any forward move yielding $\Delta w \leq 0$ is regarded as an improving solution. Once the expression for $\Delta w$ has been derived, it can be further analyzed according to the two key terms $I_1$ and $I_2$. This leads to the following four mutually exclusive cases: 

Case 1: $I_1 \leq 0,\ I_2 \geq 0$;

Case 2: $I_1 > 0,\ I_2 > 0$;

Case 3: $I_1 \leq 0,\ I_2 \leq 0$;

Case 4: $I_1 > 0,\ I_2 < 0$.

Consider an arbitrary schedule $\sigma_t$ and any four consecutive jobs $\sigma_t(i)$, $\sigma_t(i+1)$, and $\sigma_t(i+2)$, where $i\in\{1,\dots,n-2\}$. Without loss of generality, relabel them as $4,5,6$, respectively. We analyze the move in which job $4$ is moved forward to the position immediately after job $6$.

Case 2 can never reduce the total waiting time and can therefore be excluded directly. Case 4 requires $I_2<0$. A necessary condition for $I_2<0$ is
\[
w_4 < 0, \quad w_5>0,\ w_6>0, \quad p_4 - w_5<0,\ p_4 - w_6<0.
\]

Case 2 is clearly non-improving and can be excluded immediately. For case 4, the condition $I_2<0$ is required, which in turn implies the necessary conditions
\[
w_4<0,\qquad w_5>0,\ w_6>0,\qquad p_4-w_5<0,\ p_4-w_6<0.
\]
If the original sequence follows the FCFS order, then $r_4\leq r_5$. By \cref{Proposition1},
\[
\max(0,w_4)+p_4 \geq \max(0,w_5).
\]
Since $w_4<0$ and $w_5>0$, this reduces to
\[
p_4-w_5\geq 0,
\]
contradicting $p_4-w_5<0$. Hence, under the FCFS order, case (iv) cannot occur.

If, after some adjustments, the release-time order is reversed so that $r_4>r_5$, then \cref{Proposition1} yields
\[
\max(0,w_4)+p_4<\max(0,w_5),
\]
and thus $p_4<w_5$. In this situation, the queue is discontinuous at job $4$ because $w_4<0$. Such a case is handled by the Bottleneck Breakthrough Rule, and therefore need not be considered in the present analysis of forward candidate solutions.

Cases 1 and 3 can be unified into the single situation $I_1 \leq 0$. By \cref{Lemma1}, the upper bound of $\Delta w_5^-$ is determined by its initial value:
\[
\Delta w_5^- \leq p_4 - \min(0, w_4).
\]
Consider the most favorable upper bound for the propagation, namely, the extreme case in which $\Delta w_{\sigma_t(j)}^- = p_4 - \min(w_4,0)$ for all relevant $j$, corresponding to a situation where the decrease flow does not attenuate during propagation. This yields the most conservative estimate of $I_1$:
\begin{align*}
I_1
&= \sum_{j=5}^{6} \left(p_{\sigma_t(j)} - \min(0, w_{\sigma_t(j)})\right)
   - \sum_{j=6}^{7} \Delta w_{\sigma_t(j)}^- \\
&\geq \sum_{j=5}^{6} p_{\sigma_t(j)} - \sum_{j=6}^{7} \Delta w_{\sigma_t(j)}^-.
\end{align*}
Combining $\Delta w_5^- \leq p_4 - \min(0,w_4 )$ with the monotonicity in \cref{Lemma1}, we further obtain $\Delta w_7^- \leq \Delta w_5^-$. Therefore,
\begin{align*}
I_1
&\geq \sum_{j=5}^{6} p_{\sigma_t(j)} - \sum_{j=6}^{7} \Delta w_{\sigma_t(j)}^- \\
&\geq \sum_{j=5}^{6} \left(p_{\sigma_t(j)} - \Delta w_{\sigma_t(j)}^- \right) \\
&\geq \sum_{j=5}^{6} \left( p_{\sigma_t(j)} - \bigl(p_4 - \min(0,w_4 )\bigr) \right).
\end{align*}

It follows that this condition is necessary for a forward move to reduce the total waiting time. If job $\sigma_t(i)$ is moved forward to position $k$ and satisfies
\[
\sum_{j=i+1}^{k} \left( p_{\sigma_t(j)} - \bigl(p_{\sigma_t(i)} - \min(w_{0},\sigma_t(i) )\bigr) \right) \leq 0,
\]
then such a move may reduce the total waiting time. We call position $k$ a \emph{forward candidate solution} of job $\sigma_t(i)$, denoted by $\theta_{i,k}^f$.
\end{proof}

\subsection{ Proof of \cref{backward condition}}\label{app:backward condition}
\begin{proof}
To complete the structural characterization of the solution space, we further analyze the improvement mechanism induced by a backward move of a job. Consider an arbitrary sequence $\sigma_t$, and select any five consecutive jobs
\[
\sigma_t(i),\ \sigma_t(i+1),\ \sigma_t(i+2),\ \sigma_t(i+3),\ \sigma_t(i+4),
\qquad i\in\{1,\dots,n-4\}.
\]
For notational simplicity and without loss of generality, let
\[
\sigma_t(i)=2,\quad \sigma_t(i+1)=3,\quad \sigma_t(i+2)=4,\quad \sigma_t(i+3)=5,\quad \sigma_t(i+4)=6.
\]
This relabeling is introduced only for local analysis and does not affect the general validity of the argument.

Now consider moving job $4$ backward to the position immediately before job $2$. After this move, the waiting times of jobs $2$ and $3$ may increase, whereas job $4$ is processed earlier and its waiting time may decrease. Meanwhile, the waiting-time variations of job $5$ and all subsequent jobs are determined by the propagation effect. Let $\Delta_4\ge 0$ denote the improvement contributed by job $4$, that is, the change in its contribution to the objective value is $-\Delta_4$. Let $\pi_5(\Delta w_5)$ denote the propagation effect starting from job $5$. Then the total variation in the total waiting time caused by this backward move can be written as
\[
\Delta w=\Delta_2+\Delta_3+\pi_5(\Delta w_5)-\Delta_4,
\]
where $\Delta_2,\Delta_3\ge 0$ represent the nonnegative deterioration contributions of jobs $2$ and $3$, respectively.

We first show that
\[
-\Delta_4+\Delta w_5\le 0
\]
is a necessary condition for this backward move to improve the total waiting time.

\medskip
\noindent
\textbf{Step 1: Proof of the necessary condition.}

If $\Delta w_5\le 0$, then by the propagation rule in \cref{Lemma2}, the propagation contribution generated by job $5$ and all subsequent jobs satisfies
\[
\pi_5(\Delta w_5)\le 0.
\]
In this case, the downstream jobs cannot offset the improvement brought by job $4$, and therefore
\[
-\Delta_4+\Delta w_5\le 0
\]
cannot be violated. Hence, the necessary condition holds trivially in this case.

We next consider the case $\Delta w_5>0$. In this case, the propagation effect starting from job $5$ becomes a nonnegative deterioration term. By \cref{Lemma2},
\[
\Delta_3\ge 0,\qquad \pi_5(\Delta w_5)\ge 0,
\]
and thus
\[
\Delta w=\Delta_2+\Delta_3+\pi_5(\Delta w_5)-\Delta_4
\ge \Delta_2-\Delta_4.
\]
It therefore remains to derive a lower bound for $\Delta_2-\Delta_4$.

By \cref{Theorem2}, when $\Delta w_5>0$,
\[
\Delta w_5=I_4=p_2+p_3-\max(0,w_4)>0.
\]
Hence,
\[
p_2+p_3>\max(0,w_4)\ge w_4.
\]
The increase in the waiting time of job $2$ is
\begin{align*}
\Delta w_2^+
&=\max\left\{p_4,\ \sum_{j=2}^{4}p_j-\sum_{j=2}^{3}\min(0,w_j)-w_4\right\}.
\end{align*}
Since
\[
\sum_{j=2}^{4}p_j-\sum_{j=2}^{3}\min(0,w_j)-w_4
\ge p_4+\bigl(p_2+p_3-w_4\bigr)\ge p_4,
\]
the maximum is attained by the second term, and therefore
\[
\Delta w_2^+=\sum_{j=2}^{4}p_j-\sum_{j=2}^{3}\min(0,w_j)-w_4.
\]
Accordingly, the contribution of job $2$ to the total waiting time is
\begin{align*}
\Delta_2
&=\max\bigl(0,\Delta w_2^++\min(0,w_2)\bigr)\\
&=\sum_{j=2}^{4}p_j-\min(0,w_3)-w_4.
\end{align*}

On the other hand, the new waiting time of job $4$ after the move is
\[
w_4^2=w_4+\min(0,w_2)+\min(0,w_3)-p_2-p_3.
\]
From
\[
p_2+p_3-\max(0,w_4)>0
\]
it follows that
\[
p_2+p_3-w_4>0.
\]
Combining this with
\[
-\min(0,w_2)\ge 0,\qquad -\min(0,w_3)\ge 0,
\]
we obtain
\[
p_2+p_3-\min(0,w_2)-\min(0,w_3)-w_4>0,
\]
that is,
\[
w_4^2<0.
\]
Therefore, the waiting time of job $4$ after the move remains nonpositive, and its improvement in the objective value is
\[
-\Delta_4=-\max(0,w_4).
\]

Hence,
\begin{align*}
\Delta_2-\Delta_4
&=\max\bigl(0,\Delta w_2^++\min(0,w_2)\bigr)-\max(0,w_4)\\
&=\sum_{j=2}^{4}p_j-\min(0,w_3)-w_4-\max(0,w_4)\\
&\ge p_2+p_3-\max(0,w_4)-\max(0,w_4).
\end{align*}
Since
\[
\Delta w_5=p_2+p_3-\max(0,w_4),
\]
we have
\[
\Delta_2-\Delta_4\ge -\Delta_4+\Delta w_5.
\]
Therefore,
\[
\Delta w\ge \Delta_2-\Delta_4\ge -\Delta_4+\Delta w_5.
\]

Consequently, if
\[
-\Delta_4+\Delta w_5>0,
\]
then necessarily $\Delta w>0$, and the backward move cannot reduce the total waiting time. The contrapositive statement is therefore: if a backward move satisfies $\Delta w\le 0$, then it must hold that
\[
-\Delta_4+\Delta w_5\le 0.
\]
Hence, this inequality is a necessary condition for a backward move to improve the total waiting time.

\medskip
\noindent
\textbf{Step 2: Construction of the discriminant function and boundary analysis.}

Define
\[
H(p_4,w_4):=-\Delta_4+\Delta w_5.
\]
By Step 1, $H(p_4,w_4)\le 0$ is a necessary condition for the backward move to be potentially improving. We next examine the structure of $H$ according to the value of $w_4$.

\medskip
\noindent
\textbf{Case 1:} $w_4\le 0$.

By \cref{Theorem2},
\begin{align*}
w_4^2 &= \min(0,w_3)+\min(0,w_2)-p_3-p_2+w_4,\\
\Delta w_5 &= I_4=\max\bigl\{-p_4,\ \min(0,w_2)+\min(0,w_3),\ p_2+p_3-\max(0,w_4)\bigr\}.
\end{align*}
When $w_4\le 0$, we have $\max(0,w_4)=0$, so the third term becomes $p_2+p_3\ge 0$, whereas the first two terms are both no greater than $0$. Thus,
\[
\Delta w_5=p_2+p_3.
\]
Moreover, since $w_4\le 0$, it follows directly that $w_4^2\le 0$, and hence
\[
-\Delta_4=0.
\]
Therefore,
\[
H(p_4,w_4)=p_2+p_3\ge 0.
\]
Thus, when $w_4\le 0$, this backward move can never improve the total waiting time.

\medskip
\noindent
\textbf{Case 2:} $w_4>0$.

In this case,
\begin{align*}
w_4^2 &= \min(0,w_3)+\min(0,w_2)-p_3-p_2+w_4,\\
\Delta w_5 &= I_4=\max\bigl\{-p_4,\ \min(0,w_2)+\min(0,w_3),\ p_2+p_3-w_4\bigr\}.
\end{align*}

We further divide the analysis into two subcases.

\medskip
\noindent
\textbf{Subcase 2.1:} $w_4^2>0$.

In this case,
\[
\min(0,w_2)+\min(0,w_3)-p_2-p_3+w_4>0,
\]
and therefore the change in the objective contribution of job $4$ is
\[
-\Delta_4=\min(0,w_2)+\min(0,w_3)-p_2-p_3.
\]
Moreover, since
\[
\min(0,w_2)+\min(0,w_3)>p_2+p_3-w_4,
\]
we have
\[
\Delta w_5=\max\bigl\{-p_4,\ \min(0,w_2)+\min(0,w_3)\bigr\}.
\]

If
\[
\min(0,w_2)+\min(0,w_3)>-p_4,
\]
then
\[
H(p_4,w_4)=2\min(0,w_2)+2\min(0,w_3)-p_2-p_3\le 0.
\]
If
\[
\min(0,w_2)+\min(0,w_3)\le -p_4,
\]
then
\[
H(p_4,w_4)=\min(0,w_2)+\min(0,w_3)-p_2-p_3-p_4\le 0.
\]

\medskip
\noindent
\textbf{Subcase 2.2:} $w_4^2\le 0$.

In this case,
\[
\min(0,w_2)+\min(0,w_3)\le p_2+p_3-w_4,
\]
and thus
\[
-\Delta_4=-w_4,\qquad
\Delta w_5=\max\{-p_4,\ p_2+p_3-w_4\}.
\]

If
\[
p_2+p_3-w_4>-p_4,
\]
then
\[
H(p_4,w_4)=p_2+p_3-2w_4.
\]
Therefore, when
\[
w_4\ge \frac{p_2+p_3}{2},
\]
it holds that $H(p_4,w_4)\le 0$.

If
\[
p_2+p_3-w_4\le -p_4,
\]
then
\[
H(p_4,w_4)=-w_4-p_4\le 0.
\]

In summary, in the case $w_4>0$, the condition $H(p_4,w_4)\le 0$ characterizes the possibility that the backward move may still improve the total waiting time.

\medskip
\noindent
\textbf{Step 3: Extension from the local structure to a general backward candidate.}

Extending the local pair $(2,3)$ in the above analysis to a general interval $(j,\dots,i-1)$ yields the discriminant condition for a general backward move.

Furthermore, the above piecewise representation can be unified as
\begin{align*}
H(p_4,w_4)
=&\max\bigl(-w_4,\ \min(0,w_2)+\min(0,w_3),\ -p_2-p_3\bigr)\\
&+\max\bigl(-p_4,\ p_2+p_3-w_4,\ \min(0,w_2)+\min(0,w_3)\bigr).
\end{align*}
Hence, $H$ is the sum of maxima of finitely many affine functions, and is therefore a piecewise-linear convex function. Its minimum is attained at the boundaries where the dominant affine pieces switch. Accordingly, for any position $j\in[1,i-1]$, if
\begin{align*}
w_{\sigma_t(i)}
&=
\sum_{k=j}^{i-1}
\bigl(
p_{\sigma_t(k)}-\min(0,w_{\sigma_t(k)})
\bigr),\\
 p_{\sigma_t(i)}
&=
\sum_{k=j}^{i-1}
\bigl(
-\min(0,w_{\sigma_t(k)})
\bigr),
\end{align*}
then, by the above analysis, equality corresponds exactly to the critical boundary at which the backward move reaches the threshold and the function $H$ attains its minimum. Moreover, only when the corresponding cumulative quantities are below this boundary can the waiting time of the moved job remain nonnegative after the move. If the waiting time of the moved job remains negative after the move, then any further backward move cannot improve the total waiting time. Therefore, for a backward candidate to be meaningful in practice, one must further require
\begin{align*}
w_{\sigma_t(i)}
\ge
\sum_{k=j}^{i-1}
\bigl(
p_{\sigma_t(k)}-\min(0,w_{\sigma_t(k)})
\bigr),\;\land\;
p_{\sigma_t(i)}
\ge
\sum_{k=j}^{i-1}
\bigl(
-\min(0,w_{\sigma_t(k)})
\bigr).
\end{align*}
Any position $j$ satisfying the above inequalities is called a backward candidate position of job $\sigma_t(i)$, denoted by $\theta_{i,j}^b$.

This condition indicates that only when job $\sigma_t(i)$ satisfies the above inequalities can inserting $\sigma_t(i)$ backward into position $j$ potentially reduce the objective value. Therefore, the above two inequalities provide the necessary discriminant conditions for the set of backward candidate solutions.
\end{proof}

\subsection{ Proof of \cref{Lemma3}}\label{app:Lemma3}
\begin{proof}
We first show that only a job located after position $i$ can repair the queue discontinuity at position $i$. Taking position $i$ as the dividing point, we partition the prefix of the sequence into
\[
P_a=\{\sigma_t(1),\dots,\sigma_t(i-1)\}.
\]

Viewing $P_a$ as an aggregate block, its completion time $C_{P_a}$ can be written as
\[
C_{P_a}
=
C_{P_o}
+
\sum_{j\in P_a} p_j
+
\sum_{j\in P_a}\bigl(-\min(0,w_j)\bigr),
\]
where $C_{P_o}$ is the completion time immediately before $P_a$, $\sum_{j\in P_a} p_j$ is the total processing time of the jobs in $P_a$, and $\sum_{j\in P_a}\bigl(-\min(0,w_j)\bigr)$ is the total machine idle time accumulated within $P_a$. Since the set of jobs in $P_a$ is fixed, the total processing time is constant. It follows that any change in $C_{P_a}$ is entirely determined by the total machine idle time within $P_a$.

To repair the queue discontinuity at job $\sigma_t(i)$, one must increase $C_{P_a}$, that is, delay the completion time of the entire block $P_a$. Since the total processing time of $P_a$ remains unchanged, this is equivalent to shifting an equal amount of machine idle time from after position $i$ into $P_a$, thereby postponing the service start time of job $\sigma_t(i)$. Under such an adjustment, job $\sigma_t(i)$ and all subsequent jobs are affected by an increasing flow. However, for an increasing flow, a queue discontinuity does not obstruct the realization of the ideal move, and therefore no repair is needed there.

Therefore, the region that truly requires repair is the one affected by the decreasing flow. This implies that the repair must be achieved by moving some job located after position $i$ forward, so as to repair the queue discontinuity at position $i$.

Next, we prove that, among the jobs located after position $i$, only those whose release times are strictly smaller than that of job $\sigma_t(i)$ can repair the discontinuity.

At this point, we have $w_{\sigma_t(i)}\le 0$. Since \cref{Proposition1} only determines the release-time order of two adjacent jobs, we first move job $\sigma_t(j)$ to position $i+1$, and denote its waiting time there by $w_{\sigma_t(j)}'$. By \cref{Theorem2},
\begin{align*}
w_{\sigma_t(j)}'
=
w_{\sigma_t(j)}
-
\sum_{k=i+1}^{j-1}
\bigl(
p_{\sigma_t(k)}-\min(0,w_{\sigma_t(k)})
\bigr).
\end{align*}

Then move job $\sigma_t(j)$ further to position $i$, and denote its waiting time by $w_{\sigma_t(j)}^{2}$. Again by \cref{Theorem2},
\begin{align*}
w_{\sigma_t(j)}^{2}
=
w_{\sigma_t(j)}
-
\sum_{k=i}^{j-1}
\bigl(
p_{\sigma_t(k)}-\min(0,w_{\sigma_t(k)})
\bigr).
\end{align*}

Under the standing assumption $w_{\sigma_t(i)}\le 0$, we have
\[
w_{\sigma_t(i)}=\min(0,w_{\sigma_t(i)}).
\]

If the queue discontinuity at job $\sigma_t(i)$ can be repaired, then it is necessary that
\[
w_{\sigma_t(j)}^{2}>w_{\sigma_t(i)}.
\]
Substituting the above expression gives
\[
w_{\sigma_t(j)}'-p_{\sigma_t(i)}+w_{\sigma_t(i)}>w_{\sigma_t(i)},
\]
that is,
\[
w_{\sigma_t(j)}'>p_{\sigma_t(i)}.
\]
Since $p_{\sigma_t(i)}\ge 0$, it follows that
\[
w_{\sigma_t(j)}'>0.
\]

Therefore, the criterion in \cref{Proposition1} is satisfied:
\begin{align*}
w_{\sigma_t(j)}'
=
\max(0,w_{\sigma_t(j)}')
>
\max(0,w_{\sigma_t(i)})+p_{\sigma_t(i)}
=
p_{\sigma_t(i)},
\end{align*}
where the last equality holds because $w_{\sigma_t(i)}\le 0$. Hence,
\[
r_{\sigma_t(j)}<r_{\sigma_t(i)}.
\]

We thus conclude that the queue discontinuity at position $i$ can be repaired only if job $\sigma_t(j)$ satisfies
\[
r_{\sigma_t(j)}<r_{\sigma_t(i)}
\qquad\text{and}\qquad
w_{\sigma_t(j)}'>0.
\]
Equivalently, only when the release time of job $\sigma_t(j)$ is strictly smaller than that of job $\sigma_t(i)$, and its waiting time after being moved to position $i+1$ is positive, can the queue discontinuity at job $\sigma_t(i)$ be repaired.
\end{proof}

\subsection{Proof of the Covering Relation for the Repair Set $B$}
\label{app:cover-B}

\begin{lemma}
\label{LemmaCoverB}
Let $\sigma^t$ be the current schedule, and suppose that the discontinuity at $b$ corresponds to a forward-type local optimum structure. Let $\sigma^t(j)$ denote the single-job compensation selected first by the greedy rule at discontinuity $b$, and let
\[
B=\{\sigma^t(h_1),\dots,\sigma^t(h_m)\}
\]
be a better combined repair set for the same discontinuity $b$. Then the current positions of all jobs in $B$ are contained in the traversal range of the forward candidate set induced by job $\sigma^t(j)$.
\end{lemma}

\begin{proof}
First, by \cref{Lemma3}, any job that can repair discontinuity $b$ must have a release time strictly smaller than that of the job at the discontinuity. Therefore, the jobs capable of repairing discontinuity $b$ are ordered nondecreasingly by release times around the discontinuity.

Next, consider the propagation of the decrease flow induced when a compensation job moves forward. By the recursive formula for the decrease flow of waiting times, the initial decrease-flow value of job $\sigma^t(j)$ is
\[
\Delta^-_{\sigma^t(j+1)}
=
p_{\sigma^t(j)}-\min\bigl(0,w_{\sigma^t(j)}\bigr).
\]
The initial value of this theoretical propagation upper bound can be uniformly tightened to the processing time $p_{\sigma^t(j)}$. Specifically, if $w_{\sigma^t(j)}>0$, then
\[
\Delta^-_{\sigma^t(j+1)}=p_{\sigma^t(j)}.
\]
If $w_{\sigma^t(j)}\le 0$, although
\[
\Delta^-_{\sigma^t(j+1)}
=
p_{\sigma^t(j)}-w_{\sigma^t(j)}
\ge p_{\sigma^t(j)},
\]
it follows from the recursive decrease-flow formula together with the relation in \cref{Proposition1},
\[
\max\bigl(0,w_{\sigma^t(j)}\bigr)+p_{\sigma^t(j)}
\ge
\max\bigl(0,w_{\sigma^t(j+1)}\bigr),
\]
that once this initial decrease flow enters the subsequent continuous queue, it is immediately truncated by the waiting times of the subsequent jobs. Hence, the theoretical upper bound on its further propagation may still be uniformly regarded as $p_{\sigma^t(j)}$. Therefore, when constructing the forward candidate set, the initial upper bound of the compensation job's propagation can be tightened to its processing time.

Accordingly, under the nondecreasing order of release times, for any given continuous queue interval $C$, a necessary condition for the forward candidate set of job $\sigma^t(j)$ can be written as
\[
\sum_{l\in C}\bigl(p_{\sigma^t(l)}-p_{\sigma^t(j)}\bigr)\le 0.
\]
Similarly, for any job $\sigma^t(h_q)\in B$, we have
\[
\sum_{l\in C}\bigl(p_{\sigma^t(l)}-p_{\sigma^t(h_q)}\bigr)\le 0.
\]
It should be emphasized that these conditions correspond to the forward candidate sets derived from the upper bound of the decrease flow, and are used to characterize the forward range that may be covered by an ideal improvement direction. The actually reachable range is further truncated by the true waiting-time structure and is therefore usually smaller than this theoretical range.

We now compare the covering relation between job $\sigma^t(j)$ and job $\sigma^t(h_q)$ under the same discontinuity $b$.

\textbf{Case 1:} $p_{\sigma^t(h_q)}<p_{\sigma^t(j)}$.  
For the same region $C$, the necessary condition for a theoretical forward candidate can be written uniformly as
\[
\sum_{l\in C}\bigl(p_{\sigma^t(l)}-p\bigr)\le 0,
\]
where $p$ denotes the processing time of the job serving as the local compensation benchmark. This condition is equivalent to
\[
\sum_{l\in C}p_{\sigma^t(l)}\le |C|\,p.
\]
Therefore, this necessary condition is monotone nondecreasing in the benchmark processing time $p$: if $p_1\le p_2$, then the set of theoretical candidate endpoints corresponding to $p_1$ must be contained in that corresponding to $p_2$. Hence, under the same discontinuity and the same family of intervals, a compensation job with a larger processing time has a theoretical forward candidate coverage range no smaller than that of a job with a smaller processing time. Therefore, in this case, the theoretical forward candidate traversal range of $\sigma^t(j)$ covers the current position of job $\sigma^t(h_q)$.

\textbf{Case (2):} $p_{\sigma^t(h_q)}\ge p_{\sigma^t(j)}.$
In this case, under the tightened upper bound on the initial decrease flow, at the same position, the decrease flow induced by job $\sigma^t(j)$ is no greater than that induced by job $\sigma^t(h_q)$. However, by the greedy rule, job $\sigma^t(j)$ is selected because its immediate compensation at discontinuity $b$ is no smaller than that of any candidate job. Hence, under the same discontinuity, if jobs $\sigma^t(j)$ and $\sigma^t(h_q)$ are respectively moved to discontinuity $b$ and then moved forward back to their original positions, the total accumulated idle time induced by the former process is no smaller than that induced by the latter.

By the monotone decay property of the recursive formula for the decrease flow (\cref{Lemma1}), if job $\sigma^t(j)$, which has weaker initial propagation capability, is to attain a total accumulated idle time no smaller than that of job $\sigma^t(h_q)$, then it can do so only by relying on a propagation length no shorter than that of $\sigma^t(h_q)$. On the other hand, the necessary condition defining the forward candidate set is itself based on accumulated-sum tests over consecutive queue intervals, and hence its candidate endpoints advance continuously along the queue in the forward direction by definition. That is, once the accumulated idle time associated with a job is sufficient to cover a longer consecutive interval, then all positions corresponding to shorter consecutive intervals preceding it must also fall within its forward candidate traversal range.

Therefore, in the present case, as long as the propagation length of the decrease flow induced by job $\sigma^t(j)$ is no shorter than that of job $\sigma^t(h_q)$, the consecutive interval corresponding to the current position of job $\sigma^t(h_q)$ must be contained in the family of forward candidate intervals of job $\sigma^t(j)$. Hence, the forward candidate traversal range of job $\sigma^t(j)$ also covers the current position of job $\sigma^t(h_q)$.

Combining the two cases, for any $\sigma^t(h_q)\in B$, its current position is contained in the traversal range of the forward candidate set induced by job $\sigma^t(j)$. This completes the proof.
\end{proof}

\section{Omitted Proofs in Section 4}
\subsection{ Proof of \cref{Lemma4}}\label{app:Lemma4}
\begin{proof}
Consider an arbitrary sequence $\sigma_t$ and any two adjacent jobs $\sigma_t(i)$ and $\sigma_t(i+1)$. For ease of exposition, and without loss of generality, let $\sigma_t(i)=4$ and $\sigma_t(i+1)=5$. This relabeling is used only to denote the adjacent pair and has no effect on the subsequent derivation.

Swapping jobs $4$ and $5$ is equivalent to moving job $4$ to the position immediately after job $5$. After the swap, the waiting times of the two jobs are given by
\begin{align*}
    w_4^2&=\max(0,w_4)+p_4+p_5-\min(0,w_4)-\Delta w_6^{-},\\
    w_5^2&=w_5-\Delta w_5^{-}
    =w_5-\bigl(p_4-\min(0,w_4)\bigr).
\end{align*}

Let $\Delta w_{\mathrm{pos}}$ denote the change in the sum of the positive waiting times of jobs $4$ and $5$ before and after the swap. Then
\begin{align*}
\Delta w_{\mathrm{pos}}
=\max(0,w_4^2)+\max(0,w_5^2)-\max(0,w_4)-\max(0,w_5).
\end{align*}
If the swap is to reduce the total waiting time of jobs $4$ and $5$, it is necessary that
\[
\Delta w_{\mathrm{pos}}<0.
\]

Let $\Delta w_{\mathrm{neg}}$ denote the change in machine idle time before and after the swap, namely,
\begin{align*}
\Delta w_{\mathrm{neg}}
=\min(0,w_4^2)+\min(0,w_5^2)-\min(0,w_4)-\min(0,w_5).
\end{align*}
If $\Delta w_{\mathrm{neg}}<0$, then the waiting times of the subsequent jobs increase; if $\Delta w_{\mathrm{neg}}=0$, then the subsequent jobs are unaffected; and if $\Delta w_{\mathrm{neg}}>0$, then the waiting times of the subsequent jobs decrease.

Now suppose that jobs $4$ and $5$ satisfy the first-come-first-served condition in Proposition 1, namely,
\begin{align*}
    \max(0,w_4)+p_4\ge \max(0,w_5).
\end{align*}

\medskip
\noindent
\textbf{Case 1:} $w_4\le 0$.

In this case, from $\max(0,w_4)+p_4\ge \max(0,w_5)$ it follows that $w_5\le p_4$.

\medskip
\noindent
\textbf{Subcase 1.1:} $0<w_5\le p_4$.

The new waiting times of jobs $4$ and $5$ after the swap are
\begin{align*}
    w_4^2
    &=\max(0,w_4)+p_4+p_5-\min(0,w_4)-\Delta w_6^{-}\\
    &=p_4+p_5-w_5,\\
    w_5^2
    &=w_5-\bigl(p_4-\min(0,w_4)\bigr)
    =w_5-p_4+w_4.
\end{align*}
Since $0<w_5\le p_4$, we have
\[
w_4^2=p_4+p_5-w_5\ge 0,
\qquad
w_5^2=w_5-p_4+w_4\le 0.
\]

The change in machine idle time is
\begin{align*}
\Delta w_{\mathrm{neg}}
&=\min(0,w_4^2)+\min(0,w_5^2)-\min(0,w_4)-\min(0,w_5)\\
&=w_5-p_4.
\end{align*}
Because $w_5\le p_4$, it follows that $\Delta w_{\mathrm{neg}}=w_5-p_4\le 0$. Hence, if $w_5<p_4$, the waiting times of the subsequent jobs increase; if $w_5=p_4$, they remain unchanged.

The change in positive waiting time is
\begin{align*}
\Delta w_{\mathrm{pos}}
&=\max(0,w_4^2)+\max(0,w_5^2)-\max(0,w_4)-\max(0,w_5)\\
&=p_4+p_5-w_5-w_5\\
&=p_4+p_5-2w_5.
\end{align*}
Therefore,
\[
\Delta w_{\mathrm{pos}}<0
\iff
p_4+p_5-2w_5<0
\iff
w_5>\frac{1}{2}(p_4+p_5).
\]
Thus, $w_5>\frac{1}{2}(p_4+p_5)$ is the necessary and sufficient condition for $\Delta w_{\mathrm{pos}}<0$ in this subcase.

Accordingly, under $0<w_5\le p_4$, the condition $\Delta w_{\mathrm{pos}}<0$ is equivalent to
\[
w_5>\frac{1}{2}(p_4+p_5).
\]
Moreover, combining this with $w_5\le p_4$ yields $p_4>p_5$. Hence, in this subcase, the condition $p_4>p_5$ is automatically implied and need not be imposed separately.

\medskip
\noindent
\textbf{Subcase 1.2:} $w_5\le 0$.

The new waiting times are
\begin{align*}
    w_4^2
    &=\max(0,w_4)+p_4+p_5-\min(0,w_4)-\Delta w_6^{-}\\
    &=p_4+p_5-w_5,\\
    w_5^2
    &=w_5-\bigl(p_4-\min(0,w_4)\bigr)
    =w_5-p_4+w_4.
\end{align*}
Since $w_5\le 0$, we have
\[
w_4^2=p_4+p_5-w_5\ge 0,
\qquad
w_5^2=w_5-p_4+w_4\le 0.
\]

The change in machine idle time is
\begin{align*}
\Delta w_{\mathrm{neg}}
&=\min(0,w_4^2)+\min(0,w_5^2)-\min(0,w_4)-\min(0,w_5)\\
&=w_5-p_4.
\end{align*}
Since $w_5\le 0\le p_4$, it follows that $\Delta w_{\mathrm{neg}}\le 0$. Hence, if $w_5<p_4$, the waiting times of the subsequent jobs increase; if $w_5=p_4=0$, they remain unchanged.

The change in positive waiting time is
\begin{align*}
\Delta w_{\mathrm{pos}}
&=\max(0,w_4^2)+\max(0,w_5^2)-\max(0,w_4)-\max(0,w_5)\\
&=p_4+p_5-w_5.
\end{align*}
Since $w_5\le 0$, we obtain $\Delta w_{\mathrm{pos}}\ge 0$. Therefore, in this subcase,
\[
\Delta w_{\mathrm{pos}}\ge 0.
\]

\medskip
\noindent
\textbf{Case 2:} $w_4>0$.

In this case, from $\max(0,w_4)+p_4\ge \max(0,w_5)$ it follows that
\[
w_5\le p_4+w_4.
\]

\medskip
\noindent
\textbf{Subcase 2.1:} $w_5\le 0$.

The new waiting times are
\begin{align*}
    w_4^2
    &=\max(0,w_4)+p_4+p_5-\min(0,w_4)-\Delta w_6^{-}\\
    &=w_4+p_4+p_5-w_5,\\
    w_5^2
    &=w_5-\bigl(p_4-\min(0,w_4)\bigr)
    =w_5-p_4.
\end{align*}
Since $w_5\le 0$, we have
\[
w_4^2=w_4+p_4+p_5-w_5\ge 0,
\qquad
w_5^2=w_5-p_4\le 0.
\]

The change in machine idle time is
\begin{align*}
\Delta w_{\mathrm{neg}}
&=\min(0,w_4^2)+\min(0,w_5^2)-\min(0,w_4)-\min(0,w_5)\\
&=-p_4.
\end{align*}
Since $p_4\ge 0$, we have $\Delta w_{\mathrm{neg}}=-p_4\le 0$. Hence, if $p_4>0$, the waiting times of the subsequent jobs increase; if $p_4=0$, they remain unchanged.

The change in positive waiting time is
\begin{align*}
\Delta w_{\mathrm{pos}}
&=\max(0,w_4^2)+\max(0,w_5^2)-\max(0,w_4)-\max(0,w_5)\\
&=p_4+p_5-w_5.
\end{align*}
Since $w_5\le 0$, it follows that $\Delta w_{\mathrm{pos}}\ge 0$. Therefore, in this subcase,
\[
\Delta w_{\mathrm{pos}}\ge 0.
\]

\medskip
\noindent
\textbf{Subcase 2.2:} $0<w_5\le p_4$.

The new waiting times are
\begin{align*}
    w_4^2
    &=\max(0,w_4)+p_4+p_5-\min(0,w_4)-\Delta w_6^{-}\\
    &=w_4+p_4+p_5-w_5,\\
    w_5^2
    &=w_5-\bigl(p_4-\min(0,w_4)\bigr)
    =w_5-p_4.
\end{align*}
Since $0<w_5\le p_4$, we have
\[
w_4^2=w_4+p_4+p_5-w_5\ge 0,
\qquad
w_5^2=w_5-p_4\le 0.
\]

The change in machine idle time is
\begin{align*}
\Delta w_{\mathrm{neg}}
&=\min(0,w_4^2)+\min(0,w_5^2)-\min(0,w_4)-\min(0,w_5)\\
&=w_5-p_4.
\end{align*}
Because $w_5\le p_4$, it follows that $\Delta w_{\mathrm{neg}}\le 0$. Hence, if $w_5<p_4$, the waiting times of the subsequent jobs increase; if $w_5=p_4$, they remain unchanged.

The change in positive waiting time is
\begin{align*}
\Delta w_{\mathrm{pos}}
&=\max(0,w_4^2)+\max(0,w_5^2)-\max(0,w_4)-\max(0,w_5)\\
&=p_4+p_5+w_4-w_5-w_4-w_5\\
&=p_4+p_5-2w_5.
\end{align*}
Therefore,
\[
\Delta w_{\mathrm{pos}}<0
\iff
p_4+p_5-2w_5<0
\iff
w_5>\frac{1}{2}(p_4+p_5).
\]
Thus, under $0<w_5\le p_4$, the condition $\Delta w_{\mathrm{pos}}<0$ is again equivalent to
\[
w_5>\frac{1}{2}(p_4+p_5).
\]
As above, combining this inequality with $w_5\le p_4$ implies $p_4>p_5$, so the latter condition is automatically satisfied in this subcase.

\medskip
\noindent
\textbf{Subcase 2.3:} $p_4<w_5\le p_4+w_4$.

The new waiting times are
\begin{align*}
    w_4^2
    &=\max(0,w_4)+p_4+p_5-\min(0,w_4)-\Delta w_6^{-}\\
    &=w_4+p_4+p_5-p_4\\
    &=w_4+p_5,\\
    w_5^2
    &=w_5-\bigl(p_4-\min(0,w_4)\bigr)
    =w_5-p_4.
\end{align*}
Since $p_4<w_5$ and $w_4>0$, we have
\[
w_4^2=w_4+p_5\ge 0,
\qquad
w_5^2=w_5-p_4\ge 0.
\]

The change in machine idle time is
\begin{align*}
\Delta w_{\mathrm{neg}}
&=\min(0,w_4^2)+\min(0,w_5^2)-\min(0,w_4)-\min(0,w_5)\\
&=0.
\end{align*}
Hence, the waiting times of the subsequent jobs remain unchanged.

The change in positive waiting time is
\begin{align*}
\Delta w_{\mathrm{pos}}
&=\max(0,w_4^2)+\max(0,w_5^2)-\max(0,w_4)-\max(0,w_5)\\
&=(w_4+p_5)+(w_5-p_4)-w_4-w_5\\
&=p_5-p_4.
\end{align*}
Therefore,
\[
\Delta w_{\mathrm{pos}}<0
\iff
p_5-p_4<0
\iff
p_4>p_5.
\]
Thus, in this subcase, $p_4>p_5$ is the necessary and sufficient condition for $\Delta w_{\mathrm{pos}}<0$.

Moreover, under $p_4<w_5\le p_4+w_4$, if $p_4>p_5$, then automatically
\[
w_5>p_4>\frac{1}{2}(p_4+p_5),
\]
so the condition $w_5>\frac{1}{2}(p_4+p_5)$ is automatically satisfied and need not be imposed separately.

\medskip
\noindent
Combining the above cases, we see that swapping jobs $4$ and $5$ can only leave the waiting times of the subsequent jobs unchanged or increase them. Therefore, in order for the swap to improve the total waiting time, it is necessary that
\[
\Delta w_{\mathrm{pos}}<0.
\]
Collecting all cases in which $\Delta w_{\mathrm{pos}}<0$, we obtain
\[
\Delta w_{\mathrm{pos}}<0
\iff
\left(
w_5>\frac{1}{2}(p_4+p_5)
\right)
\land
\left(
p_4>p_5
\right).
\]

That is, when jobs are arranged in nondecreasing order of release times, the necessary and sufficient conditions for the adjacent swap to reduce the total waiting time are
\[
w_5>\frac{1}{2}(p_4+p_5)
\qquad\text{and}\qquad
p_4>p_5.
\]

Finally, under the condition $\Delta w_{\mathrm{pos}}<0$, we examine when $\Delta w_{\mathrm{neg}}=0$. There are three possibilities:

\begin{enumerate}
    \item If $w_4\le 0$ and $w_5=p_4$, then the waiting times of the subsequent jobs remain unchanged, and
    \[
    w_5^2=w_5-p_4+\min(0,w_4)=w_4.
    \]
    
    \item If $w_4>0$ and $w_5=p_4$, then the waiting times of the subsequent jobs remain unchanged, and
    \[
    w_5^2=w_5-p_4=0.
    \]
    
    \item If $w_4>0$ and $w_5>p_4$, then the waiting times of the subsequent jobs remain unchanged, and
    \[
    w_5^2=w_5-p_4>0.
    \]
\end{enumerate}

Hence, the above three cases can be unified as
\[
w_5^2\ge \min(0,w_4).
\]

\end{proof}

\subsection{From Over-Consumption to Under-Consumption and Their Joint Elimination}
\label{app:over under}

\begin{lemma}[Joint elimination of over-consumption and under-consumption local optima]
\label{over under}
Let the current schedule $\sigma^t$ be a backward-type local optimum, and suppose that it contains an over-consumption local optimum whose corresponding excessively forward-moved job is $\sigma^t(j)$. Let
\[
B=\{\sigma^t(h_1),\dots,\sigma^t(h_m)\}
\]
denote a set of jobs whose combined forward adjustment replaces the individual forward move of $\sigma^t(j)$. Then, after Algorithm~\ref{algorithm1} first exposes the relevant local structure explicitly, and Algorithm~\ref{algorithm3} subsequently moves $\sigma^t(j)$ backward to a position where the discontinuity induced by it can be repaired, at most the following two cases may occur:
\begin{enumerate}
    \item no under-consumption local optimum is further induced, in which case the over-consumption local optimum can be eliminated through finitely many backward candidate traversals combined with the adjacent-swap rule;
    \item an under-consumption local optimum is further induced, in which case the newly exposed local optimum can still be eliminated within the combined neighborhood enumerated by Algorithm~\ref{algorithm1} and Algorithm~\ref{algorithm3}.
\end{enumerate}
Therefore, regardless of which case occurs, all subsequent local structures generated by this over-consumption local optimum can ultimately be eliminated successively within the combined neighborhood, thereby further reducing the side effects caused by the increase flow.
\end{lemma}

\begin{proof}
Let $\sigma^t(j)$ be the excessively forward-moved job corresponding to the over-consumption local optimum. The essence of this structure is that, under the greedy rule, although the individual forward move of $\sigma^t(j)$ is accepted first at the local level, this move introduces too much machine idle time in the increase-flow region and is therefore not globally optimal. Correspondingly, there exists a set of jobs
\[
B=\{\sigma^t(h_1),\dots,\sigma^t(h_m)\},
\]
whose combined forward adjustment can replace the individual forward move of $\sigma^t(j)$, induce less total machine idle time, and hence yield a smaller additional loss caused by the increase flow, resulting in a better overall change in total waiting time.

To eliminate this over-consumption local optimum, Algorithm~\ref{algorithm1} first makes the relevant local structure explicit. Then Algorithm~\ref{algorithm3} moves $\sigma^t(j)$ backward, through backward candidate traversal, to a position where the discontinuity induced by it can be reabsorbed. We then distinguish two cases for the increase-flow region induced by this backward move.

\textbf{Case 1.} No under-consumption local optimum is further exposed in the increase-flow region.

In this case, the backward repair of $\sigma^t(j)$ has already removed the main structural imbalance originally caused by its excessive forward move, and the remaining effect in the increase-flow region consists only of additional losses that can be further compressed by local forward adjustments. One may then further realize local forward adjustments of some jobs in the set $B$ through the adjacent-swap rule, thereby replacing the original single-job structure by a combined compensation structure. Since the total machine idle time induced by this combined structure is strictly smaller than that induced by the individual forward move of $\sigma^t(j)$, the total waiting time is strictly reduced. Hence, this over-consumption local optimum is directly eliminated within the combined neighborhood.

\textbf{Case 2.} An under-consumption local optimum is further induced in the increase-flow region.

In this case, there exist some jobs, two representative ones being denoted by $\sigma^t(h_a)$ and $\sigma^t(h_b)$, such that under the current local structure a single adjacent swap is insufficient to strictly reduce the total waiting time corresponding to the increase flow, whereas a better improvement is achieved by moving $\sigma^t(h_a)$ forward across several positions beyond $\sigma^t(h_b)$. In other words, after Algorithm~\ref{algorithm3} moves $\sigma^t(j)$ backward to a suitable position, a new under-consumption local optimum becomes explicitly exposed in the increase-flow region induced by this move.

The essence of an under-consumption local optimum is precisely that an ideal forward direction admitting strict improvement has not yet been fully released. To handle this structure, Algorithm~\ref{algorithm1} again first applies a forward move together with the bottleneck breakthrough rule to make the corresponding discontinuity explicit, and Algorithm~\ref{algorithm3} then further coordinates the local state through backward candidate traversal and adjacent swaps.

Moreover, since the decrease-flow regions corresponding to $\sigma^t(j)$ and $\sigma^t(h_a)$ overlap, \cref{Lemma1} implies that when the decrease flows induced by multiple job moves act on the same region, they must compete with one another because the room for further compressing positive waiting times in that region is limited. Therefore, once $\sigma^t(j)$ and $\sigma^t(h_a)$ cannot coexist, they cannot simultaneously serve as the final dominant improvement structure in the overlapping region; one of them must be replaced by the more effective one. Hence, this newly exposed under-consumption local optimum can likewise be further eliminated within the combined neighborhood.

Combining the above arguments, when handling an over-consumption local optimum, Algorithm~\ref{algorithm1} first explicitly exposes the relevant structure, and Algorithm~\ref{algorithm3} then performs the backward repair. If no under-consumption local optimum is further induced, the structure can be eliminated directly through combined forward adjustments; if an under-consumption local optimum is further induced, the newly exposed structure can still be eliminated successively by following the same order, namely, Algorithm~\ref{algorithm1} first exposes the structure and Algorithm~\ref{algorithm3} then repairs it. Therefore, regardless of which case occurs, all subsequent local structures generated by the over-consumption local optimum can be successively eliminated within the combined neighborhood, thereby minimizing the side effects caused by the increase flow as much as possible.
\end{proof}

\bibliography{reference}

\end{document}